\documentclass[11pt]{article}
\usepackage{amssymb}
\usepackage{amsmath}
\usepackage[all]{xy}
\usepackage{times,esint}
\usepackage{graphicx}

\title{Universal and complete sets in martingale theory\indent}
\author{Dominique LECOMTE and Miroslav ZELEN\'Y$^1$}
\date{\today}

\def\ufootnote#1{\let\savedthfn\thefootnote\let\thefootnote\relax
\footnote{#1}\let\thefootnote\savedthfn\addtocounter{footnote}{-1}}

\newcommand{\Ana}{{\it\Sigma}^{1}_{1}}

\newcommand{\Ca}{{\it\Pi}^{1}_{1}}
\newcommand{\Boraone}{{\it\Sigma}^{0}_{1}}
\newcommand{\ca}{{\bf\Pi}^{1}_{1}}
\newcommand{\Borel}{{\it\Delta}^{1}_{1}}

\newcommand{\Borone}{{\it\Delta}^{0}_{1}}

\newcommand{\boratwo}{{\bf\Sigma}^{0}_{2}}
\newcommand{\borathree}{{\bf\Sigma}^{0}_{3}}

\newcommand{\Bormone}{{\it\Pi}^{0}_{1}}
\newcommand{\Bormtwo}{{\it\Pi}^{0}_{2}}
\newcommand{\bormtwo}{{\bf\Pi}^{0}_{2}}
\newcommand{\bormthree}{{\bf\Pi}^{0}_{3}}

\newtheorem{thm} {Theorem} [section]
\newtheorem{cor} [thm] {Corollary}
\newtheorem{lem} [thm] {Lemma}
\newtheorem{prop} [thm] {Proposition}
\newtheorem{defi} [thm] {Definition}

\begin{document}

\maketitle

\centerline{$\bullet$ Universit\' e Paris 6, Institut de Math\'ematiques de Jussieu, Projet Analyse Fonctionnelle}

\centerline{Couloir 16-26, 4\`eme \'etage, Case 247, 4, place Jussieu, 75 252 Paris Cedex 05, France}

\centerline{dominique.lecomte@upmc.fr}\bigskip

\centerline{$\bullet$ Universit\'e de Picardie, I.U.T. de l'Oise, site de Creil,}

\centerline{13, all\'ee de la fa\"\i encerie, 60 107 Creil, France}\bigskip

\centerline{$\bullet^1$ Charles University, Faculty of Mathematics and Physics, Department of Mathematical Analysis}

\centerline{Sokolovsk\'a 83, 186 75 Prague, Czech Republic}

\centerline{zeleny@karlin.mff.cuni.cz}\bigskip\bigskip\bigskip\bigskip\bigskip\bigskip

\ufootnote{{\it 2010 Mathematics Subject Classification.}~03E15, 60G42}

\ufootnote{{\it Keywords and phrases.}~divergence, Lebesgue measure, martingale, measure zero, $G_{\delta\sigma}$ set}

\ufootnote{{\it Acknowledgements.}~The main result was obtained during the first author's
stay at Charles University in Prague in August 2014. The first author thanks Charles University in Prague for the hospitality. We also thank A. Louveau for asking about the converse of Doob's theorem and for suggesting that the set of everywhere converging martingales is $\ca$-complete.}

\noindent {\bf Abstract.} The Doob convergence theorem implies that the set of divergence of any martingale has measure zero. We prove that, conversely, any $G_{\delta\sigma}$ subset of the Cantor space with Lebesgue-measure zero can be represented as the set of divergence of some martingale. In fact, this is effective and uniform. A consequence of this is that the set of everywhere converging martingales is $\ca$-complete, in a uniform way. We derive from this some universal and complete sets for the whole projective hierarchy, via a general method. We provide some other complete sets for the classes $\ca$ and ${\bf\Sigma}^1_2$ in the theory of martingales.

\vfill\eject

\section{$\!\!\!\!\!\!$ Introduction}\indent

 The reader should see [K2] for the notation used in this paper.

\begin{defi} We say that a map $f\! :\! 2^{<\omega}\!\rightarrow\! [0,1]$ is a {\bf martingale} if $f(s)\! =\!\frac{f(s0)\! +\! f(s1)}{2}$ for each $s\!\in\! 2^{<\omega}$. The set of martingales is denoted by $\cal M$ and is a compact subset of $[0,1]^{2^{<\omega}}$ (equipped with the usual product topology).\end{defi}

 This terminology is not the standard one, but the set $\cal M$ can be interpreted as the set of all discrete martingales (in the classical sense) taking values in [0,1], as follows. If
 $s\!\in\! 2^{<\omega}$, then
$$N_s\! :=\!\{\beta\!\in\! 2^\omega\mid s\!\subseteq\!\beta\}$$
is the usual basic clopen set. Let $f\!\in\! {\cal M}$. If $n\!\in\!\omega$, then let ${\cal S}_n$ be the
$\sigma$-algebra on $2^\omega$ generated by $\{ N_s\mid s\!\in\! 2^n\}$, and
$f_n\! :\! 2^\omega\!\rightarrow\! [0,1]$ be defined by $f_n(\beta )\! :=\! f(\beta\vert n)$. Then the sequence $(f_n)_{n\in\omega}$ is a discrete martingale taking values in [0,1] with respect to the sequence of $\sigma$-algebras $({\cal S}_n)_{n\in\omega}$ and the usual Lebesgue product measure $\lambda$ on $2^\omega$. Conversely, if $(f_n)_{n\in\omega}$ is any such martingale, it can be viewed as an element of $\cal M$ by setting $f(s)\! :=\! f_{\vert s\vert}(\alpha )$ if
$\alpha\!\in\! N_s$. This definition is correct because $f_{\vert s\vert}$, as a function measurable with respect to ${\cal S}_{\vert s\vert}$, has a constant value on $N_s$.

\begin{defi} Let $f$ be a martingale and $\beta\!\in\! 2^\omega$. The {\bf oscillation} of $f$ at
$\beta$ is the number
$$\mbox{osc}(f,\beta )\! :=\!
\mbox{inf}_{N\in\omega}~\mbox{sup}_{p,q\geq N}~\vert f(\beta\vert p)\! -\! f(\beta\vert q)\vert .$$
The {\bf set of divergence} of $f$ is
$D(f)\! :=\!\{\beta\!\in\! 2^\omega\mid\mbox{osc}(f,\beta )\! >\! 0\}$.\end{defi}

 By definition, if $f$ is a martingale, then
$$\beta\!\in\! D(f)\Leftrightarrow\exists r\!\in\!\omega ~~\forall N\!\in\!\omega ~~
\exists p,q\!\geq\! N~~\vert f(\beta\vert p)\! -\! f(\beta\vert q)\vert\! >\! 2^{-r}.$$
This shows that $D(f)\!\in\!\borathree$. Moreover, $D(f)$ has $\lambda$-measure zero, by Doob's convergence theorem (see Chapter XI, Section 14 in [D]). So it is natural to ask whether any
$\borathree$ subset of $2^\omega$ with $\lambda$-measure zero is the set of divergence of some martingale (this question was asked by Louveau). We answer positively:

\begin{thm} \label{main} Let $B$ be a subset of $2^\omega$. Then the following are equivalent:\smallskip

\noindent (a) $B$ is $\borathree$ and has $\lambda$-measure zero,\smallskip

\noindent (b) there is a martingale $f$ with $B\! =\! D(f)$.\end{thm}

\begin{defi} Let $\bf\Gamma$ be a class of subsets of Polish spaces, $X,Y$ be Polish spaces, and
${\cal U}\!\subseteq\! Y\!\times\! X$.\smallskip

\noindent (a) We say that ${\cal U}$ is $Y$-{\bf universal for the} $\bf\Gamma$ {\bf subsets of} $X$ if ${\cal U}\!\in\! {\bf\Gamma}(Y\!\times\! X)$ and ${\bf\Gamma}(X)\! =\!\{ {\cal U}_y\mid y\!\in\! Y\}$.\smallskip

\noindent (b) We say that ${\cal U}$ is {\bf uniformly} $Y$-{\bf universal for the} $\bf\Gamma$
{\bf subsets of} $X$ if ${\cal U}$ is $Y$-universal for the $\bf\Gamma$ subsets of $X$ and, for each
$S\!\in\! {\bf\Gamma}(\omega^\omega\!\times\! X)$, there is a Borel map
$b\! :\!\omega^\omega\!\rightarrow\! Y$ such that $S_\alpha\! =\! {\cal U}_{b(\alpha )}$ for each
$\alpha\!\in\!\omega^\omega$.\end{defi}

\begin{cor} \label{univ} Let  $\cal G$ be a $G_\delta$ subset of $2^\omega$ with
$\lambda ({\cal G})\! =\! 0$. Then the set
$\{ (f,\beta )\!\in\! {\cal M}\!\times\! {\cal G}\mid\beta\!\in\! D(f)\}$ is $\cal M$-universal for the
$\borathree$ subsets of $\cal G$.\end{cor}

\vfill\eject

 In fact, we prove an effective and uniform version of the implication (a) $\Rightarrow$ (b) in Theorem \ref{main}. In particular, we can associate, via a Borel map $F$, a martingale to a code
$\alpha$ of an arbitrary $G_\delta$ subset $G$ of $\cal G$ (as in the previous corollary), in such a way that  $G\! =\! D\big( F(\alpha )\big)$. A consequence of this is the following:

\begin{thm} \label{completepi11} The set $\cal P$ of everywhere converging martingales is $\ca$-complete.\end{thm}

 These statements are in the spirit of some results concerning the differentiability of functions due to Zahorski and Mazurkiewicz (see Section 4 for details). In fact, $\cal P$ is $\ca$-complete in a uniform way, which allows to derive some universal and complete sets for the whole projective hierarchy, in spaces of continous functions, starting from $\cal P$. More precisely, let
$P_1\! :=\! [0,1]^{2^{<\omega}}$ and $C_1\! :=\! {\cal P}$. We define, for each natural number $n\!\geq\! 1$,\bigskip

\noindent $\bullet$ the space $P_{n+1}\! :=\! {\cal C}(2^\omega ,P_n)$ of continuous functions from $2^\omega$ into $P_n$, equipped with the topology of uniform convergence (inductively),\bigskip

\noindent $\bullet$ $C_{n+1}\! :=\!
\{ h\!\in\! P_{n+1}\mid\forall\beta\!\in\! 2^\omega ~~h(\beta )\!\notin\! C_n\}$ (inductively),\bigskip

\noindent $\bullet$ $U_n\! :=\!\{ (h,\beta )\!\in\! P_{n+1}\!\times\! 2^\omega\mid h(\beta )\!\in\! C_n\}$.\bigskip

 We prove the following:

\begin{thm} \label{projapp} Let $n\!\geq\! 1$ be a natural number. Then\smallskip

\noindent (a) the set $U_n$ is uniformly $P_{n+1}$-universal for the ${\bf\Pi}^1_n$ subsets of
$2^\omega$,\smallskip

\noindent (b) the set $C_n$ is ${\bf\Pi}^1_n$-complete.\end{thm}

 In fact, our method is more general and works if we start with a $\ca$ set which is complete in a uniform way.\bigskip

 Let $f$ be a martingale. As $D(f)$ has $\lambda$-measure zero, we can associate to $f$ the partial function $\psi (f)$ defined $\lambda$-almost everywhere by
$\psi (f)(\beta )\! :=\!\mbox{lim}_{l\rightarrow\infty}~f(\beta\vert l)$. The partial function $\psi (f)$ will be called the {\bf associated partial function}. The martingale $f$ is in $\cal P$ if and only if
$\psi (f)$ is total, in which case $\psi (f)$ is called the {\bf associated function}. Using the work in
[B-Ka-L] and [K2] about spaces of continuous functions, we prove the following:

\begin{thm} \label{completes} (a) The set of sequences of everywhere converging martingales whose associated functions converge pointwise is
${\bf\Pi}^1_1$-complete.\smallskip

\noindent (b) The set of sequences of everywhere converging martingales whose associated functions converge pointwise to zero is ${\bf\Pi}^1_1$-complete.\smallskip

\noindent (c) The set of sequences of everywhere converging martingales having a subsequence whose associated functions converge pointwise to zero is ${\bf\Sigma}^1_2$-complete.\end{thm}

\section{$\!\!\!\!\!\!$ $\borathree$ sets of measure zero}

\bf Notation.\rm\ In the sequel, $B$ will be a Borel subset of $2^\omega$, and $M$ will be a 
$\lambda$-measurable subset of $2^\omega$. If $\beta\!\in\! 2^\omega$, then the {\bf density of} $M$ {\bf at} $\beta$ is the number $d(M,\beta )\! :=\!
\mbox{lim}_{l\rightarrow\infty}~\frac{\lambda (M\cap N_{\beta\vert l})}{\lambda (N_{\beta\vert l})}$  when it is defined. Note that $d(B,\beta )\! =\! 1$ if $\beta\!\in\! B$ and $B$ is open. We first recall the Lebesgue density theorem (see 17.9 in [K2]).

\begin{thm} \label{densityone} (Lebesgue) The equality 
$\lambda (M)\! =\!\lambda\big(\{\beta\!\in\! M\mid d(M,\beta )\! =\! 1\}\big)$ holds for any 
$\lambda$-measurable subset $M$ of $2^\omega$.\end{thm}

 The reader should see [C] for the next lemma. We include a proof to be self-contained and also because we will prove an effective and uniform version of it later.

\begin{lem} \label{LM} (Lusin-Menchoff) Let $F$ be a closed subset of $2^\omega$, and
$M\!\supseteq\! F$ be a $\lambda$-measurable subset of $2^\omega$ such that 
$d(M,\beta )\! =\! 1$ for each $\beta\!\in\! F$.  Then there is a closed subset $C$ of $2^\omega$ such that\smallskip

(1) $F\!\subseteq\! C\!\subseteq\! M$,\smallskip

(2) $d(M,\beta )\! =\! 1$ for each $\beta\!\in\! C$,\smallskip

(3) $d(C,\beta )\! =\! 1$ for each $\beta\!\in\! F$.\end{lem}

\noindent\bf Proof.\rm\ If $F$ is $2^\omega$, then we can take $C\! :=\! F$. So we may assume that $F$ is not $2^\omega$. We set $s^-\! :=\! s\vert (\vert s\vert\! -\! 1)$ if 
$\emptyset\!\not=\! s\!\in\! 2^{<\omega}$. Note that $\neg F$ is the disjoint union of the elements of a sequence $(N_{s_n})_{n\in\omega}$, where $N_{s_n^-}\cap F\!\not=\!\emptyset$ for each 
$n\!\in\!\omega$. Fix $n\!\in\!\omega$. By Theorem \ref{densityone},
$$\lambda (M\cap N_{s_n})\! =\!
\lambda\big(\{\beta\!\in\! M\cap N_{s_n}\mid d(M\cap N_{s_n},\beta )\! =\! 1\}\big) .$$
The regularity of $\lambda$ gives a closed subset $F_n$ of $2^\omega$ contained in
$\{\beta\!\in\! M\cap N_{s_n}\mid d(M\cap N_{s_n},\beta )\! =\! 1\}$ such that
$\lambda (F_n)\!\geq\! (1\! -\! 2^{-n})\lambda (M\cap N_{s_n})$. We set
$C\! :=\! F\cup\bigcup_{n\in\omega}~F_n$, which is closed since $\vert s_n\vert\rightarrow\infty$.\bigskip

 As Conditions (1) and (2) are clearly satisfied, pick $\beta\!\in\! F$. Note that
$$\begin{array}{ll}
\lambda (N_{\beta\vert l}\!\setminus\! C)\!\!\!\!\!
& =\!\Sigma_{s_n\supseteq\beta\vert l}~\lambda (N_{s_n}\!\setminus\! C)\cr
& \leq\!\Sigma_{s_n\supseteq\beta\vert l}~\lambda (N_{s_n}\!\setminus\! F_n)\cr
& \leq\!\Sigma_{s_n\supseteq\beta\vert l}~2^{-n}\lambda (M\cap N_{s_n})\! +\!
\Sigma_{s_n\supseteq\beta\vert l}~\lambda (N_{s_n}\!\setminus\! M)\cr
& \leq\!\Sigma_{s_n\supseteq\beta\vert l}~2^{-n}\lambda (N_{s_n})\! +\!
\lambda (N_{\beta\vert l}\!\setminus\! M).
\end{array}$$
This implies that the limit of
$\frac{\lambda (N_{\beta\vert l}\setminus C)}{\lambda (N_{\beta\vert l})}$ is zero since
$d(M,\beta )\! =\! 1$.\hfill{$\square$}\bigskip

 The next topology is considered in [Lu-Ma-Z], see Chapter 6.

\begin{defi} The $\tau$-{\bf topology} on $2^\omega$ is generated by
$${\cal F}\! :=\!\{ M\!\subseteq\! 2^\omega\mid M\mbox{ is }\lambda\mbox{-measurable}\wedge
\forall\beta\!\in\! M~~d(M,\beta )\! =\! 1\} .$$
\end{defi}

 The next result is proved in [Lu-Ma-Z], but in a much more abstract way. This is the reason why we include a much more direct proof here, since it is not too long.

\begin{lem} \label{measurable} The family $\cal F$ is a topology. In particular, any $\tau$-open set is $\lambda$-measurable.\end{lem}

\noindent\bf Proof.\rm\ Note first that $\cal F$ is closed under finite intersections, so that it is a basis for the $\tau$-topology. Indeed, let $M,M'$ be in $\cal F$, and $\beta\!\in\! M\cap M'$. Then we use the facts that
$$\lambda (M\cap M'\cap N_{\beta\vert l})\! =\!\lambda (M\cap N_{\beta\vert l})\! -\!
\lambda\big( (M\cap N_{\beta\vert l})\!\setminus\! M'\big)$$
and $\lambda\big( (M\cap N_{\beta\vert l})\!\setminus\! M'\big)\!\leq\!
\lambda (N_{\beta\vert l}\!\setminus\! M')$.\bigskip

 Let $\cal H$ be a subfamily of $\cal F$, and $H\! :=\!\cup\cal H$. We claim that there is a countable subfamily $\cal C$ of $\cal H$ such that $m\! :=\!\mbox{sup}\{\lambda (\cup {\cal D})\mid {\cal D}\!\subseteq\! {\cal H}\mbox{ countable}\}\! =\!\lambda (\cup {\cal C})$. Indeed, for each $n\!\in\!\omega$ there is ${\cal D}_n\!\subseteq\! {\cal H}$ countable such that
 $\lambda (\cup {\cal D}_n)\! >\! m\! -\! 2^{-n}$, and
 ${\cal C}\! :=\!\bigcup_{n\in\omega}\ {\cal D}_n$ is suitable. Let $C\! :=\!\cup {\cal C}$.\bigskip

 Let $\beta\!\in\! H$, and $M$ in $\cal H$ with $\beta\!\in\! M$. Note that 
$\lambda (M\cup C)\! =\!\lambda (C)$ (consider the family ${\cal C}\cup\{ M\}$). Thus 
$\lambda (M\!\setminus\! C)\! =\! 0$. As $d(M,\beta )\! =\! 1$, the equality $d(M\cap C,\beta )\! =\! 1$ holds, and 
$d(\neg C,\beta )\! =\! 0$. This implies that $H\!\setminus\! C$ is contained in 
$\{\beta\!\notin\! C\mid d(\neg C,\beta )\! <\! 1\}$, which has $\lambda$-measure zero by Theorem \ref{densityone}. Therefore $H\!\setminus\! C$ has $\lambda$-measure zero and 
$H\! =\! C\cup (H\!\setminus\! C)$ is $\lambda$-measurable.\bigskip

 Pick $\beta\!\in\! H$, and $M\!\in\! {\cal H}$ with $\beta\!\in\! M$. Then $d(M,\beta )\! =\! 1$, and thus $d(H,\beta )\! =\! 1$. Therefore $H\!\in\! {\cal F}$. This finishes the proof.\hfill{$\square$}\bigskip

 The next lemma is in the style of Urysohn's theorem (see [Lu] for its version on the real line). We include a proof to be self-contained and also because we will prove an effective and uniform version of it later.

\begin{lem} \label{map} Let $C$ be a closed subset of $2^\omega$, and $G$ be a
$G_\delta$ subset of $2^\omega$ disjoint from $C$ such that $\lambda (G)\! =\! 0$. Then there is a $\tau$-continuous map $h\! :\! 2^\omega\!\rightarrow\! [0,1]$ such that $h_{\vert C}\!\equiv\! 0$ and $h_{\vert G}\!\equiv\! 1$.\end{lem}

\noindent\bf Proof.\rm\ Let $(F_n)_{n\in\omega}$ be an increasing sequence of closed subsets of $2^\omega$ with union $\neg G$ and $F_0\! =\! C$. We first construct a sequence
$(C_{\frac{1}{2^n}})_{n\in\omega}$ of closed subsets of $2^\omega$ with $F_n\!\subseteq\!
C_{\frac{1}{2^n}}\!\subseteq\!\neg G$, $C_{\frac{1}{2^n}}\!\subseteq\! C_{\frac{1}{2^{n+1}}}$, and
$d(C_{\frac{1}{2^{n+1}}},\beta )\! =\! 1$ for each $\beta\!\in\! C_{\frac{1}{2^n}}$. We first apply Lemma \ref{LM} to $F\! :=\! F_0$ and $M\! :=\!\neg G$, which gives
$F_0\!\subseteq\! C_1\!\subseteq\!\neg G$. Then, inductively, we apply Lemma \ref{LM} to
$F\! :=\! C_{\frac{1}{2^n}}\cup F_{n+1}$ and $M\! :=\!\neg G$, which gives
$C_{\frac{1}{2^n}}\cup F_{n+1}\!\subseteq\! C_{\frac{1}{2^{n+1}}}\!\subseteq\!\neg G$ such that
$d(C_{\frac{1}{2^{n+1}}},\beta )\! =\! 1$ for each $\beta\!\in\! C_{\frac{1}{2^n}}$.\bigskip

 Then we construct $C_{\frac{2k+1}{2^n}}$, for $0\! <\! k\! <\! 2^{n-1}$ and $n\!\geq\! 2$. This will give us a family $(C_{\frac{k}{2^n}})_{n\in\omega ,0<k\leq 2^n}$ of closed subsets of $2^\omega$. We want to ensure that $C_\zeta\!\subseteq\! C_{\zeta'}$ and
$d(C_{\zeta'},\beta )\! =\! 1$ for each $\beta\!\in\! C_\zeta$ if $\zeta'\! <\!\zeta$. We proceed by induction on $n$. We apply Lemma \ref{LM} to $F\! :=\! C_{\frac{k+1}{2^{n-1}}}$ and
$M\! :=\! C_{\frac{k}{2^{n-1}}}$, which gives $C_{\frac{2k+1}{2^n}}$ such that
$C_{\frac{k+1}{2^{n-1}}}\!\subseteq\! C_{\frac{2k+1}{2^n}}\!\subseteq\! C_{\frac{k}{2^{n-1}}}$,
$d(C_{\frac{k}{2^{n-1}}},\beta )\! =\! 1$ for each $\beta\!\in\! C_{\frac{2k+1}{2^n}}$, and
$d(C_{\frac{2k+1}{2^n}},\beta )\! =\! 1$ for each $\beta\!\in\! C_{\frac{k+1}{2^{n-1}}}$. This allows us to define $\tilde h$ by
$$\tilde h(\beta )\! :=\!\left\{\!\!\!\!\!\!
\begin{array}{ll}
& 0\mbox{ if }\beta\!\in\! G\mbox{,}\cr
& \mbox{sup}\{\zeta\mid\beta\!\in\! C_\zeta\}\mbox{ if }\beta\!\notin\! G.
\end{array}
\right.$$
It remains to see that $\tilde h$ is $\tau$-continuous (and then we will set
$h(\beta )\! :=\! 1\! -\!\tilde h(\beta )$). So let $b\!\in\! (0,1]$, and $\beta\!\in\! 2^\omega$ with $\tilde h(\beta )\! <\! b$. Note that there is $\zeta\! <\! b$ with
$\tilde h(\beta )\! <\!\zeta$, so that $\beta\!\notin\! C_\zeta$. If $\gamma\!\notin\! C_\zeta$, then $\tilde h(\gamma )\!\leq\!\zeta\! <\! b$, so that $\neg C_\zeta$ is an open (and thus $\tau$-open since the $\tau$-topology is finer than the usual one) neighborhood of $\beta$ on which
$\tilde h\! <\! b$. In particular, $\tilde h$ is Borel.

\vfill\eject

 Now let $a\!\in\! [0,1)$. It is enough to see that
$B\! :=\!\{\gamma\!\in\! 2^\omega\mid\tilde h(\gamma )\! >\! a\}$ is $\tau$-open. So assume that
$\tilde h(\gamma )\! >\! a$. Note that there are $\zeta\! >\!\zeta'\! >\! a$ with
$\tilde h(\gamma )\! >\!\zeta$, so that $\gamma\!\in\! C_\zeta\!\subseteq\! C_{\zeta'}\!\subseteq\! B$. Thus $d(C_{\zeta'},\gamma )\! =\! 1$, by construction of the family. As $\tilde h$ is Borel, $B$ is Borel, $d(B,\gamma )$ is defined and equal to $1$.\hfill{$\square$}\bigskip

\noindent\bf Remark.\rm\ We in fact proved that $h$ is lower semi-continuous.\bigskip

\noindent\bf Notation.\rm\ If $h\! :\! 2^\omega\!\rightarrow\! [0,1]$ is a $\lambda$-measurable map and $s\!\in\! 2^{<\omega}$, then we set
$\fint_{N_s}h~d\lambda\! :=\!\frac{\int_{N_s}h~d\lambda}{\lambda (N_s)}$.

\begin{lem} \label{moy} Let $h\! :\! 2^\omega\!\rightarrow\! [0,1]$ be a $\tau$-continuous map, and $\beta\!\in\! 2^\omega$. Then
$$\lim_{l\rightarrow\infty}~\fint_{N_{\beta\vert l}}h~d\lambda\! =\! h(\beta ).$$\end{lem}

\noindent\bf Proof.\rm\ Let $\varepsilon\! >\! 0$, and 
$\beta\!\in\! M\! :=\! h^{-1}\Big( B\big( h(\beta),\varepsilon\big)\Big)$. Note that 
$d(M,\gamma )\! =\! 1$ for each $\gamma\!\in\! M$ since $h$ is $\tau$-continuous. As $h$ is $\lambda$-measurable, we can write
$$\int_{N_{\beta\vert l}}h~d\lambda\! =\!\int_{M\cap N_{\beta\vert l}}h~d\lambda\! +\!
\int_{N_{\beta\vert l}\setminus M}h~d\lambda .$$
Note that $0\!\leq\!\int_{N_{\beta\vert l}\setminus M}h~d\lambda\!\leq\!
\lambda (N_{\beta\vert l}\!\setminus\! M)$, so that
$0\!\leq\!\fint_{N_{\beta\vert l}\setminus M}h~d\lambda\!\leq\!
\frac{\lambda (N_{\beta\vert l}\setminus M)}{\lambda (N_{\beta\vert l})}\rightarrow 0$. Similarly,
$$\fint_{M\cap N_{\beta\vert l}}h~d\lambda\!\in\!\big[\big( h(\beta )\! -\!\varepsilon\big)
\frac{\lambda (M\cap N_{\beta\vert l})}{\lambda (N_{\beta\vert l})},
\big( h(\beta )\! +\!\varepsilon\big)
\frac{\lambda (M\cap N_{\beta\vert l})}{\lambda (N_{\beta\vert l})}\big]\mbox{,}$$ 
and we are done since $\frac{\lambda (M\cap N_{\beta\vert l})}{\lambda (N_{\beta\vert l})}$ tends to $1$ as $l$ tends to $\infty$.\hfill{$\square$}\bigskip

 Now we come to our main lemma, inspired by Zahorski (see [Za]).

\begin{lem} \label{Gdelta} Let $G$ be a $G_\delta$ subset of $2^\omega$ with $\lambda$-measure zero. Then there is a martingale $f$ with $G\! =\! D(f)$ and
$\{\mbox{osc}(f,\beta )\mid\beta\!\in\! 2^\omega\}\!\subseteq\!\{ 0\}\cup [\frac{1}{2},1]$.\end{lem}

\noindent\bf Proof.\rm\ Let $(G_n)_{n\in\omega}$ be a decreasing sequence of open subsets of $2^\omega$ with intersection $G$ and $G_0\! =\! 2^\omega$.\bigskip

\noindent $\bullet$ We construct $g_n\! :\! 2^\omega\!\rightarrow\! [0,1]$, open subsets
$G^*_n,G^{**}_n$ of $2^\omega$, and a sequence $(s^n_j)_{j\in I_n}$ of pairwise incompatible finite binary sequences, by induction on $n\!\in\!\omega$, such that, if
$S_n\! :=\!\Sigma_{j\leq n}~(-1)^jg_j$,
$$\begin{array}{ll}
& (1)~G\!\subseteq\! G^*_{n+1}\!\subseteq\! G^{**}_n\! =\!\bigcup_{j\in I_n}~N_{s^n_j}\!\subseteq\! G^*_n\!\subseteq\! G_n\ \wedge\ G^*_0\! =\! 2^\omega\mbox{,}\cr
& (2)~{g_n}_{\vert G}\!\equiv\! 1\ \wedge\ {g_n}_{\vert\neg G^*_n}\!\equiv\! 0\mbox{,}\cr
& (3)~g_n\mbox{ is }\tau\mbox{-continuous}\mbox{,}\cr
& (4)~g_{n+1}\!\leq\! g_n\mbox{,}\cr
& (5)~\lambda (G^*_{n+1}\cap N_{s^n_j})\! <\! 2^{-n-3}\lambda (N_{s^n_j})\mbox{,}\cr
& (6)~\vert\fint_{N_{s^n_j}}S_n~d\lambda\! -\! S_n(\beta )\vert\! <\! 2^{-3}\mbox{ if }
\beta\!\in\! G\cap N_{s^n_j}.
\end{array}$$
We set $g_0\! :\equiv\! 1$, $G^*_0,G^{**}_0\! :=\! 2^\omega$, $I_0\! :=\!\{ 0\}$ and
$s^0_0\! :=\!\emptyset$. Assume that our objects are constructed up to $n$. We first construct an open subset $G^*_{n+1}$ of $2^\omega$ with
$G\!\subseteq\! G^*_{n+1}\!\subseteq\! G^{**}_n\cap G_{n+1}$ such that
$$\lambda (G^*_{n+1}\cap N_{s^n_j})\! <\! 2^{-n-3}\lambda (N_{s^n_j})$$
if $j\!\in\! I_n$. For each $j\!\in\! I_n$, there is an open set $O_j$ with
$G\cap N_{s^n_j}\!\subseteq\! O_j\!\subseteq\! G_{n+1}\cap N_{s^n_j}$ such that
$\lambda (O_j)\! <\! 2^{-n-3}\lambda (N_{s^n_j})$. We then set
$G^*_{n+1}\! :=\!\bigcup_{j\in I_n}~O_j$.

\vfill\eject

 We now apply Lemma \ref{map} to $C\! :=\!\neg G^*_{n+1}$ and $G$, which gives a $\tau$-continuous map $h\! :\! 2^\omega\!\rightarrow\! [0,1]$ with
$h_{\vert\neg G^*_{n+1}}\!\equiv\! 0$ and $h_{\vert G}\!\equiv\! 1$. We set
$g_{n+1}\! :=\!\mbox{min}(g_n,h)$, so that $g_{n+1}$ satisfies (2)-(4).\bigskip

 By Lemma \ref{moy},
$\mbox{lim}_{l\rightarrow\infty}~\fint_{N_{\beta\vert l}}S_{n+1}~d\lambda\! =\! S_{n+1}(\beta )$ for each $\beta\!\in\! G$. This gives $l(\beta )\!\in\!\omega$ minimal with
$\vert\fint_{N_{\beta\vert l(\beta )}}S_{n+1}~d\lambda\! -\! S_{n+1}(\beta )\vert\! <\! 2^{-3}$ and 
$N_{\beta\vert l(\beta )}\!\subseteq\! G^*_{n+1}$. The set $G^{**}_{n+1}$ is the union of the 
$N_{\beta\vert l(\beta )}$'s, which defines $I_{n+1}$ and $(s^{n+1}_j)_{j\in I_{n+1}}$ 
($S_{n+1}(\beta )$ is 0 if $n$ is even and 1 otherwise when $\beta\!\in\! G$).\bigskip

\noindent $\bullet$ We then define a partial map $f_\infty\! :\! 2^\omega\!\rightarrow\! [0,1]$ by
$f_\infty\! :=\!\Sigma_{j\in\omega}~(-1)^jg_j$. If $\beta\!\in\! G$, then $S_n(\beta )$ takes alternatively the values $1$ and $0$, depending on the parity of $n$, so that $f_\infty (\beta )$ is not defined. If $\beta\!\notin\! G$, then there is $n$ such that
$\beta\!\in\!\neg G^*_{n+1}\!\subseteq\!\neg G^*_{n+2}\!\subseteq\! ...$ This implies that
$f_\infty (\beta )$ is defined and equal to $S_n(\beta )$. As
$0\!\leq\!\Sigma_{p\leq q}~(g_{2p}\! -\! g_{2p+1})\! =\! S_{2q+1}\!\leq\! S_{2q}\! =\! 
g_0\! +\!\Sigma_{1\leq p\leq q}~(g_{2p}\! -\! g_{2p-1})\!\leq\! g_0$, $f_\infty$ takes values in
$[0,1]$. So $f_\infty$ is a partial $\lambda$-measurable map defined $\lambda$-almost everywhere since $\lambda (G)\! =\! 0$ (we use Lemma \ref{measurable}).\bigskip

\noindent $\bullet$ This allows us to define $f\! :\! 2^{<\omega}\!\rightarrow\! [0,1]$ by
$f(s)\! :=\!\fint_{N_s}f_\infty ~d\lambda$. As $\lambda (N_s)\! =\! 2\lambda (N_{s\varepsilon})$ for each $\varepsilon\!\in\! 2$,
$f(s)\! =\!\fint_{N_s}f_\infty ~d\lambda\! =\!\frac{\int_{N_{s0}}f_\infty ~d\lambda +
\int_{N_{s1}}f_\infty ~d\lambda}{\lambda (N_s)}\! =\!\frac{f(s0)}{2}\! +\!\frac{f(s1)}{2}$
and $f$ is a martingale.\bigskip

\noindent $\bullet$ If $\beta\!\notin\! G$, then there is $n$ with
$\beta\!\in\! G^*_n\!\setminus\! G^*_{n+1}$, so that $f_\infty (\beta )\! =\! S_n(\beta )$. By Lemma
\ref{moy}, $k\!\geq\! n$ implies that
$\mbox{lim}_{l\rightarrow\infty}~\fint_{N_{\beta\vert l}}S_{k+1}~d\lambda\! =\!
S_{k+1}(\beta )\! =\! S_n(\beta )$ since $S_{k+1}$ is $\tau$-continuous. Note that, for each $k\!\geq\! n$,
$$\begin{array}{ll}
\big\vert\int_{N_{\beta\vert l}}(f_\infty\! -\! S_{k+1}) ~d\lambda\big\vert\!\!\!\!
& \leq\!\lambda (G^*_{k+2}\cap N_{\beta\vert l})\cr
& \leq\!\Sigma_{\beta\vert l\subseteq s^{k+1}_j}~\lambda (G^*_{k+2}\cap N_{s^{k+1}_j})\cr
& \leq\!\Sigma_{\beta\vert l\subseteq s^{k+1}_j}~2^{-k-4}\lambda (N_{s^{k+1}_j})\cr
& \leq\!\lambda (N_{\beta\vert l})2^{-k-4}.
\end{array}$$
Moreover,
$$\begin{array}{ll}
\vert f(\beta\vert l)\! -\! f_\infty (\beta )\vert\! =\!
\vert\fint_{N_{\beta\vert l}}f_\infty ~d\lambda\! -\! f_\infty (\beta )\vert\!\!\!\!
& =\!\vert\fint_{N_{\beta\vert l}}\big( f_\infty\! -\! S_{k+1}\big) ~d\lambda
\! +\!\fint_{N_{\beta\vert l}}S_{k+1}~d\lambda\! -\! S_{k+1}(\beta )\vert\cr
& \leq\! 2^{-k-4}\! +\!
\vert\fint_{N_{\beta\vert l}}S_{k+1}~d\lambda\! -\! S_{k+1}(\beta )\vert\mbox{,}
\end{array}$$
so that $\mbox{lim}_{l\rightarrow\infty}~f(\beta\vert l)\! =\! f_\infty (\beta )$,
$\mbox{osc}(f,\beta )\! =\! 0$ and $\beta\!\notin\! D(f)$.\bigskip

\noindent $\bullet$ If $\beta\!\in\! G$ and $n\!\in\!\omega$, then there is $j\!\in\!\omega$ with
$\beta\!\in\! N_{s^n_j}$. Note that
$$f(s^n_j)\! =\!\fint_{N_{s^n_j}}f_\infty ~d\lambda\! =\!\fint_{N_{s^n_j}}S_n~d\lambda\! +\!
\fint_{N_{s^n_j}}(f_\infty\! -\! S_n)~d\lambda$$
and $\vert\int_{N_{s^n_j}}(f_\infty\! -\! S_n)~d\lambda\vert\!\leq\!
\lambda (G^*_{n+1}\cap N_{s^n_j})\! <\!\frac{1}{8}\lambda (N_{s^n_j})$, so that
$\vert\fint_{N_{s^n_j}}(f_\infty\! -\! S_n)~d\lambda\vert\! <\!\frac{1}{8}$. By (6), $\vert f(s^n_j)\! -\! S_n(\beta )\vert\! <\!\frac{1}{8}\! +\!\frac{1}{8}\! =\!\frac{1}{4}$. As
$S_n(\beta )$ takes infinitely often the values $1$ and $0$, $\mbox{osc}(f,\beta )\!\geq\!\frac{1}{2}$ and $\beta\!\in\! D(f)$.\hfill{$\square$}\bigskip

 The main result will be a consequence of the main lemma and the following.

\begin{lem} \label{union} Let $(f_n)_{n\in\omega}$ be a sequence of martingales such that
$$\{\mbox{osc}(f_n, \beta )\mid (n,\beta )\!\in\!\omega\!\times\! 2^\omega\}\!\subseteq\!
\{ 0\}\cup [\frac{1}{2}, 1].$$
Then there is a martingale $f$ with $D(f)\! =\!\bigcup_{n\in\omega}~D(f_n)$.\end{lem}

\noindent\bf Proof.\rm\ We first observe the following facts. Let
$g,h\! :\! 2^{<\omega}\!\rightarrow\!\mathbb{R}$ be bounded, $\beta\!\in\! 2^\omega$ and $a\!\in\!\mathbb{R}$.\bigskip

\noindent (1) $\mbox{osc}(g\! +\! h,\beta )\!\leq\!\mbox{osc}(g,\beta )\! +\!
\mbox{osc}(h,\beta )$.\bigskip

 This comes from the triangle inequality.\bigskip

\noindent (2) $\mbox{osc}(ag,\beta )\! =\!\vert a\vert\!\cdot\!\mbox{osc}(g,\beta )$.\bigskip

\noindent (3) $\mbox{osc}(g\! +\! h,\beta )\! =\!\mbox{osc}(h,\beta )$ if 
$\mbox{osc}(g,\beta )\! =\! 0$.\bigskip

 By (1), $\mbox{osc}(h,\beta )\!\leq\!\mbox{osc}(g\! +\! h,\beta )\! +\!
\mbox{osc}(-g,\beta )\! =\!\mbox{osc}(g\! +\! h,\beta )\!\leq\!\mbox{osc}(g,\beta )\! +\!
\mbox{osc}(h,\beta )\! =\!\mbox{osc}(h,\beta )$, so that 
$\mbox{osc}(h,\beta )\! =\!\mbox{osc}(g\! +\! h,\beta )$.\bigskip

\noindent (4) $\mbox{osc}(g,\beta )\!\leq\! a$ if $g(\beta\vert l)\!\in\! [0, a]$ for each $l\!\in\!\omega$.\bigskip

\noindent $\bullet$ We set $D_n\! :=\! D(f_n)$ for each $n\!\in\!\omega$, and
$f\! :=\!\Sigma_{n\in\omega}~4^{-n}f_n$. Note that $f$ is defined and a martingale.\bigskip

\noindent $\bullet$ If $\beta\!\notin\!\bigcup_{n\in\omega}~D_n$, then 
$\mbox{osc}(f_n,\beta )\! =\! 0$ for each $n\!\in\!\omega$. In particular, 
$\mbox{osc}(4^{-n}f_n,\beta )\! =\! 0$ for each $n\!\in\!\omega$, by (2). Let $\varepsilon\! >\! 0$, and $M\!\in\!\omega$ with $\Sigma_{n>M}~4^{-n}\!\leq\!\varepsilon$. By (1), 
$\mbox{osc}(\Sigma_{n\leq M}~4^{-n}f_n,\beta )\! =\! 0$. By (3) and (4), 
$\mbox{osc}(f,\beta )\! =\!\mbox{osc}(\Sigma_{n>M}~4^{-n}f_n,\beta )\!\leq\!
\Sigma_{n>M}~4^{-n}\!\leq\!\varepsilon$. As $\varepsilon$ is arbitrary,
$\mbox{osc}(f,\beta )\! =\! 0$, $\beta\!\notin\! D(f)$, which shows that
$D(f)\!\subseteq\!\bigcup_{n\in\omega}~D_n$.\bigskip

\noindent $\bullet$ If $\beta\!\in\!\bigcup_{n\in\omega}~D_n$, then let $m$ be minimal such that 
$\beta\!\in\! D_m$. Note that
$$f\! =\!\Sigma_{n<m}~4^{-n}f_n\! +\! 4^{-m}f_m\! +\!\Sigma_{n>m}~4^{-n}f_n.$$ 
By (2) and (3), 
$\mbox{osc}(f,\beta )\! =\!\mbox{osc}(4^{-m}f_m\! +\!\Sigma_{n>m}~4^{-n}f_n,\beta )$. By (1), (2) and (4),
$$\mbox{osc}(f,\beta )\!\geq\!\mbox{osc}(4^{-m}f_m,\beta )\! -\!
\mbox{osc}(\Sigma_{n>m}~4^{-n}f_n,\beta )\!\geq\! 
4^{-m}\frac{1}{2}\! -\! 4^{-m}\frac{1}{3}\! >\! 0.$$
Thus $\beta\!\in\! D(f)$.\hfill{$\square$}

\section{$\!\!\!\!\!\!$ Effectivity and uniformity}

- We refer to [M] for the basic notions of effective descriptive set theory. We first recall some material present in it.\bigskip

$\bullet$ Let $(p_n)_{n\in\omega}$ be the sequence of prime numbers $2,3,...$\bigskip

$\bullet$ If $l\!\in\!\omega$ and $s\!\in\!\omega^l$, then
$\overline{s}\! :=<s(0),...,s(l\! -\! 1)>:=\! p_0^{s(0)+1}...p_{l-1}^{s(l-1)+1}\!\in\!\omega$ codes $s$ (if
$l\! =\! 0$, then $<>:=\! 1$).\bigskip

$\bullet$ If $\alpha\!\in\!\omega^\omega$ and $l\!\in\!\omega$, then
$\overline{\alpha}(l)\! :=<\alpha (0),...,\alpha (l\! -\! 1)>\in\omega$ codes
$\alpha\vert l\!\in\!\omega^l$, and
$\alpha^*$ is defined by removing the first coordinate:
$\alpha^*\! :=\!\big(\alpha (1),\alpha (2),...\big)$.

\vfill\eject

$\bullet$ If $\kappa\!\in\!\{ 2,\omega\}$, then
$<.,.>:\! (\kappa^\omega )^2\!\rightarrow\!\kappa^\omega$ is a recursive homeomorphism with inverse map $\alpha\!\mapsto\!\big( (\alpha )_0,(\alpha )_1\big)$ defined for example by
$(\alpha )_\varepsilon (n)\! :=\!\alpha (2n\! +\!\varepsilon )$ if
$(n,\varepsilon )\!\in\!\omega\!\times\! 2$ (we will also consider recursive homeomorphisms
$<.,.,.>:\! (\kappa^\omega )^3\!\rightarrow\!\kappa^\omega$ and
$<.,.,...>:\! (\kappa^\omega )^\omega\!\rightarrow\!\kappa^\omega$).\bigskip

$\bullet$ If $u\!\in\!\omega$, then $\mbox{Seq}(u)$ means that there are $l\!\in\!\omega$ and
$s\!\in\!\omega^l$ (denoted by $s(u)$) such that $u\! =<s(0),...,s(l\! -\! 1)>$. The natural number
$(u)_i$ is $s(i)$ if $i\! <\! l$, and $0$ otherwise. The number $l$ is the {\bf length} of $u$ and is denoted by $\mbox{lh}(u)$. If $k\!\leq\! l$, then $\underline{u}(k)\! :=<s(0),...,s(k\! -\! 1)>$, so that
$\underline{u}(l)\! =\! u$. The standard basic clopen set is
$N^u\! :=\!\{\beta\!\in\! 2^\omega\mid\forall i\! <\!\mbox{lh}(u)~~\beta (i)\! =\! (u)_i\}$. We set
$u^-\! :=<(u)_0,...,(u)_{\mbox{lh}(u)-2}>$ ($u^-\! :=<>$ if $\mbox{lh}(u)\!\leq\! 1$).\bigskip

$\bullet$ Let $X$ be a recursively presented Polish space. Then we will consider the effective basic open set $N(X,u)\! =\! B_X(r_{((u)_1)_0},\frac{((u)_1)_1}{((u)_1)_2+1})$.\bigskip

$\bullet$ Let $n\!\geq\! 1$ be a natural number. A subset $T$ of $\omega^n$ is a {\bf tree} if
$\mbox{Seq}(u_i)$ and $\mbox{lh}(u_i)\! =\!\mbox{lh}(u_0)$ for each $(u_0,...,u_{n-1})\!\in\! T$ and each $i\! <\! n$, and $\big(\underline{u_0}(k),...,\underline{u_{n-1}}(k)\big)\!\in\! T$ if $(u_0,...,u_{n-1})\!\in\! T$ and $k\!\leq\!\mbox{lh}(u_0)$.\bigskip

$\bullet$ The next result is a part of 4A.1 in [M].

\begin{thm} \label{repsemi} Let $m\!\geq\! 1$ be a natural number, and $B\!\in\!\Boraone\big(\omega^\omega\!\times\! (\omega^\omega )^m\big)$. Then there is a recursive subset $T$ of $\omega^\omega\!\times\!\omega^m$ such that 
$(\alpha ,\alpha_1,...,\alpha_m)\!\in\! B\Leftrightarrow\exists l\!\in\!\omega ~~\big(\alpha ,\underline{\alpha_1}(l),...,\underline{\alpha_m}(l)\big)\!\notin\! T$, and $T_\alpha\! :=\!\{ (u_0,...,u_{m-1})\!\in\!\omega^m\mid (\alpha ,u_0,...,u_{m-1})\!\in\! T\}$ is a tree for each $\alpha\!\in\!\omega^\omega$.\end{thm}

$\bullet$ The next result is a part of 4A.7 in [M].

\begin{thm} \label{repboreff} Let $X$ be a recursively presented Polish space and
$B\!\in\!\Borel (X)$. Then we can find a recursive function $\pi\! :\!\omega^\omega\!\rightarrow\! X$ and $C\!\in\!\Bormone (\omega^\omega )$ such that $\pi$ is injective on $C$ and $\pi [C]\! =\! B$.
\end{thm}

\noindent - We then recall some material from [L].\bigskip

\noindent\bf Notation.\rm\ Let $X$ be a recursively presented Polish space. Recall that there is a pair
$({\cal W}^X,{\cal C}^X)$ such that\bigskip

$\bullet$ ${\cal W}^X\!\subseteq\!\omega$ is a $\Ca$ set of codes for the $\Borel$ subsets of $X$,\smallskip

$\bullet$ ${\cal C}^X\!\subseteq\!\omega\!\times\! X$ is $\Ca$ and
$\Borel (X)\! =\!\{ {\cal C}^X_n\mid n\!\in\! {\cal W}^X\}$, which means that ${\cal C}^X$ is ``universal" for the $\Borel$ subsets of $X$,\smallskip

$\bullet$ the relation ``$n\!\in\! {\cal W}^X~\wedge ~(n,x)\!\notin\! {\cal C}^X$" is $\Ca$ in $(n,x)$.\bigskip

\noindent If $X\! =\!\omega^\omega\!\times\! 2^\omega$, then we simply write
$({\cal W},{\cal C})\! :=\! ({\cal W}^X,{\cal C}^X)$.\bigskip

 The next result will be extremely useful in the sequel.\bigskip

\noindent\bf The uniformization lemma.~\it Let $X,Y$ be recursively presented Polish spaces, and
$P\!\in\!\Ca (X\!\times\! Y)$. Then the set
$P^+\! :=\!\{ x\!\in\! X\mid\exists y\!\in\!\Borel (x)~~(x,y)\!\in\! P\}$ is $\Ca$, and there is a partial
$\Ca$-recursive map $f\! :\! X\!\rightarrow\! Y$ such that $\big( x,f(x)\big)\!\in\! P$ for each $x\!\in\! P^+$. If moreover $S\!\subseteq\! P^+$ is a $\Ana$ subset of $X$, then there is a total $\Borel$-recursive map $g\! :\! X\!\rightarrow\! Y$ such that $\big( x,g(x)\big)\!\in\! P$ for each $x\!\in\! S$.\rm\bigskip

\noindent - The following definition is inspired by 3H.1 in [M].

\begin{defi} (a) Let $\Gamma$ be a class of subsets of recursively presented Polish spaces, and
$\bf\Gamma$ be the associated boldface class. A system of sets
${\cal U}^X\!\in\!\Gamma (\omega^\omega\!\times\! X)$, where is $X$ is a recursively presented Polish space, is a {\bf nice parametrization} in $\Gamma$ for $\bf\Gamma$ if the following hold:\smallskip

(1) ${\bf\Gamma}(X)\! =\!\{ {\cal U}^X_\alpha\mid\alpha\!\in\!\omega^\omega\}$,\smallskip

(2) $\Gamma (X)\! =\!\{ {\cal U}^X_\alpha\mid\alpha\!\in\!\omega^\omega\mbox{ recursive}\}$,\smallskip

(3) if $X$ is a recursively presented Polish space, then there is
${\cal R}\! :\!\omega^\omega\!\times\!\omega^\omega\!\rightarrow\!\omega^\omega$ recursive such that $(\alpha ,\gamma ,x)\!\in\! {\cal U}^{\omega^\omega\times X}\Leftrightarrow
\big( {\cal R}(\alpha ,\gamma ),x\big)\!\in\! {\cal U}^X$ if
$(\alpha ,\gamma ,x)\!\in\!\omega^\omega\!\times\!\omega^\omega\!\times\! X$.\smallskip

\noindent (b) If $\cal U$ belongs to a nice parametrization, then we will say that $\cal U$ is a
{\bf good universal set }.\smallskip

\noindent (c) If $\cal U$ satisfies all these properties except maybe (3), then we will say that
$\cal U$ is a {\bf suitable universal set }.\end{defi}

 By 3E.2, 3F.6 and 3H.1 in [M], there is a nice parametrization in ${\it\Pi}^1_n$ for ${\bf\Pi}^1_n$, for each natural number $n\!\geq\! 1$.\bigskip

\noindent - We now recall two results that can essentially be found in [K1]. The first one is Theorem 2.2.3.(a) (see also [T1]).

\begin{thm} \label{hyperar} (Tanaka) Let $U\!\in\!\Ana (\omega^\omega\!\times\!\omega^\omega )$ be $\omega^\omega$-universal for the analytic subsets of $\omega^\omega$. Then
$L(U)\! :=\!\big\{ (\alpha ,p)\!\in\!\omega^\omega\!\times\!\omega\mid
\lambda (U_\alpha\cap 2^\omega )\! >\!\frac{(p)_0}{(p)_1+1}\big\}$ is $\Ana$.\end{thm}

\begin{cor} \label{Delta} Let $B\!\in\!\Borel (\omega^\omega\!\times\! 2^\omega )$.\smallskip

\noindent (a) The map $\lambda_B\! :\!\omega^\omega\!\rightarrow\!\mathbb{R}$ defined by
$\lambda_B(\alpha )\! :=\!\lambda (B_\alpha )$ is $\Borel$-recursive, and the partial function
${(n,\alpha )\!\mapsto\!\lambda ({\cal C}_{n,\alpha})}$ is $\Ca$-recursive on its domain
${\cal W}\!\times\!\omega^\omega$.\smallskip

\noindent (b) Let $D\!\subseteq\!\omega$, $O_0\!\in\!\Ana (\omega\!\times\!\omega^\omega\!\times\! 2^\omega )$, and 
$O_1\!\in\!\Ca (\omega\!\times\!\omega^\omega\!\times\! 2^\omega )$ be such that $\lambda\big( (O_0)_{n,\alpha}\big)\! =\!\lambda\big( (O_1)_{n,\alpha}\big)$ if $n\!\in\! D$. Then the partial map $\lambda_O\! :\! D\!\times\!\omega^\omega\!\rightarrow\!\mathbb{R}$ defined by 
$\lambda_O(n,\alpha )\! :=\!\lambda\big( (O_0)_{n,\alpha}\big)$ is $\Ana$-recursive and 
$\Ca$-recursive on its domain.\smallskip

\noindent (c) The partial map $d_B\! :\!\omega^\omega\!\times\! 2^\omega\!\rightarrow\!\mathbb{R}$ defined by $d_B(\alpha ,\beta )\! :=\! d(B_\alpha ,\beta )$ is $\Borel$-recursive, and the partial map $(n,\alpha ,\beta )\!\mapsto\! d({\cal C}_{n,\alpha},\beta )$ is $\Ca$-recursive on its $\Ca$ domain
$$\{ (n,\alpha ,\beta )\!\in\! {\cal W}\!\times\!\omega^\omega\!\times\! 2^\omega\mid
d({\cal C}_{n,\alpha},\beta )\mbox{ exists}\} .$$
(d) Let $h\! :\!\omega^\omega\!\times\! 2^\omega\!\rightarrow\!\mathbb{R}$ be $\Borel$-recursive taking values in $[0,1]$. Then the partial map 
$i_h\! :\!\omega^\omega\!\times\!\omega\!\rightarrow\!\mathbb{R}$ defined by
$i_h(\alpha ,u)\! :=\!\int_{N^u}h(\alpha ,.)~d\lambda$ is $\Borel$-recursive on its $\Borone$ domain
$\omega^\omega\!\times\!\{ u\!\in\!\omega\mid\mbox{Seq}(u)\}$.\end{cor}

\noindent\bf Proof.\rm\ (a) It is enough to see that the relations
$P_B(\alpha ,p)\Leftrightarrow\lambda (B_\alpha )\! >\!
r_p\! :=\! (-1)^{(p)_0}\!\cdot\!\frac{(p)_1}{(p)_2+1}$ and
$$Q_B(\alpha ,p)\Leftrightarrow\lambda (B_\alpha )\! <\! r_p$$
are $\Borel$ to see that $\lambda_B$ is $\Borel$-recursive. Note that there is 
$\phi\! :\!\omega^2\!\rightarrow\!\omega$ recursive with $r_{\phi (p,l)}\! =\! r_p\! -\!\frac{1}{l+1}$. Thus
$$\begin{array}{ll}
Q_B(\alpha ,p)\!\!\!\!
& \Leftrightarrow\exists l\!\in\!\omega ~~
\lambda (B_\alpha )\!\leq\! r_p\! -\!\frac{1}{l+1}\cr
& \Leftrightarrow\exists l\!\in\!\omega ~~
\neg\big(\lambda (B_\alpha )\! >\! r_p\! -\!\frac{1}{l+1}\big)\cr
& \Leftrightarrow\exists l\!\in\!\omega ~~
\neg P_B\big(\alpha ,\phi (p,l)\big)\mbox{,}
\end{array}$$
so that it is enough to see that $P_B$ is $\Borel$.

\vfill\eject

\noindent $\bullet$ Now let $S\!\in\!\Ana\big(\omega^\omega\!\times\! (\omega^\omega )^2\big)$ be a good $\omega^\omega$-universal for the analytic subsets of $(\omega^\omega )^2$. We set
$$U(\alpha ,\gamma )\Leftrightarrow S\big( (\alpha )_0,(\alpha )_1,\gamma\big)\mbox{,}$$
so that $U\!\in\!\Ana (\omega^\omega\!\times\!\omega^\omega )$ is $\omega^\omega$-universal for the analytic subsets of $\omega^\omega$. Let $A$ be a $\Ana$ subset of 
$\omega^\omega\!\times\! 2^\omega$. Then there is $\alpha_0\!\in\!\omega^\omega$ recursive with $A\! =\! S_{\alpha_0}$, so that
$$\gamma\!\in\! A_\alpha\Leftrightarrow (\alpha_0,\alpha ,\gamma )\!\in\! S\Leftrightarrow
(<\alpha_0,\alpha >,\gamma )\!\in\! U.$$
This implies that the relation $R_A(\alpha ,p)\Leftrightarrow\lambda (A_\alpha )\! >\! r_p$, equivalent to
$$\big( (p)_0\mbox{ is odd }\wedge ~(p)_1\! >\! 0\big) ~\vee ~\big( (p)_0\mbox{ is even }~\wedge ~
(<\alpha_0,\alpha >,<(p)_1,(p)_2>)\!\in\! L(U)\big)\mbox{,}$$
is $\Ana$, by Theorem \ref{hyperar}.\bigskip

\noindent $\bullet$ In particular, this applies to $A\! :=\! B$, so that $P_B$ is $\Ana$. Now note that
$$P_B(\alpha ,p)\Leftrightarrow\lambda\big( (\neg B)_\alpha\big)\! <\! 1\! -\! r_p
\Leftrightarrow Q_{\neg B}\big(\alpha ,\phi'(p)\big)\mbox{,}$$
for some $\phi'\! :\!\omega\!\rightarrow\!\omega$ is recursive, so that $P_B$ is $\Ca$ by the previous computation.\bigskip

\noindent $\bullet$ We set ${\cal C}'\! :=\!\big\{ (\gamma ,\beta )\!\in\!\omega^\omega\!\times\! 2^\omega\mid \gamma (0)\!\in\! {\cal W}~\wedge ~\big(\gamma (0),\gamma^*,\beta\big)\!\in\! {\cal C}\big\}$. As
${\cal C}'$ is $\Ca$,
$${\cal A}\! :=\!\big\{ (\alpha ,p)\!\in\!\omega^\omega\!\times\!\omega\mid
\lambda\big( (\neg {\cal C}')_\alpha\big)\! >\! r_p\big\}$$
is $\Ana$, by the previous discussion. So let $n\!\in\! {\cal W}$. Note that
$$\begin{array}{ll}
\lambda ({\cal C}_{n,\alpha})\! >\! r_p\!\!\!
& \Leftrightarrow\lambda (\neg {\cal C}_{n,\alpha})\! <\! 1\! -\! r_p\Leftrightarrow
\lambda\big( (\neg {\cal C}')_{n\alpha}\big)\! <\! 1\! -\! r_p\cr
& \Leftrightarrow\exists l\!\in\!\omega ~~
\lambda\big( (\neg {\cal C}')_{n\alpha}\big)\!\leq\! 1\! -\! r_p\! -\!\frac{1}{l+1}\Leftrightarrow
\exists l\!\in\!\omega ~~\big( n\alpha ,\phi''(p,l)\big)\!\notin\! {\cal A}\mbox{,}
\end{array}$$
for some recursive $\phi''\! :\!\omega^2\!\rightarrow\!\omega$. Similarly, the relation
``$\lambda ({\cal C}_{n,\alpha})\! <\! r_p$" is $\Ca$ in $(n,\alpha ,p)$ since the relation
``$n\!\in\! {\cal W}~\wedge ~(n,\alpha ,\beta )\!\notin\! {\cal C}$" is $\Ca$, so that
$(n,\alpha )\!\mapsto\!\lambda ({\cal C}_{n,\alpha})$ is $\Ca$-recursive on
${\cal W}\!\times\!\omega^\omega$.\bigskip

\noindent (b) Let $A\! :=\!\big\{ (\alpha ,\beta)\!\in\!\omega^\omega\!\times\! 2^\omega\mid
\big(\alpha (0),\alpha^*,\beta\big)\!\in\! O_0\big\}$. Note that $A$ is $\Ana$. By (a), the relation $R_A(\alpha ,p)\Leftrightarrow\lambda (A_\alpha )\! >\! r_p$ is $\Ana$. Therefore the relation 
$R_{O_0}(n,\alpha ,p)\Leftrightarrow R_A(n\alpha ,p)$ is $\Ana$ too. Moreover, 
$R_{O_0}(n,\alpha ,p)\Leftrightarrow\lambda\big( (O_0)_{n,\alpha}\big)\! >\! r_p\Leftrightarrow
\lambda_O(n,\alpha )\! >\! r_p$.\bigskip

\noindent $\bullet$ Assume now that $n\!\in\! D$. Then as above there is 
$\phi''\! :\!\omega^2\!\rightarrow\!\omega$ recursive such that
$$\begin{array}{ll}
\lambda_O(n,\alpha)\! >\! r_p\!\!\!
& \Leftrightarrow\lambda\big( (O_1)_{n,\alpha}\big)\! >\! r_p\Leftrightarrow
\lambda\big( (\neg O_1)_{n,\alpha}\big)\! <\! 1\! -\! r_p\cr
& \Leftrightarrow\exists l\!\in\!\omega ~~
\lambda\big( (\neg O_1)_{n,\alpha}\big)\!\leq\! 1\! -\! r_p\! -\!\frac{1}{l+1}\Leftrightarrow
\exists l\!\in\!\omega ~~\neg\Big(\lambda\big( (\neg O_1)_{n,\alpha}\big)\! >\! r_{\phi''(p,l)}\Big)\cr
& \Leftrightarrow\exists l\!\in\!\omega ~~\neg R_{\neg O_1}\big( n,\alpha ,\phi''(p,l)\big)
\mbox{,}
\end{array}$$
which shows the existence of $R'_{O_0}\!\in\!\Ca$ such that 
$\lambda_O(n,\alpha)\! >\! r_p\Leftrightarrow R'_{O_0}(n,\alpha ,p)$ if $n\!\in\! D$.\bigskip

\noindent $\bullet$ Assume that $n\!\in\! D$. Then there is 
$\phi'\! :\!\omega\!\rightarrow\!\omega$ recursive such that 
$$\lambda_O(n,\alpha)\! <\! r_q
\Leftrightarrow\lambda\big( (O_1)_{n,\alpha}\big)\! <\! r_q
\Leftrightarrow\lambda\big( (\neg O_1)_{n,\alpha}\big)\! >\! 1\! -\! r_q
\Leftrightarrow R_{\neg O_1}\big( n,\alpha ,\phi'(q)\big)\mbox{,}$$
which shows the existence of $R''_{O_0}\!\in\!\Ana$ such that 
$\lambda_O(n,\alpha)\! <\! r_q\Leftrightarrow R''_{O_0}(n,\alpha ,q)$ if $n\!\in\! D$.

\vfill\eject

\noindent $\bullet$ Assume that $n\!\in\! D$. Then there is 
$\phi''\! :\!\omega^2\!\rightarrow\!\omega$ recursive such that 
$$\begin{array}{ll}
\lambda_O(n,\alpha)\! <\! r_q\!\!\!
& \Leftrightarrow\lambda\big( (O_0)_{n,\alpha}\big)\! <\! r_q\Leftrightarrow
\exists l\!\in\!\omega ~~
\lambda\big( (O_0)_{n,\alpha}\big)\!\leq\! 1\! -\! r_q\! -\!\frac{1}{l+1}\cr
& \Leftrightarrow
\exists l\!\in\!\omega ~~\neg\Big(\lambda\big( (O_0)_{n,\alpha}\big)\! >\! r_{\phi''(q,l)}\Big)\Leftrightarrow\exists l\!\in\!\omega ~~\neg R_{O_0}\big( n,\alpha ,\phi''(q,l)\big)
\mbox{,}
\end{array}$$
which shows the existence of $R'''_{O_0}\!\in\!\Ca$ such that 
$\lambda_O(n,\alpha)\! <\! r_q\Leftrightarrow R'''_{O_0}(n,\alpha ,q)$ if $n\!\in\! D$.\bigskip

\noindent $\bullet$ Finally, 
$r_p\! <\!\lambda_O(n,\alpha )\! <\! r_q\Leftrightarrow R_{O_0}(n,\alpha ,p)\wedge 
R''_{O_0}(n,\alpha ,q)$ and 
$$r_p\! <\!\lambda_O(n,\alpha )\! <\! r_q\Leftrightarrow R'_{O_0}(n,\alpha ,p)\wedge 
R'''_{O_0}(n,\alpha ,q)$$ 
if $n\!\in\! D$, which shows that $\lambda_O$ is $\Ana$-recursive and $\Ca$-recursive on 
$D\!\times\!\omega$.\bigskip

\noindent (c) We first prove the following. Let $X,Y$ be a recursively presented Polish spaces and
${g\! :\! X\!\times\!\omega\!\rightarrow\! Y}$ be a $\Borel$-recursive map. Then the partial map
$h\! :\! X\!\rightarrow\! Y$ defined by
$$h(x)\! :=\!\mbox{lim}_{l\rightarrow\infty}~g(x,l)$$
when this limit exists is $\Borel$-recursive.\bigskip

 Indeed, the domain $D$ of $h$ is $\{ x\!\in\! X\mid
\forall r\!\in\!\omega ~~\exists L\!\in\!\omega ~~\forall k,l\!\geq\! L~~d_Y\big( g(x,k),g(x,l)\big)\! <\! 2^{-r}\}$, so that $D$ is $\Borel$. If $x\!\in\! D$, then $h(x)\!\in\! N(Y,u)$ is equivalent to
$$\exists p,q\!\in\!\omega ~~\frac{p}{q+1}\! <\!\frac{\big( (u)_1\big)_1}{\big( (u)_1\big)_2+1}~\wedge ~
\exists L\!\in\!\omega ~~\forall l\!\geq\! L~~
g(x,l)\!\in\! N\big( Y,\big<~0,<\big( (u)_1\big)_0,p,q>\!\big>\big)\mbox{,}$$
and we are done.\bigskip

\noindent $\bullet$ We set $B'\! :=\!\big\{ (\alpha ,\gamma )\!\in\!\omega^\omega\!\times\! 2^\omega\mid
\big( (\alpha )_0,\gamma\big)\!\in\! B~\wedge ~\gamma\!\in\!
N_{(\alpha )_1^*\vert (\alpha )_1(0)}\big\}$, so that
$B_\alpha\cap N_{\beta\vert l}\! =\! B'_{<\alpha ,l\beta >}$ and $B'$ is $\Borel$. By (a), the map
$g\! :\!\omega^\omega\!\times\! 2^\omega\!\times\!\omega\!\rightarrow\! [0,1]$ defined by
$g(\alpha ,\beta ,l)\! :=\! 2^{-l}\lambda(B_\alpha\cap N_{\beta\vert l})$ is $\Borel$-recursive. By the previous point, the partial map $h\! :\!\omega^\omega\!\times\! 2^\omega\!\rightarrow\! [0,1]$ defined by
$$h(\alpha ,\beta )\! :=\!\mbox{lim}_{l\rightarrow\infty}~2^{-l}\lambda(B_\alpha\cap N_{\beta\vert l})$$ when it exists is also $\Borel$-recursive. But $h\! =\! d_B$.\bigskip

\noindent $\bullet$ Fix $n\!\in\! {\cal W}$. Then there is $q(n)\!\in\! {\cal W}$ such that
$${\cal C}_{q(n)}\! =\!\big\{ (\gamma ,\delta )\!\in\!\omega^\omega\!\times\! 2^\omega\mid
\big( n,(\gamma )_0,\delta )\!\in\! {\cal C}~\wedge ~
(\gamma )_1^*\vert (\gamma )_1(0)\!\subseteq\!\delta\big\} .$$
Moreover, we may assume that $q$ is $\Ca$-recursive on $\cal W$, by the uniformization lemma. As $\Ca$ has the substitution property, the map
$g'\! :\! (n,\alpha ,\beta ,l)\!\mapsto\! 2^{-l}\lambda ({\cal C}_{q(n),<\alpha ,l\beta >})\! =\!
2^{-l}\lambda ({\cal C}_{n,\alpha}\cap N_{\beta\vert l})$ is $\Ca$-recursive on
${\cal W}\!\times\!\omega^\omega\!\times\! 2^\omega\!\times\!\omega$. As above, the map
$$h'\! :\! (n,\alpha ,\beta )\!\mapsto\!
\mbox{lim}_{l\rightarrow\infty}~2^{-l}\lambda ({\cal C}_{n,\alpha}\cap N_{\beta\vert l})\! =\!
d({\cal C}_{n,\alpha},\beta )$$
is $\Ca$-recursive on the $\Ca$ set
$\{ (n,\alpha ,\beta )\!\in\! {\cal W}\!\times\!\omega^\omega\!\times\! 2^\omega\mid
d({\cal C}_{n,\alpha},\beta )\mbox{ exists}\}$.\bigskip

\noindent (d) The argument here is partly similar to 11.6 and 17.25 in [K2]. We set, for
$(k,l)\!\in\!\omega^2$,
$$A_{k,l}\! :=\! h^{-1}\big( [\frac{k}{2^l},\frac{k\! +\! 1}{2^l})\big)$$
and define $h_l\! :\!\omega^\omega\!\times\! 2^\omega\!\rightarrow\! [0,1]$ by
$h_l\! =\!\Sigma_{k\leq 2^l}~\frac{k}{2^l}\chi_{A_{k,l}}$. We also define
$R\!\subseteq\!\omega^\omega\!\times\! 2^\omega\!\times\!\omega^3$ by
$$R(\alpha ,\beta ,u,k,l)\Leftrightarrow\frac{k}{2^l}\!\leq\! h(\alpha ,\beta )\! <\!\frac{k\! +\! 1}{2^l}
~\wedge ~\mbox{Seq}(u)~\wedge ~\beta\!\in\! N^u\mbox{,}$$
so that $R$ is $\Borel$. Then we define $O\!\subseteq\!\omega^\omega\!\times\! 2^\omega$ by
$$O(\alpha ,\beta )\Leftrightarrow\mbox{Seq}\big(\alpha (0)\big) ~\wedge ~
\mbox{lh}\big(\alpha (0)\big)\! =\! 3~\wedge ~
R\Big(\alpha^*,\beta ,\big(\alpha (0)\big)_0,\big(\alpha (0)\big)_1,\big(\alpha (0)\big)_2\Big)
\mbox{,}$$
so that $O$ is $\Borel$.\bigskip

 Note that $(h_l)$ is a sequence of Borel functions pointwise converging to $h$. By Lebesgue's dominated convergence theorem, $\int_{N^u}h(\alpha ,.)~d\lambda\! =\!
\mbox{lim}_{l\rightarrow\infty}~\int_{N^u}h_l(\alpha ,.)~d\lambda$ if $\mbox{Seq}(u)$. Note that
$$\begin{array}{ll}
\int_{N^u}h_l(\alpha ,.)~d\lambda\!\!\!
& \! =\!\int_{N^u}\Sigma_{k\leq 2^l}~\frac{k}{2^l}\chi_{A_{k,l}}(\alpha ,.)~d\lambda\! =\!
\Sigma_{k\leq 2^l}~\frac{k}{2^l}\lambda\big( (A_{k,l})_\alpha\cap N^u\big)\cr\cr
& \! =\!\Sigma_{k\leq 2^l}~\frac{k}{2^l}\lambda (R_{\alpha ,u,k,l})
\! =\!\Sigma_{k\leq 2^l}~\frac{k}{2^l}\lambda (O_{<u,k,l>\alpha}).
\end{array}$$
Using (a), this implies that the map $(\alpha ,u,l)\!\mapsto\!\int_{N^u}h_l(\alpha ,.)~d\lambda$ is
$\Borel$-recursive on its $\Borone$ domain
$\omega^\omega\!\times\!\{ u\!\in\!\omega\mid\mbox{Seq}(u)\}\!\times\!\omega$. As in the proof of (c), $i_h$ is $\Borel$-recursive on its domain.\hfill{$\square$}\bigskip

 We now prove a uniform version of Theorem 4.3.2 in [K1] (due to Tanaka, see [T2]).

\begin{thm} \label{unifreg} Let $B\!\in\!\Borel (\omega^\omega\!\times\! 2^\omega )$, and
$\epsilon\! :\!\omega^\omega\!\rightarrow\!\mathbb{R}$ be $\Borel$-recursive such that
$\epsilon (\alpha )\!\in\! (0,1]$ for each $\alpha\!\in\!\omega^\omega$. Then there is
$T\!\in\!\Borel (\omega^\omega\!\times\!\omega )$ such that\smallskip

\noindent (a) $T_\alpha$ is a tree for each $\alpha\!\in\!\omega^\omega$,\smallskip

\noindent (b) if $K\! =\!\big\{ (\alpha ,\beta )\!\in\!\omega^\omega\!\times\! 2^\omega\mid
\forall l\!\in\!\omega ~~\big(\alpha ,\overline{\beta}(l)\big)\!\in\! T\big\}$, then
$K_\alpha\!\subseteq\! B_\alpha$ and
$\lambda (K_\alpha)\!\geq\!\lambda (B_\alpha )\! -\!\epsilon (\alpha )$ for each
$\alpha\!\in\!\omega^\omega$.\end{thm}

\noindent\bf Proof.\rm\ Theorem \ref{repboreff} gives
$\pi\! :\!\omega^\omega\!\rightarrow\!\omega^\omega\!\times\! 2^\omega$ recursive and
$C\!\in\!\Bormone (\omega^\omega )$ such that $\pi$ is injective on $C$ and $\pi [C]\! =\! B$. We set
$Q\! :=\!\{ (\alpha ,\beta ,\gamma )\!\in\! (\omega ^\omega )^3\mid\gamma\!\in\! C\wedge
\pi (\gamma )\! =\! (\alpha ,\beta )\}$. As $Q\!\in\!\Bormone$, Theorem \ref{repsemi} gives a recursive subset $\overline{T}$ of $\omega^\omega\!\times\!\omega^2$ such that
$(\alpha ,\beta ,\gamma )\!\in\! Q\Leftrightarrow\forall l\!\in\!\omega ~~
\big(\alpha ,\overline{\beta}(l),\overline{\gamma}(l)\big)\!\in\!\overline{T}$ and $\overline{T}_\alpha$ is a tree for each $\alpha\!\in\!\omega^\omega$.\bigskip

\noindent $\bullet$ We set, for $u,v\!\in\!\omega$,
$$u\leq^av~\Leftrightarrow ~\mbox{Seq}(u),\mbox{Seq}(v)~\wedge ~\mbox{lh}(u)\! =\!\mbox{lh}(v)
~\wedge ~\forall i\! <\!\mbox{lh}(u)~~(u)_i\!\leq\! (v)_i.$$
$\bullet$ Then we set, for $u\!\in\!\omega$ with $\mbox{Seq}(u)$ and
$\alpha\!\in\!\omega^\omega$,
$$B^u_\alpha\! :=\!\big\{\beta\!\in\! 2^\omega\mid\exists\gamma\!\in\!\omega^\omega ~~
\overline{\gamma}\big(\mbox{lh}(u)\big)\leq^au~\wedge ~\forall l\!\in\!\omega ~~
\big(\alpha ,\overline{\beta}(l),\overline{\gamma}(l)\big)\!\in\!\overline{T}\big\}$$
and $B'\! :=\!\{ (\alpha ,\beta )\!\in\!\omega^\omega\!\times\! 2^\omega\mid
\mbox{Seq}\big(\alpha (0)\big)~\wedge ~\beta\!\in\! B^{\alpha (0)}_{\alpha^*}\}$. Note that $B'$ is
$\Ana$. In fact, $B'$ is $\Borel$ by uniqueness of the witness $\gamma$.

\vfill\eject

\noindent $\bullet$ We now define $\delta_\alpha\!\in\!\omega^\omega$ as follows. We define
$\delta_\alpha (i)$ by induction on $i$. We first set
$$\delta_\alpha (0)\! :=\!\mbox{min}\{ k\!\in\!\omega\mid
\lambda (B_\alpha^{<k>})\! >\!\lambda (B_\alpha )\! -\!\frac{\epsilon (\alpha )}{2}\} .$$
This number exists since $B_\alpha$ is the increasing union of the $B_\alpha^{<k>}$'s. Then
$$\delta_\alpha (i\! +\! 1)\! :=\!\mbox{min}\{ k\!\in\!\omega\mid
\lambda (B_\alpha^{<\delta_\alpha (0),...,\delta_\alpha (i),k>})\! >\!\lambda (B_\alpha )\! -\!
\frac{\epsilon (\alpha )}{2}\! -\! ...\! -\!\frac{\epsilon (\alpha )}{2^{i+2}}\} .$$
Note that $\delta_\alpha\!\in\!\Borel (\alpha )$, by Corollary \ref{Delta}.(a).\bigskip

\noindent $\bullet$ We set $T\! :=\!\{ (\alpha ,v)\!\in\!\omega^\omega\!\times\!\omega\mid
\mbox{Seq}(v)\wedge\exists u\leq^a\overline{\delta_\alpha}(\mbox{lh}(v))~~
(\alpha ,v,u)\!\in\!\overline{T}\}$, so that $T\!\in\!\Borel (\omega^\omega\!\times\!\omega )$ and 
$T_\alpha$ is a tree for each $\alpha\!\in\!\omega^\omega$.\bigskip

\noindent $\bullet$ We set $K\! :=\!\{ (\alpha ,\beta )\!\in\!\omega^\omega\!\times\! 2^\omega\mid
\forall l\!\in\!\omega ~~\beta\!\in\! B^{\overline{\delta_\alpha}(l)}_\alpha\}$, so that
$K_\alpha\!\subseteq\! B_\alpha$ and
$$\lambda (K_\alpha )\! =\!\mbox{lim}_{l\rightarrow\infty}~
\lambda (B^{\overline{\delta_\alpha}(l)}_\alpha )\!\geq\!\lambda (B_\alpha )\! -\!\epsilon (\alpha )$$
for each $\alpha\!\in\!\omega^\omega$ since
$(B^{\overline{\delta_\alpha}(l)}_\alpha)_{l\in\omega}$ is decreasing. It remains to apply K\"onig's lemma to see that $K\! =\!\big\{ (\alpha ,\beta )\!\in\!\omega^\omega\!\times\! 2^\omega\mid
\forall l\!\in\!\omega ~~\big(\alpha ,\overline{\beta}(l)\big)\!\in\! T\big\}$ since
$$\big\{ s\!\in\!\omega^{<\omega}\mid
<s(0),...,s(\vert s\vert\! -\! 1)>\leq^a\overline{\delta_\alpha}\big(\vert s\vert\big)~\wedge ~
\big(\alpha ,\overline{\beta }(\vert s\vert ),<s(0),...,s(\vert s\vert\! -\! 1)>\!\!\big)\!\in\!\overline{T}\big\}$$
is a finitely splitting tree.\hfill{$\square$}\bigskip

\noindent - We want to prove an effective and uniform version of the Lusin-Menchoff lemma. We first need the following result, which slightly and uniformly refines Theorem A in [L] at the first level of the Borel hierarchy.

\begin{lem} \label{opensect} Let $O$ be a $\Borel$ subset of $\omega^\omega\!\times\! 2^\omega$ with open vertical sections. Then there is a $\Borel$-recursive map $f\! :\!\omega^\omega\!\rightarrow\!\omega^\omega$ such that $O_\alpha$ is the disjoint union 
$\bigcup\big\{ N^{f(\alpha )(u)}\mid u\!\in\!\omega\wedge\mbox{Seq}\big( f(\alpha )(u)\big)\big\}$, for each $\alpha\!\in\!\omega^\omega$.\end{lem}

\noindent\bf Proof.\rm\ Let $P\! :=\!\big\{ (\alpha ,u)\!\in\!\omega^\omega\!\times\!\omega\mid
\mbox{Seq}(u)~\wedge ~\big( \mbox{lh}(u)\! =\! 0~\vee ~(N^u\!\subseteq\! O_\alpha ~\wedge ~
N^{u^-}\!\not\subseteq\! O_\alpha )\big)\big\}$. Note that $P$ is $\Ca$, since a nonempty
$\Borel (\alpha )$ closed subset of $2^\omega$ contains a $\Borel (\alpha )$ point, by 4F.15 in [M]. We then define a relation $R$ on $\omega^\omega\!\times\! 2^\omega\!\times\!\omega$ by
$R(\alpha ,\beta ,u)\Leftrightarrow P(\alpha ,u)~\wedge ~\beta\!\in\! N^u$, so that $R$ is $\Ca$. Note that, for each $(\alpha ,\beta )\!\in\! O$ there is $u$ with $R(\alpha ,\beta ,u)$. By 4B.5 in [M], there is a $\Borel$-recursive map $g\! :\!\omega^\omega\!\times\! 2^\omega\!\rightarrow\!\omega$ such that $R\big(\alpha ,\beta ,g(\alpha ,\beta )\big)$ for each $(\alpha ,\beta )\!\in\! O$. Fix
$\alpha\!\in\!\omega^\omega$. Note that $S^\alpha\! :=\!\{ g(\alpha ,\beta )\mid\beta\!\in\! O_\alpha\}$ is a $\Ana (\alpha )$ subset of $\omega$ contained in the $\Ca (\alpha )$ set $P_\alpha$. By 
4B.11 and 4C in [M], there is $D^\alpha\!\in\!\Borel (\alpha )$ with
$S^\alpha\!\subseteq\! D^\alpha\!\subseteq\! P_\alpha$. Note that
$O_\alpha\!\subseteq\!\bigcup_{u\in D^\alpha}~N^u\!\subseteq\! O_\alpha$, so that
$O_\alpha$ is the disjoint union of the sequence $(N^u)_{u\in D^\alpha}$. We define
$\delta_\alpha\!\in\!\omega^\omega$ by
$$\delta_\alpha (u)\! :=\!\left\{\!\!\!\!\!\!\!
\begin{array}{ll}
& u\mbox{~if }u\!\in\! D_\alpha\mbox{,}\cr
& 0\mbox{~otherwise.}
\end{array}
\right.$$
Note that $\delta_\alpha\!\in\!\Borel (\alpha )$ and $O_\alpha$ is the disjoint union 
$\bigcup\big\{ N^{\delta_\alpha (u)}\mid u\!\in\!\omega\wedge
\mbox{Seq}\big(\delta_\alpha (u)\big)\big\}$. As the set 
$\Big\{ (\alpha ,\delta )\!\in\!\omega^\omega\!\times\!\omega^\omega\mid
\delta\!\in\!\Borel (\alpha )\wedge O_\alpha\mbox{ is the disjoint union }\bigcup
\big\{ N^{\delta (u)}\mid u\!\in\!\omega\wedge
\mbox{Seq}\big(\delta (u)\big)\big\}\Big\}$ 
is $\Ca$, it remains to apply the uniformization lemma to get the desired map $f$.
\hfill{$\square$}

\vfill\eject

\noindent\bf Notation.\rm\ We set ${\cal W}_1\! :=\!\{ n\!\in\! {\cal W}\mid
\forall\alpha\!\in\!\omega^\omega ~~\exists\gamma_n\!\in\!\Borel (\alpha )~~
{\cal C}_{n,\alpha}\! =\!\bigcup\big\{ N^{\gamma_n(u)}\mid u\!\in\!\omega\wedge
\mbox{Seq}\big(\gamma_n(u)\big)\big\}$, so that, by Lemma \ref{opensect}, ${\cal W}_1$ is a $\Ca$ set of codes for the $\Borel$ subsets of $\omega^\omega\!\times\! 2^\omega$ with open vertical sections.

\begin{lem} \label{LMeff} Let $F$ be a $\Borel$ subset of $\omega^\omega\!\times\! 2^\omega$ with closed vertical sections, and $B$ be a $\Borel$ subset of $\omega^\omega\!\times\! 2^\omega$ such that $B\!\supseteq\! F$ and $d(B_\alpha ,\beta )\! =\! 1$ for each $(\alpha ,\beta )\!\in\! F$.  Then there is a $\Borel$ subset $C$ of $\omega^\omega\!\times\! 2^\omega$ with closed vertical sections such that\smallskip

(1) $F\!\subseteq\! C\!\subseteq\! B$,\smallskip

(2) $d(B_\alpha ,\beta )\! =\! 1$ for each $(\alpha ,\beta )\!\in\! C$,\smallskip

(3) $d(C_\alpha ,\beta )\! =\! 1$ for each $(\alpha ,\beta )\!\in\! F$.\end{lem}

\noindent\bf Proof.\rm\ Lemma \ref{opensect} gives a $\Borel$-recursive map
$f\! :\!\omega^\omega\!\rightarrow\!\omega^\omega$ such that $(\neg F)_\alpha$ is the disjoint union $\bigcup\big\{ N^{f(\alpha )(u)}\mid u\!\in\!\omega\wedge\mbox{Seq}\big( f(\alpha )(u)\big)\big\}$, for each $\alpha\!\in\!\omega^\omega$.  We set
$$B'\! :=\!\Big\{ (\alpha ,\gamma )\!\in\!\omega^\omega\!\times\! 2^\omega\mid
\big( (\alpha )_0,\gamma\big)\!\in\! B~\wedge ~
\mbox{Seq}\Big( f\big( (\alpha )_0\big)\big( (\alpha )_1(0)\big)\Big)~\wedge ~
\gamma\!\in\! N^{f((\alpha )_0)((\alpha )_1(0))}\Big\}\mbox{,}$$ 
so that $B'$ is $\Borel$ and $B_\alpha\cap N^{f(\alpha )(u)}\! =\! B'_{<\alpha ,u^\infty >}$ if $\mbox{Seq}\big( f(\alpha )(u)\big)$. By Corollary \ref{Delta}.(c), the partial map
$(\alpha ,\beta ,u)\!\mapsto\! d(B_\alpha\cap N^{f(\alpha )(u)},\beta )$ is $\Borel$-recursive. We then set
$$B''\! :=\!\{ (\alpha ,\gamma )\!\in\! B'\mid
d(B_{(\alpha )_0}\cap N^{f((\alpha )_0)((\alpha )_1(0))},\gamma )\! =\! 1\}\mbox{,}$$
so that $B''$ is $\Borel$ and $\{\beta\!\in\! B_\alpha\cap N^{f(\alpha )(u)}\mid d(B_\alpha\cap N^{f(\alpha )(u)},\beta )\! =\! 1\}
\! =\! B''_{<\alpha ,u^\infty >}$ if $\mbox{Seq}\big( f(\alpha )(u)\big)$. We define 
$\epsilon\! :\!\omega^\omega\!\rightarrow\!\mathbb{R}$ by
$$\varepsilon (\alpha )\! :=\!\left\{\!\!\!\!\!\!
\begin{array}{ll}
& 2^{-(\alpha )_1(0)}\lambda (B'_\alpha )\mbox{ if }\lambda (B'_\alpha )\!\not=\! 0\mbox{,}\cr
& 1\mbox{ otherwise,}
\end{array}
\right.$$
so that $\epsilon$ is $\Borel$-recursive by Corollary \ref{Delta}.(a), and
$\epsilon (\alpha )\!\in\! (0,1]$ for each $\alpha\!\in\!\omega^\omega$. Theorem \ref{unifreg} gives $T\!\in\!\Borel (\omega^\omega\!\times\!\omega )$ such that\bigskip

(a) $T_\alpha$ is a tree for each $\alpha\!\in\!\omega^\omega$,\bigskip

(b) if $K\! =\!\big\{ (\alpha ,\beta )\!\in\!\omega^\omega\!\times\! 2^\omega\mid
\forall l\!\in\!\omega ~~\big(\alpha ,\overline{\beta}(l)\big)\!\in\! T\big\}$, then
$K_\alpha\!\subseteq\! B''_\alpha$ and
$\lambda (K_\alpha)\!\geq\!\lambda (B''_\alpha )\! -\!\epsilon (\alpha )$ for each
$\alpha\!\in\!\omega^\omega$.\bigskip

 We set, for $u\!\in\!\omega$,
$$F^u\! :=\!\big\{ (\alpha ,\beta )\!\in\!\omega^\omega\!\times\! 2^\omega\mid
\mbox{Seq}\big( f(\alpha )(u)\big) ~\wedge ~(<\alpha ,u^\infty >,\beta )\!\in\! K~\wedge ~
\lambda (B'_{<\alpha ,u^\infty >})\!\not=\! 0\big\} .$$
As $K$ is $\Borel$ with closed vertical sections, so is $F^u$. If $\mbox{Seq}\big( f(\alpha )(u)\big)$ and $\lambda (B'_{<\alpha ,u^\infty >})\! =\! 0$, then
$\lambda (B_\alpha\cap N^{f(\alpha )(u)})\! =\! 0$ and $F^u_\alpha\! =\!\emptyset$, so that
$F^u_\alpha\!\subseteq\!\{\beta\!\in\! B_\alpha\cap N^{f(\alpha )(u)}\mid
d(B_\alpha\cap N^{f(\alpha )(u)},\beta )\! =\! 1\}$ and
$\lambda (F^u_\alpha )\!\geq\! (1\! -\! 2^{-u})\lambda (B_\alpha\cap N^{f(\alpha )(u)})$. If
$\mbox{Seq}\big( f(\alpha )(u)\big)$ and $\lambda (B'_{<\alpha ,u^\infty >})\!\not=\! 0$, then
$$F^u_\alpha\! =\! K_{<\alpha ,u^\infty >}\!\subseteq\! B''_{<\alpha ,u^\infty >}\! =\!
\{\beta\!\in\! B_\alpha\cap N^{f(\alpha )(u)}\mid d(B_\alpha\cap N^{f(\alpha )(u)},\beta )\! =\! 1\} .$$

 Moreover,
$$\begin{array}{ll}
\lambda (F^u_\alpha )\!\!\!\!
& \! =\!\lambda (K_{<\alpha ,u^\infty >})\!\geq\!
\lambda (B''_{<\alpha ,u^\infty >})\! -\!\epsilon (<\alpha ,u^\infty >)\! =\!
\lambda (B''_{<\alpha ,u^\infty >})\! -\! 2^{-u}\lambda (B'_{<\alpha ,u^\infty >})\cr
& \! =\! (1\! -\! 2^{-u})\lambda (B_\alpha\cap N^{f(\alpha )(u)})
\end{array}$$
since $\lambda (B_\alpha\cap N^{f(\alpha )(u)})\! =\!
\lambda (\{\beta\!\in\! B_\alpha\cap N^{f(\alpha )(u)}
\mid d(B_\alpha\cap N^{f(\alpha )(u)},\beta )\! =\! 1\} )$, by Theorem \ref{densityone}. It remains to set $C\! :=\! F\cup\bigcup_{u\in\omega}~F^u$. We conclude as in the proof of Lemma \ref{LM}.
\hfill{$\square$}\bigskip

\noindent - We now want to prove an effective and uniform version of Lemma \ref{map}.

\begin{lem} \label{mapeff} Let $C$ be a $\Borel$ subset of $\omega^\omega\!\times\! 2^\omega$ with closed vertical sections, $\cal G$ be a Borel subset of $2^\omega$ with $\lambda ({\cal G})\! =\! 0$, and $G$ be a $\Borel$ subset of $\omega^\omega\!\times\! 2^\omega$ with $G_\delta$ vertical sections, contained in $\omega^\omega\!\times\! {\cal G}$ and disjoint from $C$. Then there is a $\Borel$-recursive map
$h\! :\!\omega^\omega\!\times\! 2^\omega\!\rightarrow\!\mathbb{R}$ such that
$h(\alpha ,\cdot )\! :\! 2^\omega\!\rightarrow\! [0,1]$ is $\tau$-continuous for each $\alpha\!\in\!\omega^\omega$, $h_{\vert C}\!\equiv\! 0$ and
$h_{\vert G}\!\equiv\! 1$.\end{lem}

\noindent\bf Proof.\rm\ By Theorem 3.5 in [L], there is a $\Borel$ subset $F$ of
$\omega\!\times\!\omega^\omega\!\times\! 2^\omega$ such that $F_{n,\alpha}$ is closed for each
$(n,\alpha )\!\in\!\omega\!\times\!\omega^\omega$ and $\neg G\! =\!\bigcup_{n\in\omega}~F_n$. Moreover, we may assume that $(F_n)_{n\in\omega}$ is increasing and $F_0\! =\! C$.\bigskip

\noindent $\bullet$ We will define, by primitive recursion, a partial map
$f\! :\!\omega\!\rightarrow\!\omega$ which is $\Ca$-recursive on its domain such that $f(n)$ essentially codes the set $C_{\frac{1}{2^n}}$ constructed in the proof of Lemma \ref{map}. As this map will in fact be total, it will be $\Borel$-recursive by the uniformization lemma.\bigskip

 We first apply Lemma \ref{LMeff} to $F\! :=\! F_0$ and $B\! :=\!\neg G$. This is possible because
$G_\alpha\!\subseteq\! {\cal G}$, so that $(\neg G)_\alpha$ has $\lambda$-measure one and therefore density one at any point of $2^\omega$, for each $\alpha\!\in\!\omega^\omega$. Lemma \ref{LMeff} gives $C_1\!\in\!\Borel$ with closed vertical sections such that
$\neg G\!\supseteq\! C_1\!\supseteq\! F_0$. Let $f(0)\!\in\! {\cal W}_1$ with
${\cal C}_{f(0)}\! =\!\neg C_1$.\bigskip

 More generally, we will have ${\cal C}_{f(n)}\! =\!\neg C_{\frac{1}{2^n}}$. As mentioned above, $f$ will be defined by primitive recursion, which means that there will be a partial map
$g\! :\!\omega^2\!\rightarrow\!\omega$ such that $f(n\! +\! 1)\! =\! g\big( f(n),n\big)$. This partial map $g$ will be $\Ca$-recursive on its $\Ca$ domain
$\{ m\!\in\! {\cal W}_1\mid\neg {\cal C}_m\!\subseteq\!\neg G\}\!\times\!\omega$, so that $f$ will be
$\Ca$-recursive on its domain by 7A.5 in [M]. The map $g$ will take values in ${\cal W}_1$, and is constructed in such a way that, if $A\! :=\!\neg {\cal C}_m\!\subseteq\!\neg G$ and
$A'\! :=\!\neg {\cal C}_{g(m,n)}$, then
$$\begin{array}{ll}
& (1)~A\cup F_{n+1}\!\subseteq\! A'\!\subseteq\!\neg G\mbox{,}\cr
& (2)~\forall (\alpha ,\beta )\!\in\! A'~~d\big( (\neg G)_\alpha ,\beta\big)\! =\! 1\mbox{,}\cr
& (3)~\forall (\alpha ,\beta )\!\in\! A\cup F_{n+1}~~d(A'_\alpha ,\beta )\! =\! 1.
\end{array}$$
Lemma \ref{LMeff} ensures that such a $g(m,n)\!\in\!\omega$ exists if
$(m,n)\!\in\! \{ q\!\in\! {\cal W}_1\mid\neg {\cal C}_q\!\subseteq\!\neg G\}\!\times\!\omega$. As the properties (1)-(3) are $\Ca$ by Corollary \ref{Delta}, the uniformization lemma ensures the existence of $g$. So we constructed a $\Borel$-recursive map $f\! :\!\omega\!\rightarrow\!\omega$, taking values in ${\cal W}_1$, such that $C_{\frac{1}{2^n}}\! :=\!\neg {\cal C}_{f(n)}$ is a $\Borel$ subset of $\omega^\omega\!\times\! 2^\omega$ with closed vertical sections,
$F_n\!\subseteq\! C_{\frac{1}{2^n}}\!\subseteq\!\neg G$,
$C_{\frac{1}{2^n}}\!\subseteq\! C_{\frac{1}{2^{n+1}}}$, and
$$d\big( (C_{\frac{1}{2^{n+1}}})_\alpha ,\beta )\! =\! 1$$
if $(\alpha ,\beta )\!\in\! C_{\frac{1}{2^n}}$.\bigskip

\noindent $\bullet$ Similarly, we construct a $\Borel$-recursive map
$\tilde F\! :\!\omega\!\rightarrow\!\omega$ satisfying the following properties, if
$$D\! :=\!\{ p\!\in\!\omega\mid\mbox{Seq}(p)~\wedge ~
\mbox{lh}(p)\! =\! 2~\wedge ~0\! <\! (p)_1\!\leq\! 2^{(p)_0}\} .$$
$$\begin{array}{ll}
& (a)~\tilde F(p)\!\in\! {\cal W}_1\mbox{ if }p\!\in\! D\mbox{, in which case we set }
C_p\! :=\!\neg {\cal C}_{\tilde F(p)}\mbox{,}\cr
& (b)~C_p\!\subseteq\! C_{p'}\mbox{ if }p,p'\!\in\! D~\wedge ~
\frac{(p')_1}{2^{(p')_0}}\!\leq\!\frac{(p)_1}{2^{(p)_0}}\mbox{,}\cr
& (c)~d\big( (C_{p'})_\alpha ,\beta\big)\! =\! 1\mbox{ if }p,p'\!\in\! D~\wedge ~
\frac{(p')_1}{2^{(p')_0}}\! <\!\frac{(p)_1}{2^{(p)_0}}~\wedge ~(\alpha ,\beta )\!\in\! C_p.
\end{array}$$
$\bullet$ This allows us to define $h$ by
$$1\! -\! h(\alpha ,\beta )\! :=\!\left\{\!\!\!\!\!\!
\begin{array}{ll}
& 0\mbox{ if }(\alpha ,\beta )\!\in\! G\mbox{,}\cr
& \mbox{sup}\{\frac{(p)_1}{2^{(p)_0}}\mid p\!\in\! D~\wedge ~(\alpha ,\beta )\!\in\! C_p\}
\mbox{ if }(\alpha ,\beta )\!\notin\! G.
\end{array}
\right.$$
Note that $h$ is $\Borel$-recursive since $D\!\in\!\Borone$, so that the relation
``$p\!\in\! D~\wedge ~(\alpha ,\beta )\!\in\! C_p$" is $\Borel$ in $(p,\alpha ,\beta )$. We conclude as in the proof of Lemma \ref{map}.\hfill{$\square$}\bigskip

\noindent - We are now ready to prove the main lemma in this section. We equip the space
$[0,1]^{2^{<\omega}}$ with the distance defined by $d(f,g)\! :=\!
\Sigma_{u\in\omega ,\mbox{Seq}(u)}~\frac{\vert f(s(u))-g(s(u))\vert}{2^{u+1}}$. We give a recursive presentation of $([0,1]^{2^{<\omega}},d)$. We set
$$f_n(s)\! :=\!\left\{\!\!\!\!\!\!
\begin{array}{ll}
& \frac{((n)_{\overline{s}})_0}{((n)_{\overline{s}})_0+((n)_{\overline{s}})_1+1}
\mbox{ if }\mbox{Seq}(n)~\wedge ~\forall k\! <\!\mbox{lh}(n) ~~
\big( ~\mbox{Seq}\big( (n)_k\big) ~\wedge ~\mbox{lh}\big( (n)_k\big)\! =\! 2~\big)~\wedge ~
\overline{s}\! <\!\mbox{lh}(n)\mbox{,}\cr
& 0\mbox{ otherwise,}
\end{array}
\right.$$
so that $(f_n)$ is dense in $[0,1]^{2^{<\omega}}$. It is now routine to check that the relations
``$d(f_m,f_n)\!\leq\!\frac{p}{q+1}$" and
``$d(f_m,f_n)\! <\!\frac{p}{q+1}$" are recursive in $(m,n,p,q)$. It is also routine to check that $F\! :\!\omega^\omega\!\rightarrow\! [0,1]^{2^{<\omega}}$ is $\Borel$-recursive if the map 
$F'\! :\!\omega\!\times\!\omega^\omega\!\rightarrow\!\mathbb{R}$ defined by
$F'(u,\alpha )\!\ :=\! F(\alpha )\big( s(u)\big)$ if $\mbox{Seq}(u)$, $0$ otherwise, is $\Borel$-recursive ($s(u)$ was defined at the beginning of Section 3).

\begin{lem} \label{Borel} Let
${\cal V}\! :=\!\{ (f,\beta )\!\in\! {\cal M}\!\times\! 2^\omega\mid\mbox{osc}(f,\beta )\! >\! 0\}$, $\cal G$ be a nonempty $G_\delta\cap\Borel$ subset of $2^\omega$ with $\lambda ({\cal G})\! =\! 0$, and $G$ be a $\Borel$ subset of $\omega^\omega\!\times\! 2^\omega$, contained in $\omega^\omega\!\times\! {\cal G}$, and with $G_\delta$ vertical sections. Then there is a $\Borel$-recursive map $F\! :\!\omega^\omega\!\rightarrow\! [0,1]^{2^{<\omega}}$, taking values in $\cal M$, and such that $G_\alpha\! =\! {\cal V}_{F(\alpha )}$ for each $\alpha\!\in\!\omega^\omega$.\end{lem}

\noindent\bf Proof.\rm\ We will define, by primitive recursion,
$f\! :\!\omega\!\rightarrow\!\omega^4$ coding $g_n$, $S_n$, $G^*_n$, and $(s^n_j)_{j\in I_n}$ defining $G^{**}_n$ considered in the proof of the Lemma \ref{Gdelta}. We must find
$r\! :\!\omega^4\!\times\!\omega\!\rightarrow\!\omega^4$ with $f(n\! +\! 1)\! =\! r\big( f(n),n\big)$. In practice,\bigskip

(1) $f_0(n)\!\in\! {\cal W}_1$ codes $G^*_n\!\subseteq\!\omega^\omega\!\times\! 2^\omega$,\smallskip

(2) $f_1(n)\!\in\! {\cal W}^{\omega^\omega\times 2^\omega\times\mathbb{R}}$ codes the graph of $g_n\! :\!\omega^\omega\!\times\! 2^\omega\!\rightarrow\!\mathbb{R}$,\smallskip

(3) $f_2(n)\!\in\! {\cal W}^{\omega^\omega\times 2^\omega\times\mathbb{R}}$ codes the graph of $S_n\! :\!\omega^\omega\!\times\! 2^\omega\!\rightarrow\!\mathbb{R}$,\smallskip

(4) $f_3(n)\!\in\! {\cal W}^{\omega^\omega\times\omega^\omega}$ codes the graph of the function
$\alpha\!\mapsto\! (s^{n,\alpha}_j)_{j\in I_{n,\alpha}}$.\bigskip

\noindent $\bullet$ By Theorem 3.5 in [L], there is a $\Borel$ subset $O$ of
$\omega\!\times\!\omega^\omega\!\times\! 2^\omega$ such that $O_{n,\alpha}$ is open for each
$(n,\alpha )\!\in\!\omega\!\times\!\omega^\omega$ and $G\! =\!\bigcap_{n\in\omega}~O_n$. Moreover, we may assume that $(O_n)_{n\in\omega}$ is decreasing and
$O_0\! =\!\omega^\omega\!\times\! 2^\omega$.\bigskip

\noindent $\bullet$ Let $n_0\!\in\! {\cal W}_1$ with ${\cal C}_{n_0}\! =\!\omega^\omega\!\times\! 2^\omega$,
$n_1\!\in\! {\cal W}^{\omega^\omega\times 2^\omega\times\mathbb{R}}$ with
${\cal C}^{\omega^\omega\times 2^\omega\times\mathbb{R}}_{n_1}\! =\!
\{ (\alpha ,\beta ,r)\!\in\!\omega^\omega\!\times\! 2^\omega\!\times\!\mathbb{R}\mid r\! =\! 1\}$, and
$n_3\!\in\! {\cal W}^{\omega^\omega\times\omega^\omega}$ with
${\cal C}^{\omega^\omega\times\omega^\omega}_{n_3}\! =\!
\{ (\alpha ,\gamma )\!\in\!\omega^\omega\!\times\!\omega^\omega\mid\gamma\! =\! 10^\infty\}$. We set $f(0)\! :=\! (n_0,n_1,n_1,n_3)$, so that ${\cal C}_{n_0}\! =\! G^*_0$,
${\cal C}^{\omega^\omega\times 2^\omega\times\mathbb{R}}_{n_1}\! =\!\mbox{Gr}(g_0)
\! =\!\mbox{Gr}(S_0)$, ${\cal C}^{\omega^\omega\times\omega^\omega}_{n_3}\! =\!
\mbox{Gr}(\alpha\!\mapsto\! 10^\infty )$,
$$\big\{ u\!\in\!\omega\mid\mbox{Seq}\big( (10^\infty )(u)\big)\big\}\! =\!\{ 0\}\! =\! I_0$$
and $(10^\infty )(0)\! =\! 1\! =<>=\! s^0_0$. So $f(0)$ is as desired.\bigskip

\noindent $\bullet$ We now study the induction step. This means that we must define
$r(n_0,n_1,n_2,n_3,n)\!\in\!\omega^4$.\bigskip

\noindent (1) We first define $r_0(n_0,n_1,n_2,n_3,n)$ coding $G^*_{n+1}$. Fix
$n_3\!\in\! {\cal W}^{\omega^\omega\times\omega^\omega}$ coding the graph of a
$\Borel$-recursive function $\phi\! :\!\omega^\omega\!\rightarrow\!\omega^\omega$ such that the sequences $s\big(\phi (\alpha )(u)\big)$ coded by the $u$'s with $\mbox{Seq}\big(\phi (\alpha )(u)\big)$ are pairwise incompatible and $G_\alpha\!\subseteq\!\bigcup\big\{ N^{\phi (\alpha )(u)}\mid u\!\in\!\omega\wedge
\mbox{Seq}\big(\phi (\alpha )(u)\big)\big\}$ (we call $P_3$ the $\Ca$ set of such $n_3$'s). Let $\alpha\!\in\!\omega^\omega$. Assume that $\mbox{Seq}\big(\phi (\alpha )(u)\big)$ (which intuitively means that $u\!\in\! I_{n,\alpha}$ and $s^{n,\alpha}_u$ is coded by $\phi (\alpha )(u)$). By continuity of $\lambda$, 
$$0\! =\!\lambda (G_\alpha\cap N^{\phi (\alpha )(u)})\! =\!
\mbox{lim}_{j\rightarrow\infty}~\lambda (O_{j,\alpha }\cap N^{\phi (\alpha )(u)}).$$
This gives $j(n,\alpha ,u)\! >\! n$ minimal with $\lambda (O_{j(n,\alpha ,u),\alpha }\cap
N^{\phi (\alpha )(u)})\! <\! 2^{-n-3-\mbox{lh}(\phi (\alpha )(u))}$ (note that 
$2^{-\mbox{lh}(\phi (\alpha )(u))}\! =\!\lambda (N^{\phi (\alpha )(u)})$). Moreover, 
$G_\alpha\cap N^{\phi (\alpha )(u)}\!\subseteq\!
O_{j(n,\alpha ,u),\alpha }\cap N^{\phi (\alpha )(u)}\!\subseteq\!
O_{n+1,\alpha }\cap N^{\phi (\alpha )(u)}$, so that
$O_{j(n,\alpha ,u),\alpha }\cap N^{\phi (\alpha )(u)}$ satisfies the properties of the set $O_j$ in the proof of Lemma \ref{Gdelta}. We will have
$G^*_{n+1,\alpha}\! =\!\bigcup_{\mbox{Seq}(\phi (\alpha )(u))}~O_{j(n,\alpha ,u),\alpha}
\cap N^{\phi (\alpha )(u)}$. By Corollary \ref{Delta} and the uniformization lemma, we may assume that the map $j$ is $\Borel$-recursive on its $\Borel$ domain
$$\big\{ (n,\alpha ,u)\!\in\!\omega\!\times\!\omega^\omega\!\times\!\omega\mid
\mbox{ Seq}\big(\phi (\alpha )(u)\big)\big\} .$$ 
Note that $G^*_{n+1}$ is a $\Borel$ subset of $\omega^\omega\!\times\! 2^\omega$ with open vertical sections, which gives $m\!\in\! {\cal W}_1$ such that ${\cal C}_m\! =\! G^*_{n+1}$. By incompatibility,
$G^*_{n+1,\alpha}\cap N^{\phi (\alpha )(u)}\! =\! 
O_{j(n,\alpha ,u),\alpha}\cap N^{\phi (\alpha )(u)}$. So we proved that, for each 
$(n_3,n)\!\in\! P_3\!\times\!\omega$, there is $m\!\in\! {\cal W}_1$ such that, for each 
$\alpha\!\in\!\omega^\omega$, 
$$\begin{array}{ll}
& (1)~G_\alpha\!\subseteq\! {\cal C}_{m,\alpha}\!\subseteq\! O_{n+1,\alpha }\cap
\bigcup\big\{ N^{\phi (\alpha )(u)}\mid u\!\in\!\omega\wedge
\mbox{Seq}\big(\phi (\alpha )(u)\big)\big\}\mbox{,}\cr
& (5)~\lambda ({\cal C}_{m,\alpha}\cap N^{\phi (\alpha )(u)})\! <\! 
2^{-n-3-\mbox{lh}(\phi (\alpha )(u))}\mbox{ if }u\!\in\!\omega\wedge
\mbox{Seq}\big(\phi (\alpha )(u)\big).
\end{array}$$
By Corollary \ref{Delta} and the uniformization lemma, we may assume that the map
$\tilde r_0\! :\! (n_3,n)\!\mapsto\! m$ is $\Ca$-recursive on $P_3\!\times\!\omega$. We set
$r_0(n_0,n_1,n_2,n_3,n)\! :=\!\tilde r_0(n_3,n)$, which defines a partial map $r_0$ which is $\Ca$-recursive on its $\Ca$ domain $\omega^3\!\times\! P_3\!\times\!\omega$.\bigskip

\noindent (2) We now define $r_1(n_0,n_1,n_2,n_3,n)$ coding $g_{n+1}$. We use Lemma \ref{mapeff} and its proof. Note that
$r_0(n_0,n_1,n_2,n_3,n)\!\in\! D_0\! :=\!\{ m\!\in\! {\cal W}_1\mid G\!\subseteq\! {\cal C}_m\}$. The proof of Lemma \ref{mapeff} shows that for any $m\!\in\! D_0$ there is
$\tilde F_m\!\in\!\omega^\omega\cap\Borel$ satisfying the conditions (a), (b), (c) and
$$(d)~\forall p\!\in\! D~~\neg (0\! <\! (p)_1\! =\! 2^{(p)_0})~\vee ~
{\cal C}_{\tilde F_m(p)}\!\subseteq\! {\cal C}_m.$$
The uniformization lemma shows that we may assume that the partial map
$\tilde F\! :\! m\!\mapsto\!\tilde F_m$ is $\Ca$-recursive on $D_0$.

\vfill\eject

 The definition of $h$ in the proof of Lemma \ref{mapeff} and the uniformization lemma show the existence of a partial map $\tilde H\! :\!\omega\!\rightarrow\!\omega$, which is $\Ca$-recursive on $D_0$, and such that $\tilde H(m)$ is in 
${\cal W}^{\omega^\omega\times 2^\omega\times\mathbb{R}}$ and codes the graph of a $\Borel$-recurive map $h\! :\!\omega^\omega\!\times\! 2^\omega\!\rightarrow\!\mathbb{R}$ with
$$1\! -\! h(\alpha ,\beta )\! :=\!\left\{\!\!\!\!\!\!
\begin{array}{ll}
& 0\mbox{ if }(\alpha ,\beta )\!\in\! G\cr
& \mbox{sup}\{\frac{(p)_1}{2^{(p)_0}}\mid p\!\in\! D~\wedge ~
(\alpha ,\beta )\!\notin\! {\cal C}_{\tilde F(m)(p)}\}
\mbox{ if }(\alpha ,\beta )\!\notin\! G
\end{array}
\right.$$
 if $m\!\in\! D_0$. We set $P_1\! :=\!\{ c\!\in\! {\cal W}^{\omega^\omega\times 2^\omega\times\mathbb{R}}\mid  {\cal C}_c\mbox{ is the graph of a function }\zeta_c\}$. It is routine to check that there is a $\Ca$-recursive partial map $I\! :\!\omega^2\!\rightarrow\!\omega$ on its domain $P_1^2$ such that
 $I(c,c')\!\in\! {\cal W}^{\omega^\omega\times 2^\omega\times\mathbb{R}}$ is the graph of the function
 $\mbox{min}(\zeta_c,\zeta_{c'})$ if $c,c'\!\in\! P_1$. We set
 $$r_1(n_0,n_1,n_2,n_3,n)\! :=\! I\Big( n_1,\tilde H\big( r_0(n_0,n_1,n_2,n_3,n)\big)\Big)\mbox{,}$$
 so that $r_1$ is $\Ca$-recursive on its $\Ca$ domain
 $\omega\!\times\! P_1\!\times\!\omega\!\times\! P_3\!\times\!\omega$.\bigskip

\noindent (3) We now define $r_2(n_0,n_1,n_2,n_3,n)$ coding
$$S_{n+1}\! =\!\left\{\!\!\!\!\!\!\!
\begin{array}{ll}
& S_n\! +\! g_{n+1}\mbox{ if }n\mbox{ is odd,}\cr
& S_n\! -\! g_{n+1}\mbox{ if }n\mbox{ is even.}
\end{array}
\right.$$
It is routine to check that there is a $\Ca$-recursive partial map $S\! :\!\omega^3\!\rightarrow\!\omega$ on its domain $P_1^2\!\times\!\omega$ such that
$S(c,c',n)\!\in\! {\cal W}^{\omega^\omega\times 2^\omega\times\mathbb{R}}$ codes the graph of the function
$$(\alpha ,\beta )\!\mapsto\!\left\{\!\!\!\!\!\!\!
\begin{array}{ll}
& \zeta_c(\alpha ,\beta )\! +\!\zeta_{c'}(\alpha ,\beta )\mbox{ if }n\mbox{ is odd}\cr
& \zeta_c(\alpha ,\beta )\! -\!\zeta_{c'}(\alpha ,\beta )\mbox{ if }n\mbox{ is even}
\end{array}
\right.$$
if $(c,c',n)\!\in\! P_1^2\!\times\!\omega$. We set
$r_2(n_0,n_1,n_2,n_3,n)\! :=\! S\big( n_2,r_1(n_0,n_1,n_2,n_3,n),n\big)$, so that $r_2$ is $\Ca$-recursive on its $\Ca$ domain $\omega\!\times\! P_1^2\!\times\! P_3\!\times\!\omega$.\bigskip

\noindent (4) We now define $r_3(n_0,n_1,n_2,n_3,n)$ coding the graph of the function
$\alpha\!\mapsto\! (s^{n+1,\alpha}_j)_{j\in I_{n+1,\alpha}}$. We want to ensure the two following conditions:
$$\begin{array}{ll}
& (1)~G_\alpha\!\subseteq\!\bigcup_{j\in I_{n+1,\alpha}}~N_{s^{n+1,\alpha}_j}\!\subseteq\! G^*_{n+1,\alpha}\cr
& (6)~\vert\fint_{N_{s^{n+1,\alpha}_j}}S_{n+1}(\alpha ,.)~d\lambda\! -\! S_{n+1}(\alpha ,\beta )\vert\! <\! 2^{-3}\mbox{ if }j\!\in\! I_{n+1,\alpha}\ \wedge\
\beta\!\in\! G_\alpha\cap N_{s^{n+1,\alpha}_j}
\end{array}$$
Note first that in practice
$$S_{n+1}(\alpha ,\beta )\! =\!\left\{\!\!\!\!\!\!\!
\begin{array}{ll}
& 0\mbox{ if }n\mbox{ is even}\cr
& 1\mbox{ if }n\mbox{ is odd}
\end{array}
\right.$$
if $(\alpha ,\beta )\!\in\! G$ since $g_p(\alpha ,\beta )\! =\! 1$ for each $p$ in this case. So there is
$\psi\! :\!\omega\!\rightarrow\!\mathbb{R}^2$ recursive with
$$\vert\fint_{N_{s^{n+1,\alpha}_j}}S_{n+1}(\alpha ,.)~d\lambda\! -\! S_{n+1}(\alpha ,\beta )\vert\! <\! 2^{-3}\Leftrightarrow
\psi_0(n)\! <\!\fint_{N_{s^{n+1,\alpha}_j}}S_{n+1}(\alpha ,.)~d\lambda\! <\! \psi_1(n)$$
if $(\alpha ,\beta )\!\in\! G$. We use Corollary \ref{Delta} and its proof. Note that $r_2(n_0,n_1,n_2,n_3,n)\!\in\! P_1$.

\vfill\eject

 We first consider $n'_0\!\in\! {\cal W}_1$ and $n'_2\!\in\! P_1$ (coding $G^*_{n+1}$ and $S_{n+1}$ respectively) as variables. We define 
$R_0,R_1\!\subseteq\!\omega\!\times\!\omega^\omega\!\times\! 2^\omega\!\times\!\omega^3$ by
$$\begin{array}{ll}
R_0(n'_2,\alpha ,\beta ,u,k,l)\Leftrightarrow
& \exists r\!\in\!\mathbb{R}~~\neg\big( n'_2\!\in\! {\cal W}^{\omega^\omega\times 2^\omega\times\mathbb{R}}
~\wedge~(n'_2,\alpha ,\beta ,r)\!\notin\! {\cal C}^{\omega^\omega\times 2^\omega\times\mathbb{R}}\big) ~\wedge\cr
& \big(\frac{k}{2^l}\!\leq\! r\! <\!\frac{k\! +\! 1}{2^l}~\wedge ~\mbox{Seq}(u)~\wedge ~\beta\!\in\! N^u\big)\cr\cr
R_1(n'_2,\alpha ,\beta ,u,k,l)\Leftrightarrow
& \forall r\!\in\!\mathbb{R}~~\big( n'_2\!\in\! {\cal W}^{\omega^\omega\times 2^\omega\times\mathbb{R}}
~\wedge~ (n'_2,\alpha ,\beta ,r)\!\notin\! {\cal C}^{\omega^\omega\times 2^\omega\times\mathbb{R}}\big) ~\vee\cr
& \big(\frac{k}{2^l}\!\leq\! r\! <\!\frac{k\! +\! 1}{2^l}~\wedge ~\mbox{Seq}(u)~\wedge ~\beta\!\in\! N^u\big)\mbox{,}
\end{array}$$
so that $R_0$ is $\Ana$, $R_1$ is $\Ca$, and $R_0(n'_2,\alpha ,\beta ,u,k,l)\Leftrightarrow R_1(n'_2,\alpha ,\beta ,u,k,l)$ if $n'_2\!\in\! P_1$. Then, as in the proof of Corollary \ref{Delta}.(d), we define $O_0,O_1\!\subseteq\!\omega\!\times\!\omega^\omega\!\times\! 2^\omega$ by
$$O_\varepsilon (n'_2,\alpha ,\beta )\Leftrightarrow\mbox{Seq}\big(\alpha (0)\big) ~\wedge ~\mbox{lh}\big(\alpha (0)\big)\! =\! 3~\wedge ~
R_\varepsilon\Big( n'_2,\alpha^*,\beta ,\big(\alpha (0)\big)_0,\big(\alpha (0)\big)_1,\big(\alpha (0)\big)_2\Big)$$
if $\varepsilon\!\in\! 2$, so that $O_0$ is $\Ana$, $O_1$ is $\Ca$, and $O_0(n'_2,\alpha ,\beta )\Leftrightarrow O_1(n'_2,\alpha ,\beta )$ if $n'_2\!\in\! P_1$. In particular, $n'_2\!\in\! P_1$ and $\mbox{Seq}(u)$ imply that 
$$\int_{N^u}S_{n+1}(\alpha ,.)~d\lambda\! =\!
\mbox{lim}_{l\rightarrow\infty}~\Sigma_{k\leq 2^l}~\frac{k}{2^l}\lambda\big( (O_\varepsilon )_{n'_2,<u,k,l>\alpha}\big)$$
for each $\varepsilon\!\in\! 2$. Thus $a\! <\!\int_{N^u}S_{n+1}(\alpha ,.)~d\lambda\! <\! b$ is in this case equivalent to
$$\exists p_0,p_1,q_0,q_1,N\!\in\!\omega ~a\! <\!\frac{p_0}{p_1\! +\! 1}\wedge 
\frac{q_0}{q_1\! +\! 1}\! <\! b\wedge\forall l\!\geq\! N~
\frac{p_0}{p_1\! +\! 1}\!\leq\!\Sigma_{k\leq 2^l}\frac{k}{2^l}\lambda\big( (O_\varepsilon )_{n'_2,<u,k,l>\alpha}\big)\!\leq\!\frac{q_0}{q_1\! +\! 1}.$$
By Corollary \ref{Delta}.(b) applied to $D\! :=\! P_1$, the partial map $\lambda_O\! :\! P_1\!\times\!\omega^\omega\!\rightarrow\!\mathbb{R}$ defined by 
$$\lambda_O(n'_2,\alpha )\! :=\!\lambda\big( (O_0)_{n'_2,\alpha}\big)$$ 
is $\Ana$-recursive and $\Ca$-recursive on its domain. By 3E.2, 3G.1 and 3G.2 in [M], these two classes of functions are closed under composition. In particular, the partial map 
$$(n'_2,\alpha ,u,l)\!\mapsto\!\Sigma_{k\leq 2^l}\frac{k}{2^l}\lambda\big( (O_\varepsilon )_{n'_2,<u,k,l>\alpha}\big)$$ 
is $\Ana$-recursive and $\Ca$-recursive on $P_1\!\times\!\omega^\omega\!\times\!\omega^2$. This shows the existence of 
$Q_0\!\in\!\Ana (\omega^2\!\times\!\omega^\omega\!\times\!\omega )$ and $Q_1\!\in\!\Ca (\omega^2\!\times\!\omega^\omega\!\times\!\omega )$ such that 
$$Q_0(n'_2,n,\alpha ,u)\Leftrightarrow Q_1(n'_2,n,\alpha ,u)\Leftrightarrow
\mbox{Seq}(u)\ \wedge\ \psi_0(n)\! <\!\fint_{N^u}S_{n+1}(\alpha ,.)~d\lambda\! <\! \psi_1(n)$$ 
if $n'_2\!\in\! P_1$. We now consider $n'_0\!\in\! {\cal W}_1$ and $n'_2\!\in\! P_1$ as parameters. We set\bigskip

\leftline{$P_{n'_0,n'_2}(n,\alpha ,u)\Leftrightarrow$}\smallskip

\rightline{$Q_1(n'_2,n,\alpha ,u)~\wedge ~
N^u\!\subseteq\! {\cal C}_{n'_0,\alpha}~\wedge ~\forall k\! <\!\mbox{lh}(u)~~
\big(\neg Q_0\big( n'_2,n,\alpha ,\underline{u}(k)\big) ~\vee ~
N^{\underline{u}(k)}\!\not\subseteq\! {\cal C}_{n'_0,\alpha}\big).$}\bigskip

\noindent Note that for each $(\alpha ,\beta )\!\in\! G$ there is $l\!\in\!\omega$ minimal with the properties that $N_{\beta\vert l}\!\subseteq\! {\cal C}_{n'_0,\alpha}$ and 
$Q_1\big( n'_2,n,\alpha ,<\beta (0),...,\beta (l\! -\! 1)>\!\!\big)$, so that
$P_{n'_0,n'_2}\big( n,\alpha ,<\beta (0),...,\beta (l\! -\! 1)>\!\!\big)$ since $n'_0\!\in\! {\cal W}_1$ and $n'_2\!\in\! P_1$. As $n'_0\!\in\! {\cal W}_1$, 
$N^{\underline{u}(k)}\!\setminus\! {\cal C}_{n'_0,\alpha}$ is a $\Borel (\alpha )$ compact subset of $2^\omega$, so that it contains a $\Borel (\alpha )$ point if it is not empty (see 4F.15 in [M]). This shows that $P_{n'_0,n'_2}$ is $\Ca$.

\vfill\eject

 The uniformization lemma provides a $\Borel$-recursive map $L\! :\!\omega\!\times\!\omega^\omega\!\times\! 2^\omega\!\rightarrow\!\omega$ such that 
$$P_{n'_0,n'_2}\Big(n,\alpha ,<\beta (0),...,\beta\big( L(n,\alpha ,\beta )\! -\! 1\big)\! >\!\!\Big)$$ 
if $(\alpha ,\beta )\!\in\! G$. Note that the $\Ana$ set
$$\sigma\! :=\!\big\{ (n,\alpha ,u)\!\in\!\omega\!\times\!\omega^\omega\!\times\!\omega\mid\exists\beta\!\in\! G_\alpha ~~
u\! =<\beta (0),...,\beta\big( L(n,\alpha ,\beta )\! -\! 1\big)\! >\!\!\big\}$$ 
is contained in the $\Ca$ set $\pi\! :=\!\{ (n,\alpha ,u)\!\in\!\omega\!\times\!\omega^\omega\!\times\!\omega\mid P_{n'_0,n'_2}(n,\alpha ,u)\}$. By 7B.3 in [M], there is a $\Borel$ subset $\delta$ of $\omega\!\times\!\omega^\omega\!\times\!\omega$ such that $\sigma\!\subseteq\!\delta\!\subseteq\!\pi$. We now also consider $n$ as a parameter and define
$\varphi\! :\!\omega^\omega\!\rightarrow\!\omega^\omega$ by
$$\varphi (\alpha )(u)\! :=\!\left\{\!\!\!\!\!\!\!
\begin{array}{ll}
& u\mbox{ if }(n,\alpha ,u)\!\in\!\delta\mbox{,}\cr
& 0\mbox{ otherwise.}
\end{array}
\right.$$
Note that $\varphi$ is $\Borel$-recursive, and that $\mbox{Seq}\big(\varphi (\alpha )(u)\big)$ is equivalent to $(n,\alpha ,u)\!\in\!\delta$. In particular,
$$\begin{array}{ll}
& (1)~G_\alpha\!\subseteq\!\bigcup\big\{ N^{\varphi (\alpha )(u)}\mid u\!\in\!\omega\wedge\mbox{Seq}\big(\varphi (\alpha )(u)\big)\big\}\!\subseteq\! 
{\cal C}_{n'_0,\alpha}\cr
& (6)~\vert\fint_{N^{\varphi (\alpha )(u)}}S_{n+1}(\alpha ,.)~d\lambda\! -\! S_{n+1}(\alpha ,\beta )\vert\! <\! 2^{-3}\mbox{ if }\mbox{Seq}\big(\varphi (\alpha )(u)\big)\ \wedge\ \beta\!\in\! G_\alpha\cap N^{\varphi (\alpha )(u)}
\end{array}$$
for each $\alpha\!\in\!\omega^\omega$. Let $k\!\in\! {\cal W}^{\omega^\omega\times\omega^\omega}$ such that
${\cal C}^{\omega^\omega\times\omega^\omega}_k\! =\!\mbox{Gr}(\varphi )$. We now consider $n'_0$, $n'_2$ and $n$ as variables again. Note that for each 
$(n'_0,n'_2,n)\!\in\! {\cal W}_1\!\times\! P_1\!\times\!\omega$ there is $k\!\in\!\omega$ such that 
$$R(n'_0,n'_2,n,k)\Leftrightarrow
\left\{\!\!\!\!\!\!
\begin{array}{ll}
& k\!\in\! {\cal W}^{\omega^\omega\times\omega^\omega}\ \wedge\cr
& \Big(\forall\alpha\!\in\!\omega^\omega ~~\forall\gamma\!\in\!\omega^\omega ~~\big( k\!\in\! {\cal W}^{\omega^\omega\times\omega^\omega}\wedge\neg 
{\cal C}^{\omega^\omega\times\omega^\omega}(k,\alpha ,\gamma )\big)\vee\cr
& \big( (1)~G_\alpha\!\subseteq\!\bigcup\big\{ N^{\gamma (u)}\mid u\!\in\!\omega\wedge\mbox{Seq}\big(\gamma (u)\big)\big\}\!\subseteq\! 
{\cal C}_{n'_0,\alpha}\cr 
& \wedge\ (6)~\forall u\!\in\!\omega ~~\neg\mbox{Seq}\big(\gamma (u)\big)\ \vee\ 
Q_1(n'_2,n,\alpha ,u)\big)\Big)
\end{array}
\right.$$
Note that $R\!\in\!\Ca (\omega^4)$. The uniformization lemma provides a partial map 
$K\! :\! \omega^3\!\mapsto\!\omega$ which is $\Ca$-recursive on its 
$\Ca$ domain ${\cal W}_1\!\times\! P_1\!\times\! \omega$, and 
$R\big( n'_0,n'_2,n,K(n'_0,n'_2,n)\big)$ if 
$$(n'_0,n'_2,n)\!\in\! {\cal W}_1\!\times\! P_1\!\times\! \omega .$$ 
It remains to set $r_3(n_0,n_1,n_2,n_3,n)\! :=\! K(n'_0,n'_2,n)$ if 
$n'_0\! =\! r_0(n_0,n_1,n_2,n_3,n)$ and 
$$n'_2\! =\! r_2(n_0,n_1,n_2,n_3,n)\mbox{,}$$ 
so that $r_3$ is $\Ca$-recursive on its $\Ca$ domain 
${\cal W}_1\!\times\! P_1^2\!\times\! P_3\!\times\!\omega$.\bigskip

 Finally, $r$ is $\Ca$-recursive on ${\cal W}_1\!\times\! P_1^2\!\times\! P_3\!\times\!\omega$, $f$ is $\Ca$-recursive on $\omega$, and thus $f$ is $\Borel$-recursive by the uniformization lemma since it is total.\bigskip

\noindent $\bullet$ We are now ready to define the dimension two versions of $G^*_n$, $g_n$, $S_n$, and $(s^n_j)_{j\in I_n}$:
$$\begin{array}{ll}
& (1)~G^*_n\! :=\! {\cal C}_{f_0(n)}\mbox{,}\cr
& (2)~g_n(\alpha ,\beta )\! =\!\rho\Leftrightarrow\big( f_1(n),\alpha ,\beta ,\rho\big)\!\in\!
{\cal C}^{\omega^\omega\times 2^\omega\times\mathbb{R}}\mbox{,}\cr
& (3)~S_n(\alpha ,\beta )\! =\!\rho\Leftrightarrow\big( f_2(n),\alpha ,\beta ,\rho\big)\!\in\!
{\cal C}^{\omega^\omega\times 2^\omega\times\mathbb{R}}\mbox{,}\cr
& (4)~\left\{\!\!\!\!\!\!\!
\begin{array}{ll}
& (i)~j\!\in\! I_{n,\alpha}\Leftrightarrow\exists\delta\!\in\!\omega^\omega ~~
\big( f_3(n),\alpha ,\delta \big)\!\in\! {\cal C}^{\omega^\omega\times\omega^\omega}~\wedge ~
\mbox{Seq}\big(\delta (j)\big)\mbox{,}\cr
& (ii)~s^{n,\alpha}_j\! =\!\delta (j)\mbox{ if }j\!\in\! I_{n,\alpha}.
\end{array}
\right.
\end{array}$$
By construction of $r$, these objects satisfy the conditions (1)-(6) of the proof of Lemma 
\ref{Gdelta}.

\vfill\eject

\noindent $\bullet$ Consequently, the martingale $F(\alpha )$ will be defined in such a way that if
$u\!\in\!\omega$ codes $s\!\in\! 2^{<\omega}$, then
$F(\alpha )(s)\! =\!\fint_{N^u}f_\infty (\alpha ,.)~d\lambda$. Note that
$G\! =\!\bigcap_{n\in\omega}~G^*_n$, so that $\neg G$ is the disjoint union of the
$G^*_n\!\setminus\! G^*_{n+1}$'s. Thus
$$\begin{array}{ll}
\int_{N^u}f_\infty (\alpha ,.)~d\lambda\!\!\!\!
& \! =\!\int_{N^u\setminus G_\alpha}f_\infty (\alpha ,.)~d\lambda\! =\!\Sigma_{n\in\omega}~
\int_{N^u\cap (G^*_n)_\alpha\setminus (G^*_{n+1})_\alpha}f_\infty (\alpha ,.)~d\lambda\cr
& \! =\!\Sigma_{n\in\omega}~\Sigma_{j\leq n}~(-1)^j
\int_{N^u\cap (G^*_n)_\alpha\setminus (G^*_{n+1})_\alpha}g_j(\alpha ,.)~d\lambda\cr
& \! =\!\mbox{lim}_{l\rightarrow\infty}~\Sigma_{n\leq l}~\Sigma_{j\leq n}~(-1)^j
\int_{N^u\cap (G^*_n)_\alpha\setminus (G^*_{n+1})_\alpha}g_j(\alpha ,.)~d\lambda .
\end{array}$$
Consequently, in order to prove that $F$ is $\Borel$-recursive, it is enough to check that the partial map $(u,\alpha ,j,n)\!\mapsto\!
\int_{N^u\cap (G^*_n)_\alpha\setminus (G^*_{n+1})_\alpha}g_j(\alpha ,.)~d\lambda$ is $\Borel$-recursive from $\{ u\!\in\!\omega\mid\mbox{Seq}(u)\}\!\times\!\omega^\omega\!\times\!\omega^2$ into $\mathbb{R}$. By Corollary \ref{Delta}, it is enough to check that the map
$h\! :\!\omega^\omega\!\times\! 2^\omega\!\rightarrow\!\mathbb{R}$ defined by
$$h(\alpha ,\beta )\! :=\!\left\{\!\!\!\!\!\!\!
\begin{array}{ll}
& g_{(\alpha (0))_0}(\alpha^*,\beta )\mbox{ if }\mbox{Seq}\big(\alpha (0)\big) ~\wedge ~
\mbox{lh}\big(\alpha (0)\big)\! =\! 2~\wedge ~
(\alpha^*,\beta )\!\in\! G^*_{(\alpha (0))_1}\!\setminus\! G^*_{(\alpha (0))_1+1}\mbox{,}\cr
& 0\mbox{ otherwise,}
\end{array}
\right.$$
is $\Borel$-recursive. This comes from the facts that
$$(\alpha ,\beta )\!\in\! G^*_n
\Leftrightarrow\big( f_0(n),\alpha ,\beta\big)\!\in\! {\cal C}
\Leftrightarrow\neg\Big( f_0(n)\!\in\! {\cal W} ~\wedge ~\big( f_0(n),\alpha ,\beta\big)\!\notin\! {\cal C}\Big)$$
is $\Borel$ in $(\alpha ,\beta ,n)$ and
$$\begin{array}{ll}
g_n(\alpha ,\beta )\!\in\! N(\mathbb{R},p)\!\!\!
& \Leftrightarrow\exists\rho\!\in\!\mathbb{R}~~\neg
\Big( f_1(n)\!\in\! {\cal W}^{\omega\times 2^\omega\times\mathbb{R}}~\wedge ~
\big( f_1(n),\alpha ,\beta ,\rho\big)\!\notin\! {\cal C}^{\omega\times 2^\omega\times\mathbb{R}}\Big)
~\wedge\cr
& \hfill{\rho\!\in\! N(\mathbb{R},p)}\cr
& \Leftrightarrow\forall\rho\!\in\!\mathbb{R}~~
\Big( f_1(n)\!\in\! {\cal W}^{\omega\times 2^\omega\times\mathbb{R}}~\wedge ~
\big( f_1(n),\alpha ,\beta ,\rho\big)\!\notin\! {\cal C}^{\omega\times 2^\omega\times\mathbb{R}}\Big)
~\vee\cr
& \hfill{\rho\!\in\! N(\mathbb{R},p)}
\end{array}$$
is $\Borel$ in $(\alpha ,\beta ,n,p)$.\bigskip

\noindent $\bullet$ Finally, the map $F$ is $\Borel$-recursive and is as required.\hfill{$\square$}

\section{$\!\!\!\!\!\!$ First consequences}

\bf (A) Universal sets\rm\bigskip

\noindent - We first recall some material from [K2]. The first result can be found in Section 23.F (see also [Za]).

\begin{thm} (Zahorski) Let $B$ be a subset of $[0,1]$. The following are equivalent:\smallskip

\noindent (a) there are $S\!\in\!\boratwo$ and $P\!\in\!\bormthree$ with $m(P)\! =\! 1$, where $m$ is the Lebesgue measure on $[0,1]$, such that $B\! =\! S\cap P$,\smallskip

\noindent (b) there is $f\!\in\! C([0,1])$ with $B\! =\!\{ x\!\in\! [0,1]\mid f'(x)\mbox{ exists}\}$ (we consider only one-sided derivatives at the endpoints).\end{thm}

 The second result is 23.23.

\begin{thm} Let $\cal G$ be a $G_\delta$ subset of $(0,1)$ with $m({\cal G})\! =\! 0$. Then
$$\{ (f,x)\!\in\! C([0,1])\!\times\! {\cal G}\mid f'(x)\mbox{ exists}\}$$
is $C([0,1])$-universal for $\bormthree ({\cal G})$.\end{thm}

\noindent - We prove results in that spirit here.

\begin{thm} Let $B$ be a subset of $2^\omega$. Then the following are equivalent:\smallskip

\noindent (a) $B$ is $\borathree$ and has $\lambda$-measure zero,\smallskip

\noindent (b) there is $f\!\in\! {\cal M}$ with
$B\! =\!\{\beta\!\in\! 2^\omega\mid\mbox{osc}(f,\beta )\! >\! 0\}$.\end{thm}

\noindent\bf Proof.\rm\ (a) $\Rightarrow$ (b) Write $B\! =\!\bigcup_{n\in\omega}~G_n$, where the $G_n$'s are $G_\delta$. Lemma \ref{Gdelta} gives, for each $n$, a martingale $f_n$ with
$G_n\! =\! D(f_n)$ and
$\{\mbox{osc}(f_n,\beta )\mid\beta\!\in\! 2^\omega\}\!\subseteq\!\{ 0\}\cup [\frac{1}{2},1]$. Lemma
\ref{union} gives $f\!\in\! {\cal M}$ with $D(f)\! =\! B$.\bigskip

\noindent (b) $\Rightarrow$ (a) We already noticed in the introduction that $B$ is $\borathree$. By Doob's theorem, $B$ has $\lambda$-measure zero (see [D]).\hfill{$\square$}

\begin{cor} Let $\cal G$ be a $G_\delta$ subset of $2^\omega$ with $\lambda ({\cal G})\! =\! 0$. Then $\{ (f,\beta )\!\in\! {\cal M}\!\times\! {\cal G}\mid\mbox{osc}(f,\beta )\! >\! 0\}$ is $\cal M$-universal for $\borathree ({\cal G})$.\end{cor}

 For example, $\{\beta\!\in\! 2^\omega\mid\forall n\!\in\!\omega ~~\beta (2n)\! =\! 0\}$ is a $\Bormone$ copy of $2^\omega$ and has $\lambda$-measure zero.\bigskip

\noindent\bf (B) Complete sets\rm\bigskip

\noindent - By 33.G in [K2], there is a uniform version of Zahorski's theorem, which allows to prove the following result

\begin{thm} (Mazurkiewicz) The set of differentiable functions in $C([0,1])$ is $\ca$-complete.\end{thm}

\noindent - Here again, there is a result in that spirit.

\begin{thm} \label{pi11complet} The set ${\cal P}\! :=\!\{ f\!\in\! {\cal M}\mid\forall\beta\!\in\! 2^\omega ~~
\mbox{osc}(f,\beta )\! =\! 0\}$ is $\ca$-complete.\end{thm}

\noindent\bf Notation.\rm\ Let
${\cal K}\! :=\!\{\beta\!\in\! 2^\omega\mid\forall n\!\in\!\omega ~~\beta (2n)\! =\! 0\}$, which is a
$\Bormone$ copy of the Cantor space $2^\omega$ with $\lambda ({\cal K})\! =\! 0$. In particular,
$\cal K$ is a nonempty $G_\delta\cap\Borel$ subset of $2^\omega$.\bigskip

\noindent\bf Proof.\rm\ Let $U\!\in\!\Ca (\omega^\omega\!\times\! 2^\omega )$ be
$\omega^\omega$-universal for the co-analytic subsets of $2^\omega$, and
$$\Pi\! :=\!\{\alpha\!\in\!\omega^\omega\mid\big( (\alpha )_0,(\alpha )_1\big)\!\in\! U\} .$$
Note that $\Pi\!\in\!\Ca$. If $P\!\in\!\ca (2^\omega )$, then $P\! =\! U_\alpha$ for some
$\alpha\!\in\!\omega^\omega$, so that the map $\beta\!\mapsto <\alpha ,\beta >$ is a continuous reduction of $P$ to $\Pi$ and $\Pi$ is $\ca$-complete. Let
$H\!\in\!\Bormtwo (\omega^\omega\!\times\! 2^\omega )$ with $\neg\Pi\! =\!\Pi_0[H]$. We set
$G\! :=\!\big\{ (\alpha ,\beta )\!\in\!\omega^\omega\!\times\! 2^\omega\mid
\big(\alpha ,(\beta )_1\big)\!\in\! H~\wedge ~\beta\!\in\! {\cal K}\big\}$, so that
$G\!\in\!\Borel (\omega^\omega\!\times\! 2^\omega )$, has $G_\delta$ vertical sections and
$G\!\subseteq\!\omega^\omega\!\times\! {\cal K}$. Lemma \ref{Borel} gives
$F\! :\!\omega^\omega\!\rightarrow\! {\cal M}$ Borel such that
$G_\alpha\! =\! {\cal V}_{F(\alpha )}$ for each $\alpha\!\in\!\omega^\omega$.

 Thus
$$\alpha\!\notin\!\Pi\Leftrightarrow\exists\beta\!\in\! 2^\omega ~~(\alpha ,\beta )\!\in\! H
\Leftrightarrow\exists\beta\!\in\! 2^\omega ~~(\alpha ,\beta )\!\in\! G
\Leftrightarrow\exists\beta\!\in\! 2^\omega ~~\big( F(\alpha ),\beta\big)\!\in\! {\cal V}
\Leftrightarrow F(\alpha )\!\notin\! {\cal P}.$$
Thus $\Pi\! =\! F^{-1}({\cal P})$ and $\cal P$ is Borel $\ca$-complete. By 26.C in [K2], $\cal P$ is
$\ca$-complete.\hfill{$\square$}\bigskip

\noindent - We now prove Theorem \ref{completes}. Let $X$ be a metrizable compact space and $Y$ be a Polish space. We equip ${\cal C}(X,Y)$ with the topology of uniform convergence, so that it is a Polish space (see 4.19 in [K2]). We use the map $\psi$ defined before Theorem \ref{completes}.

\begin{thm} (a) The set ${\cal P}_1\! :=\!\big\{ (f_k)_{k\in\omega}\!\in\! {\cal P}^\omega\mid
\big( \psi (f_k)\big)_{k\in\omega}\mbox{ pointwise converges}\big\}$ is $\ca$-complete.\smallskip

\noindent (b) The set ${\cal P}_2\! :=\!\big\{ (f_k)_{k\in\omega}\!\in\! {\cal P}^\omega\mid
\big( \psi (f_k)\big)_{k\in\omega}\mbox{ pointwise converges to zero}\big\}$ is $\ca$-complete.\smallskip

\noindent (c) The set ${\cal S}\! :=\!\big\{ (f_k)_{k\in\omega}\!\in\! {\cal P}^\omega\mid\exists\gamma\!\in\!\omega^\omega ~~
\big( \psi (f_{\gamma (i)})\big)_{i\in\omega}\mbox{ pointwise converges to zero}\big\}$ is ${\bf\Sigma}^1_2$-complete.\end{thm}

\noindent\bf Proof.\rm\ We define $\varphi\! :\! {\cal C}(2^\omega ,[0,1])\!\rightarrow\! {\cal M}$ by $\varphi (h)(s)\! :=\!\fint_{N_s}h~d\lambda$. As in the proof of Lemma \ref{Gdelta}, $\varphi$ is well-defined. It is also continuous, and injective: if $h\!\not=\! h'$, then we can find $q\!\in\!\omega$ and
$s\!\in\! 2^{<\omega}$ such that $h(\beta )\! -\! h'(\beta )\! >\! 2^{-q}$ for each $\beta\!\in\! N_s$ or $h'(\beta )\! -\! h(\beta )\! >\! 2^{-q}$ for each $\beta\!\in\! N_s$, so that
$$\vert\varphi (h)(s)\! -\!\varphi (h')(s)\vert\! =\!\frac{1}{\lambda (N_s)}\vert\int_{N_s}h~d\lambda\! -\!\int_{N_s}h'~d\lambda\vert\!\geq\! 2^{-q}.$$
This implies that the range $\cal R$ of $\varphi$ is Borel and
$\psi\! :=\!\varphi^{-1}\! :\! {\cal R}\!\rightarrow\! {\cal C}(2^\omega ,[0,1])$ is Borel. As every continuous map $h\! :\! 2^\omega\!\rightarrow\! [0,1]$ is $\tau$-continuous,
$$\mbox{lim}_{l\rightarrow\infty}~\varphi (h)(\beta\vert l)\! =\!
\mbox{lim}_{l\rightarrow\infty}~\fint_{N_{\beta\vert l}}h~d\lambda\! =\! h(\beta )$$
for each $\beta\!\in\! 2^\omega$, by Lemma \ref{moy}. This implies that $f\!\in\! {\cal P}$ and
$\psi (f)(\beta )\! =\!\mbox{lim}_{l\rightarrow\infty}~f(\beta\vert l)$ for each $\beta\!\in\! 2^\omega$ if
$f\!\in\! {\cal R}$.\bigskip

\noindent (a) Note that the proof of 33.11 in [K2] shows that the set
$$P_1\! :=\!\big\{ (h_k)_{k\in\omega}\!\in\!\big( {\cal C}(2^\omega ,[0,1])\big)^\omega\mid (h_k)_{k\in\omega}\mbox{ pointwise converges}\big\}$$
is ${\bf\Pi}^1_1$-complete. As ${\cal E}\! :=\!\big\{ (f_k)_{k\in\omega}\!\in\! {\cal R}^\omega\mid
\big( \psi (f_k)\big)_{k\in\omega}\mbox{ pointwise converges}\big\}\! =\! (\psi^\omega )^{-1}(P_1)$, the equalities
$P_1\! =\! (\varphi^\omega )^{-1}({\cal E})\! =\! (\varphi^\omega )^{-1}({\cal P}_1)$ hold and ${\cal P}_1$ is ${\bf\Pi}^1_1$-complete.\bigskip

\noindent (b) We argue as in (a).\bigskip

\noindent (c) As in [B-Ka-L], the set
$$S\! :=\!\big\{ (h_k)_{k\in\omega}\!\in\!\big( {\cal C}(2^\omega ,[0,1])\big)^\omega\mid\exists\gamma\!\in\!\omega^\omega ~~
\big( h_{\gamma (i)}\big)_{i\in\omega}\mbox{ pointwise converges to zero}\big\}\mbox{,}$$
is ${\bf\Sigma}^1_2$-complete. Indeed, fix $Q\!\in\! {\bf\Sigma}^1_2(2^\omega )$.

\vfill\eject

 Lemma 2.2 in [B-Ka-L] gives $(g_k)_{k\in\omega}\!\in\!\big( {\cal C}(2^\omega\!\times\! 2^\omega ,2)\big)^\omega$ such that, for each $\delta\!\in\! 2^\omega$, the following are equivalent:\bigskip

(i) $\delta\!\in\! Q$,\smallskip

(ii) $\exists\gamma\!\in\!\omega^\omega ~~\forall\beta\!\in\! 2^\omega ~~\mbox{lim}_{i\rightarrow\infty}~g_{\gamma (i)}(\delta ,\beta )\! =\! 0$.\bigskip

We define, $g\! :\! 2^\omega\!\rightarrow\!\big( {\cal C}(2^\omega ,[0,1])\big)^\omega$ by
$g(\delta )(k)(\beta )\! :=\! g_k(\delta ,\beta )$. Then $g$ is continuous and reduces $Q$ to $S$. As
$${\cal E}'\! :=\!\big\{ (f_k)_{k\in\omega}\!\in\! {\cal R}^\omega\mid\exists\gamma\!\in\!\omega^\omega ~~
\big( \psi (f_{\gamma (i)})\big)_{i\in\omega}\mbox{ pointwise converges to zero}\big\}\! =\! (\psi^\omega )^{-1}(S)\mbox{,}$$
$S\! =\! (\varphi^\omega )^{-1}({\cal E}')\! =\! (\varphi^\omega )^{-1}({\cal S})$ and ${\cal S}$ is ${\bf\Sigma}^1_2$-complete.\hfill{$\square$}

\section{$\!\!\!\!\!\!$ Universal and complete sets in the spaces ${\cal C}(2^\omega ,X)$}

- It is known that if $\bf\Gamma$ is a self-dual Wadge class and $X$ is a Polish space, then there is no set which is $X$-universal for the subsets of $X$ in $\bf\Gamma$ (see 22.7 in [K2]). This is no longer the case if the space of codes is different from the space of coded sets.

\begin{prop} \label{Wadge} Let $X$ be a Polish space, $\bf\Gamma$ be a Wadge class with complete set $C\!\in\! {\bf\Gamma}(X)$, and ${\cal U}^{\bf\Gamma}\! :=\!
\{ (h,\beta )\!\in\! {\cal C}(2^\omega ,X)\!\times\! 2^\omega\mid h(\beta )\!\in\! C\}$. Then
${\cal U}^{\bf\Gamma}$ is ${\cal C}(2^\omega ,X)$-universal for the $\bf\Gamma$ subsets of
$2^\omega$.\end{prop}

\noindent\bf Proof.\rm\ As the evaluation map $(h,\beta )\!\mapsto\! h(\beta )$ is continuous,
${\cal U}^{\bf\Gamma}\!\in\! {\bf\Gamma}$. If $A\!\in\! {\bf\Gamma}(2^\omega )$, then
$A\! =\! h^{-1}(C)$ for some $h\!\in\! {\cal C}(2^\omega ,X)$, so that
$A\! =\! {\cal U}^{\bf\Gamma}_h$.\hfill{$\square$}\bigskip

 We will partially strengthen this result to get our uniform universal sets.\bigskip

\noindent - Recall that it is proved in [K3] that a Borel $\ca$-complete set is actually $\ca$-complete. In fact, Kechris's proof shows the result for the classes ${\bf\Pi}^1_n$. Our main tool is a uniform version of this. Kechris's result has recently been strengthened in [P] as follows.

\begin{thm} (Pawlikowski) \label{Paw} Let $n\!\geq\! 1$ be a natural number, and
$C\!\subseteq\! X\!\subseteq\! 2^\omega$. If Borel functions from $2^\omega$ into $X$ give as preimages of $C$ all ${\bf\Pi}^1_n$ subsets of $2^\omega$, then so do continuous injections.
\end{thm}

 The main tool mentioned above is the following:

\begin{thm} \label{tool} Let $n\!\geq\! 1$ be a natural number,
${\cal U}^{{\bf\Pi}^1_n,2^\omega}$ be a suitable $\omega^\omega$-universal set for the
${\bf\Pi}^1_n$ subsets of $2^\omega$, $X$ be a recursively presented Polish space,
$C\!\in\! {\it\Pi}^1_n(X)$,
${\cal R}\! :\!\omega^\omega\!\times\!\omega^\omega\!\rightarrow\!\omega^\omega$ be a recursive map, and $b\! :\!\omega^\omega\!\rightarrow\! X$ be a $\Borel$-recursive map such that
$$(\alpha ,\beta )\!\in\! {\cal U}^{{\bf\Pi}^1_n,2^\omega}\Leftrightarrow
b\big({\cal R}(\alpha ,\beta )\big)\!\in\! C$$ for each
$(\alpha ,\beta )\!\in\!\omega^\omega\!\times\! 2^\omega$. Then there is a $\Borel$-recursive map $f\! :\!\omega^\omega\!\rightarrow\! {\cal C}(2^\omega ,X)$ such that
$$(\alpha ,\beta )\!\in\! {\cal U}^{{\bf\Pi}^1_n,2^\omega}\Leftrightarrow f(\alpha )(\beta )\!\in\! C$$
for each $(\alpha ,\beta )\!\in\!\omega^\omega\!\times\! 2^\omega$.\end{thm}

\noindent - We first recall some material from [K3].

\begin{defi} (a) A {\bf coding system} for nonempty perfect binary trees is a pair
$({\cal D},{\cal O})$, where ${\cal D}\!\subseteq\! 2^\omega$ and
${\cal O}\! :\! {\cal D}\!\rightarrow\!\{ T\!\in\! 2^{2^{<\omega}}\mid
T\mbox{ is a nonempty perfect binary tree}\}$ is onto.\smallskip

\noindent (b) A coding system $({\cal D},{\cal O})$ is {\bf nice} if\smallskip

(i) for any $\alpha\!\in\!\omega^\omega$ and any $\Borel (\alpha )$-recursive map
$H\! :\! 2^\omega\!\times\! 2^\omega\!\rightarrow\!\omega$, we can find
$\beta\!\in\! {\cal D}\cap\Borel (\alpha )$ and $k\!\in\!\omega$ such that $H(\beta ,\delta )\! =\! k$ for each $\delta$ in the body $[{\cal O}(\beta )]$ of ${\cal O}(\beta )$,\smallskip

(ii) $\cal D$ is $\Ca$ and, for $\beta\!\in\! {\cal D}$, the relation
$$R(m,\beta )\Leftrightarrow\mbox{Seq}(m)~\wedge ~\big( (m)_0,...,(m)_{\mbox{lh}(m)-1}\big)\!\in\! {\cal O}(\beta )$$
is $\Borel$, i.e., there are $\Ca$ relations $\Pi_0,\Pi_1$ such that
$R(m,\beta )\Leftrightarrow\Pi_0(m,\beta )\Leftrightarrow\neg\Pi_1(m,\beta )$ if
$\beta\!\in\! {\cal D}$.\end{defi}

 Nice coding systems exist. If $\beta\!\in\! {\cal D}$, then there is a canonical homeomorphism
$\beta^*$ from $[{\cal O}(\beta )]$ onto $2^\omega$. We now check that the construction of
$\beta^*$ is effective.

\begin{lem} \label{star} (a) The partial function $e\! :\! (\beta ,\delta )\!\mapsto\!\beta^*(\delta)$ is $\Ca$-recursive on its $\Ca$ domain
$$\mbox{Domain}(e)\! :=\!
\{ (\beta ,\delta )\!\in\! {\cal D}\!\times\! 2^\omega\mid\delta\!\in\! [{\cal O}(\beta )]\} .$$
(b) The partial function $\iota\! :\! (\beta ,\gamma )\!\mapsto\!\mbox{ the unique }
\delta\!\in\! [{\cal O}(\beta )]\mbox{ with }\beta^*(\delta )\! =\!\gamma$ is $\Ca$-recursive on its
$\Ca$ domain ${\cal D}\!\times\! 2^\omega$.\end{lem}

\noindent\bf Proof.\rm\ (a) We define a $\Ca$ relation $\cal Q$ on
$\omega^2\!\times\! (2^\omega )^2$ by
$${\cal Q}(p,p',\beta ,\delta )\Leftrightarrow\Big(\big(\forall\varepsilon\!\in\! 2~~
\Pi_0(\overline{(\delta\vert p')\varepsilon},\beta )\big) ~\wedge ~\big(\forall p\!\leq\! p''\! <\! p'~~
\exists\varepsilon\!\in\! 2~~\Pi_1((\overline{\delta\vert p'')\varepsilon},\beta )\big)\Big) .$$
Note that
$$\beta^*(\delta )(n)\! =\!\varepsilon\Leftrightarrow\left\{\!\!\!\!\!\!
\begin{array}{ll}
& \exists l\!\in\!\omega ~~\mbox{Seq}(l)~\wedge ~\mbox{lh}(l)\! =\! n\! +\! 1~\wedge ~
\delta\big( (l)_n\big)\! =\!\varepsilon ~\wedge ~{\cal Q}\big( 0,(l)_0,\beta ,\delta\big)~\wedge\cr
& \forall m\! <\! n~~(l)_m\! <\! (l)_{m+1}~\wedge ~{\cal Q}\big( (l)_m\! +\! 1,(l)_{m+1},\beta ,\delta\big)
\end{array}
\right.$$
if $\beta\!\in\! {\cal D}$. The proof of (b) is similar.\hfill{$\square$}\bigskip

\noindent - Let $X$ be a recursively presented Polish space, and $d_X$ and
$(r^X_n)_{n\in\omega}$ be respectively a distance function and a recursive presentation of $X$. We now give a {\bf recursive presentation of} ${\cal C}(2^\omega ,X)$, equipped with the usual distance defined by
$$d(h,h')\! :=\!\mbox{sup}_{\beta\in 2^\omega}~d_X\big( h(\beta ),h'(\beta )\big)\mbox{,}$$
since this is not present in [M]. We define, by primitive recursion, a recursive map
$\nu\! :\!\omega\!\rightarrow\!\omega$ such that $\nu (i)$ enumerates
$\{ s\!\in\! 2^{<\omega}\mid\vert s\vert\! =\! i\}$. We first set $\nu (0)\! :=\! 1\! =<>$. Then
$$\nu (i\! +\! 1)\! =\! k\Leftrightarrow\mbox{Seq}(k)~\wedge ~\mbox{lh}(k)\! =\! 2^{i+1}~\wedge ~
\forall l\! <\! 2^i~~\forall\varepsilon\!\in\! 2~~(k)_{\varepsilon 2^i+l}
\! =\!\overline{s\Big(\big(\nu (i)\big)_l\Big)\varepsilon} .$$
If $\mbox{Seq}(n)$ and $\mbox{lh}(n)\! =\! 2^i$ for some $i$ ($<\! n$), then we define
$h_n\! :\! 2^\omega\!\rightarrow\! X$ by $h_n(\beta )\! :=\! r^X_{(n)_l}$ if
$$\beta\vert i\! =\! s^i_l\! :=\! s\Big(\big(\nu (i)\big)_l\Big) .$$
If $\neg\mbox{Seq}(n)$ or $\mbox{lh}(n)\!\not=\! 2^i$ for each $i$, then we define
$h_n\! :\! 2^\omega\!\rightarrow\! X$ by $h_n(\beta )\! :=\! r^X_0$ if $\beta\!\in\! 2^\omega$. In any case, $h_n\!\in\! {\cal C}(2^\omega ,X)$ and takes finitely many values.

\begin{lem} \label{pres} Let $X$ be a recursively presented Polish space. Then the sequence
$(h_n)_{n\in\omega}$ is a recursive presentation of ${\cal C}(2^\omega ,X)$, equipped with $d$.\end{lem}

\noindent\bf Proof.\rm\ We have to see that $(h_n)$ is dense in ${\cal C}(2^\omega ,X)$. So let
$h\!\in\! {\cal C}(2^\omega ,X)$, $\epsilon\! >\! 0$ and $m\!\in\!\omega$ with
$2^{-m}\! <\!\frac{\epsilon}{2}$. As $h$ is uniformly continuous, there is $i\!\in\!\omega$ such that
$d_X\big( h(\beta ),h(\delta )\big)\! <\! 2^{-m}$ if $\beta\vert i\! =\!\delta\vert i$. We choose, for each $l\! <\! 2^i$, $n_l\!\in\!\omega$ such that $d_X\big( r^X_{n_l},h(s^i_l0^\infty )\big)\! <\! 2^{-m}$. We set $n\! :=<n_0,...,n_{2^i-1}>$. If $\beta\!\in\! 2^\omega$ and $\beta\vert i\! =\! s^i_l$, then
$d_X\big( h(\beta ),h_n(\beta )\big)\!\leq\! d_X\big( h(\beta ),h(s^i_l0^\infty )\big)\! +\!
d_X\big( h(s^i_l0^\infty ),r^X_{n_l}\big)\!\leq\! 2^{-m}\! +\! 2^{-m}$, so that $d(h,h_n)\! <\!\epsilon$. It is routine to check that the relations ``$d(h_m,h_n)\!\leq\!\frac{p}{q+1}$" and
``$d(h_m,h_n)\! <\!\frac{p}{q+1}$" are recursive in $(m,n,p,q)$. \hfill{$\square$}\bigskip

 We saw in the proof of Proposition \ref{Wadge} that the evaluation map
$(h,\beta )\!\mapsto\! h(\beta )$ is continuous from ${\cal C}(2^\omega ,X)\!\times\! 2^\omega$ into $X$. We can say more if $X$ is recursively presented.

\begin{lem} \label{ev} Let $X$ be a recursively presented Polish space. Then the evaluation map is recursive.\end{lem}

\noindent\bf Proof.\rm\ Note that
$$\begin{array}{ll}
h(\beta )\!\in\! N(X,n)\!\!\!\!
& \Leftrightarrow d_X\big( h(\beta ),r^X_{((n)_1)_0}\big)\! <\!\frac{((n)_1)_1}{((n)_1)_2+1}\cr
& \Leftrightarrow\exists m,i,l\!\in\!\omega ~~\mbox{Seq}(m)~\wedge ~\mbox{lh}(m)\! =\! 2^i
~\wedge ~\beta\vert i\! =\! s^i_l~\wedge ~(m)_l\! =\!\big( (n)_1\big)_0~\wedge\cr
& \hfill{d(h,h_m)\! <\!\frac{((n)_1)_1}{((n)_1)_2+1}}\mbox{,}
\end{array}$$
which gives the result.\hfill{$\square$}\bigskip

\noindent - We then strengthen 7A.3 in [M] about {\bf primitive recursion} as follows. If
$Z,Y$ are recursively presented Polish spaces, $g\! :\! Z\!\rightarrow\! Y$ and
$h\! :\! Y\!\times\!\omega\!\times\! Z\!\rightarrow\! Y$ are $\Ca$-recursive and
$f\! :\!\omega\!\times\! Z\!\rightarrow\! Y$ is defined by
$$\left\{\!\!\!\!\!\!\!
\begin{array}{ll}
& f(0,z)\! :=\! g(z)\mbox{,}\cr
& f(n\! +\! 1,z)\! :=\! h\big( f(n,z),n,z\big)\mbox{,}
\end{array}
\right.$$
then $f$ is also $\Ca$-recursive. If $m\! :\! Z\!\rightarrow\! Z$ is $\Ca$-recursive, then the proof of 7A.3 in [M] shows that the map $f'\! :\!\omega\!\times\! Z\!\rightarrow\! Y$ defined by
$$\left\{\!\!\!\!\!\!\!
\begin{array}{ll}
& f'(0,z)\! :=\! g(z)\mbox{,}\cr
& f'(n\! +\! 1,z)\! :=\! h\Big( f'\big( n,m(z)\big) ,n,z\Big)\mbox{,}
\end{array}
\right.$$
is also $\Ca$-recursive. As in 7A.5 in [M], this can be extended to partial functions which are
$\Ca$-recursive on their domain.\bigskip

\noindent - We are ready for the proof of our main tool.\bigskip

\noindent\bf Proof of Theorem \ref{tool}.\rm\ 3E.6 in [M] provides
$\pi\! :\!\omega^\omega\!\rightarrow\! X$ recursive, ${\cal F}\!\in\!\Bormone (\omega^\omega )$ and a $\Borel$-recursive injection $\rho\! :\! X\!\rightarrow\!\omega^\omega$ such that
$\pi_{\vert {\cal F}}$ is injective, $\pi [{\cal F}]\! =\! X$ and $\rho$ is the inverse of
$\pi_{\vert {\cal F}}$. Let us show that the map $\mu\! :\! h\!\mapsto\!\pi\circ h$ is $\Borel$-recursive from ${\cal C}(2^\omega ,\omega^\omega )$ into ${\cal C}(2^\omega ,X)$. More generally, let $Y$ be a recursively presented Polish space, and $\psi\! :\! Y\!\rightarrow\! {\cal C}(2^\omega ,X)$. Note that
$$\begin{array}{ll}
\psi (y)\!\in\! N\big( {\cal C}(2^\omega ,X),n\big)\!\!\!\!
& \Leftrightarrow d\big(\psi (y),h_{((n)_1)_0}\big)\! <\!\frac{((n)_1)_1}{((n)_1)_2+1}\cr
& \Leftrightarrow\exists m\!\in\!\omega ~~\mbox{sup}_{\beta\in 2^\omega}~
d_X\big(\psi(y)(\beta ),h_{((n)_1)_0}(\beta )\big)\! <\!\frac{((m)_1)_1}{((m)_1)_2+1}\! <\!
\frac{((n)_1)_1}{((n)_1)_2+1}\cr
& \Leftrightarrow\exists m\!\in\!\omega ~~\forall\beta\!\in\! 2^\omega ~~
d_X\big(\psi(y)(\beta ),h_{((n)_1)_0}(\beta )\big)\! <\!\frac{((m)_1)_1}{((m)_1)_2+1}\! <\!
\frac{((n)_1)_1}{((n)_1)_2+1}
\end{array}$$
and $h_{((n)_1)_0}(\beta )\! =\! r^X_{g(n,\beta )}$ for some recursive map
$g\! :\!\omega\!\times\! 2^\omega\!\rightarrow\!\omega$.

\vfill\eject

 In the present case, $Y\! =\! {\cal C}(2^\omega ,\omega^\omega )$ and
$\psi (y)(\beta )\! =\!\pi\big( y(\beta )\big)$. Thus
$$\begin{array}{ll}
d_X\big(\psi(y)(\beta ),h_{((n)_1)_0}(\beta )\big)\! <\!\frac{((m)_1)_1}{((m)_1)_2+1}\!\!\!\!
& \Leftrightarrow
d_X\Big(\pi\big( y(\beta )\big),r^X_{g(n,\beta )}\Big)\! <\!\frac{((m)_1)_1}{((m)_1)_2+1}\cr
& \Leftrightarrow\pi\big( y(\beta )\big)\!\in\!
N\big( X,\big< 0,<\! g(n,\beta ),\big( (m)_1\big)_1,\big( (m)_1\big)_2\! >\!\big>\big)\cr
& \Leftrightarrow
\big( y(\beta ),\big< 0,<\! g(n,\beta ),\big( (m)_1\big)_1,\big( (m)_1\big)_2\! >\!\big>\big)\!\in\! G^\pi\mbox{,}
\end{array}$$
where $G^\pi$ is the $\Boraone$ neighborhood diagram of $\pi$. As the evaluation map is recursive, $h\!\mapsto\!\pi\circ h$ is $\Ca$-recursive and total, and thus $\Borel$-recursive.\bigskip

\noindent $\bullet$ Let us show that there is a $\Borel$-recursive map
$f\! :\!\omega^\omega\!\rightarrow\! {\cal C}(2^\omega ,X)$ such that
${\cal U}^{{\bf\Pi}^1_n,2^\omega}_\alpha\! =\!\big( f(\alpha )\big)^{-1}(C)$ for each
$\alpha\!\in\!\omega^\omega$. We adapt the proof of the main result in [K3]. We set
$A\! :=\!\pi^{-1}(C)$. As $C\!\in\! {\it\Pi}^1_n(X)$, $A\!\in\! {\it\Pi}^1_n(\omega^\omega )$. If
$<\beta^0,\delta^0>\in\! 2^\omega$, then we inductively define, for $i\!\in\!\omega$,
$m_i,\beta^{i+1},\delta^{i+1}$ as follows. If $(\beta^i,\delta^i)$ is given and in
$\mbox{Domain}(e)$, then $(\beta^i)^*(\delta^i)\! =<x_i,\beta^{i+1},\delta^{i+1}>$ and
$$m_i\! :=\!\left\{\!\!\!\!\!\!
\begin{array}{ll}
& \mbox{the location of the first 0 in }x_i\mbox{ if it exists,}\cr
& 2\mbox{ otherwise.}
\end{array}
\right.$$
We then set $Q\! :=\!\big\{ (\alpha ,<\beta^0,\delta^0>)\in\!\omega^\omega\!\times\! 2^\omega\mid
\forall i\!\in\!\omega ~~(\beta^i,\delta^i)\in\!\mbox{Domain}(e)~\wedge ~
\big(\alpha ,(m_i)\big)\!\in\! {\cal U}^{{\bf\Pi}^1_n,2^\omega}\big\}$ and $B^*\! :=\! Q_\alpha$, so that
$Q\!\in\! {\it\Pi}^1_n(\omega^\omega\!\times\! 2^\omega )$ and 
$\beta\!\in\! B^*\Leftrightarrow (\alpha ,\beta )\!\in\! Q$ for each
${(\alpha ,\beta )\!\in\!\omega^\omega\!\times\! 2^\omega}$ (note that $B^*$ depends on $\alpha$, but we denote it like this to keep the notation of [K3]). We define
$I\! :\!\omega^\omega\!\rightarrow\! 2^\omega$ by $I(\alpha )\! :=\! 0^{\alpha (0)}10^{\alpha (1)}1...$ Note that $I$ a $\Borel$-recursive injection onto the $\Bormtwo$ set
$$\mathbb{P}_\infty\! :=\!\{\beta\!\in\! 2^\omega\mid\forall p\!\in\!\omega ~~
\exists q\!\geq\! p~~\beta (q)\! =\! 1\}\mbox{,}$$ 
so that there is a $\Borel$-recursive map $\phi\! :\! 2^\omega\!\rightarrow\!\omega^\omega$ which is the inverse of $I$ on $\mathbb{P}_\infty$. We set
$$Q'\! :=\!\Big\{\delta\!\in\! 2^\omega\mid (\delta )_0\!\in\!\mathbb{P}_\infty ~\wedge~
\Big(\phi\big( (\delta )_0\big),(\delta )_1\Big)\!\in\! Q\Big\}\mbox{,}$$
so that $Q'\!\in\! {\it\Pi}^1_n(2^\omega )$. As ${\cal U}^{{\bf\Pi}^1_n,2^\omega}$ is suitable, there is $\alpha_Q\!\in\!\omega^\omega$ recursive with
$Q'\! =\! {\cal U}^{{\bf\Pi}^1_n,2^\omega}_{\alpha_Q}$. Note that
$$\begin{array}{ll}
\beta\!\in\! B^*\!\!\!\!
& \Leftrightarrow (\alpha ,\beta )\!\in\! Q\Leftrightarrow <I(\alpha ),\beta >\in\! Q'
\Leftrightarrow (\alpha_Q,<I(\alpha ),\beta >)\!\in\! {\cal U}^{{\bf\Pi}^1_n,2^\omega}\cr
& \Leftrightarrow b\big( {\cal R}(\alpha_Q,<I(\alpha ),\beta >)\big)\!\in\! C
\Leftrightarrow\rho\Big( b\big( {\cal R}(\alpha_Q,<I(\alpha ),\beta >)\big)\Big)\!\in\! A.
\end{array}$$
We set $G\! :=\! \rho\Big( b\big( {\cal R}(\alpha_Q,<I(\alpha ),.>)\big)\Big)$, so that
$G\! :\! 2^\omega\!\rightarrow\!\omega^\omega$ is $\Borel (\alpha )$-recursive and
$<\beta^0,\delta^0>$ is in $B^*$ if and only if $G(<\beta^0,\delta^0>)\!\in\! A$.\bigskip

\noindent $\bullet$ As in [K3], we can find
$F\! :\! 2^{<\omega}\!\rightarrow\! (2^\omega\!\times\!\omega )^{<\omega}$ satisfying the following properties:
$$\begin{array}{ll}
(1)\!\!\!\! & t\!\subseteq\! t'\Rightarrow F(t)\!\subseteq\! F(t')\cr
(2)\!\!\!\! & \vert F(t)\vert\! =\!\vert t\vert\! +\! 1\cr
(3)\!\!\!\! & (i)~\mbox{if }F(\emptyset )\! =\! (\beta^0,k_0)\mbox{, then }\beta^0\!\in\! {\cal D}~\wedge ~
\forall\delta^0\!\in\! [{\cal O}(\beta^0)]~~G(<\beta^0,\delta^0>)(0)\! =\! k_0\cr
& (ii)~\mbox{if }F(\varepsilon_0,...,\varepsilon_n)\! =\! 
(\beta^0,k_0,\beta^1,k_1,...,\beta^{n+1},k_{n+1})
\mbox{, then}\cr
& ~~~~~~(a)~\forall i\!\leq\! n\! +\! 1~~\beta^i\!\in\! {\cal D}\cr
& ~~~~~~(b)~\mbox{for all }\delta^{n+1}\!\in\! [{\cal O}(\beta^{n+1})]\mbox{, if }\delta^n,...,\delta^0\mbox{ are the uniquely determined members of}\cr
& ~~~~~~[{\cal O}(\beta^n)],...,[{\cal O}(\beta^0)]\mbox{ such that }
\forall i\!\leq\! n~~(\beta^i)^*(\delta^i)\! =<\overline{\varepsilon_i},\beta^{i+1},\delta^{i+1}>
\mbox{, where}\cr
& ~~~~~~\overline{\varepsilon_i}\! =\! 1^{\varepsilon_i}01^\infty\mbox{, then }
\forall i\!\leq\! n\! +\! 1~~G(<\beta^0,\delta^0>)(i)\! =\! k_i.
\end{array}$$

 We will need an effective version of this, so that we give the details of the construction of $F$. In fact, the $\beta^i$'s involved in the definition of $F$ can be $\Borel (\alpha )$. In order to see this, we first define
$$H_0\! :\! 2^\omega\!\times\! 2^\omega\!\rightarrow\!\omega$$
by $H_0(\beta ,\delta )\! :=\! G(<\beta ,\delta >)(0)$. As $G$ is $\Borel (\alpha )$-recursive, $H_0$ too, and the niceness of the coding system gives $\beta^0\!\in\! {\cal D}\cap \Borel (\alpha )$ and
$k_0\!\in\!\omega$ such that $G(<\beta^0,\delta^0>)(0)\! =\! k_0$ for each
$\delta^0\!\in\! [{\cal O}(\beta^0)]$. Now suppose that $n\!\in\!\omega$, 
$(\varepsilon_0,...,\varepsilon_n)$ and
$F(\varepsilon_0,...,\varepsilon_{n-1})\! =\! (\beta^0,k_0,...,\beta^n,k_n)$ are given. We define
$$H_{n+1}\! :\! 2^\omega\!\times\! 2^\omega\!\rightarrow\!\omega$$
as follows. Given $(\beta ,\delta )\!\in\! 2^\omega\!\times\! 2^\omega$, let $\delta^n,...,\delta^0$ be the uniquely determined members of $[{\cal O}(\beta^n)],...,$ $[{\cal O}(\beta^0)]$ resp., such that
$(\beta^n)^*(\delta^n)\! =<\overline{\varepsilon_n},\beta ,\delta >$, and 
$(\beta^i)^*(\delta^i)\! =<\overline{\varepsilon_i},\beta^{i+1},\delta^{i+1}>$ if $i\! <\! n$. Put
$H_{n+1}(\beta ,\delta )\! :=\! G(<\beta^0,\delta^0>)(n\! +\! 1)$. As $H_{n+1}$ is $\Borel (\alpha )$ (it is total and $\Ca (\alpha )$-recursive since $\iota$ is $\Ca$-recursive), the niceness of the coding system gives $\beta^{n+1}\!\in\! {\cal D}\cap \Borel (\alpha )$ and $k_{n+1}\!\in\!\omega$ such that ${G(<\beta^0,\delta^0>)(n\! +\! 1)\! =\! k_{n+1}}$ for each
$\delta^{n+1}\!\in\! [{\cal O}(\beta^{n+1})]$. Then 
$$F(\varepsilon_0,...,\varepsilon_n)\! :=\! (\beta^0,k_0,...,\beta^{n+1},k_{n+1})\mbox{,}$$ 
so that $F$ is as desired. So we can assume that the $\beta^i$'s are $\Borel (\alpha )$ in the conditions required for $F$.\bigskip

\noindent $\bullet$ By [K3] again, the map $h_\alpha\! :\! (\varepsilon_i)\!\mapsto\! (k_i)$ is continuous and ${\cal U}^{{\bf\Pi}^1_n,2^\omega}_\alpha\! =\! h_\alpha^{-1}(A)$. As this is not too long to prove, we give the details for completeness. The map $h_\alpha$ is in fact more than continuous: it is Lipschitz, by definition. Fix $(\varepsilon_i)$. We apply $F$ to the initial segments of $(\varepsilon_i)$, which gives $(\beta^i)$. For each $n$, we define perfect sets $C^n_0$, $C^n_1$, ..., $C^n_n\!\subseteq\! 2^\omega$ with $C^n_i\!\subseteq\! [{\cal O}(\beta^i)]$ if 
$i\!\leq\! n$, as follows:
$$\begin{array}{ll}
& C^n_n\! :=\!\{\delta^n\!\in\! [{\cal O}(\beta^n)]\mid\exists\delta^{n+1}\!\in\! 2^\omega ~~
(\beta^n)^*(\delta^n)\! =<\overline{\varepsilon_n},\beta^{n+1},\delta^{n+1}>\}\mbox{,}\cr
& C^n_{n-1}\! :=\!\{\delta^{n-1}\!\in\! [{\cal O}(\beta^{n-1})]\mid\exists\delta^n\!\in\! C^n_n~~
(\beta^{n-1})^*(\delta^{n-1})\! =<\overline{\varepsilon_{n-1}},\beta^n,\delta^n>\}\mbox{,}\cr
& ...\cr
& C^n_0\! :=\!\{\delta^0\!\in\! [{\cal O}(\beta^0)]\mid\exists\delta^1\!\in\! C^n_1~~
(\beta^0)^*(\delta^0)\! =<\overline{\varepsilon_0},\beta^1,\delta^1>\} .
\end{array}$$
Note that\bigskip

$(4)$ $\delta^0\!\in\! C^n_0\Rightarrow\ <\beta^i,\delta^i>\in\!\mbox{Domain}(e)$ for each 
$i\!\leq\! n$, where $\delta^1$, ..., $\delta^n$ are computed according to the formula in 
$(3).(ii).(b)$,\smallskip

$(5)$ $n'\!\geq\! n\Rightarrow\forall i\!\leq\! n~~C^{n'}_i\!\subseteq\! C^n_i$.\bigskip

\noindent This implies that $[{\cal O}(\beta^0)]\!\supseteq\! C^0_0\!\supseteq\! C^1_0\!\supseteq\! C^2_0\!\supseteq\! ...$ and $\bigcap_{n\in\omega}~C^n_0$ contains some $\delta^0$. Note that 
${<\beta^i,\delta^i>}$ is in Domain$(e)$, and 
$(\beta^i)^*(\delta^i)\! =<\overline{\varepsilon_i},\beta^{i+1},\delta^{i+1}>$ for each 
$i\!\in\!\omega$. By $(3).(ii).(b)$, 
$$G(<\beta^0,\delta^0>)\! =\! k_i$$ 
for each $i\!\in\!\omega$. As 
$<\beta^0,\delta^0>\in\! B^*\Leftrightarrow G(<\beta^0,\delta^0>)\!\in\! A$, 
$$\big(\forall i\!\in\!\omega ~<\beta^i,\delta^i>\in\!\mbox{Domain}(e)\wedge 
(\varepsilon_i)\!\in\! {\cal U}^{{\bf\Pi}^1_n,2^\omega}_\alpha\big)\Leftrightarrow (k_i)\!\in\! A.$$ 
As $<\beta^i,\delta^i>$ is in Domain$(e)$ for each $i\!\in\!\omega$, 
$(\varepsilon_i)\!\in\! {\cal U}^{{\bf\Pi}^1_n,2^\omega}_\alpha\Leftrightarrow 
h_\alpha\big( (\varepsilon_i)\big)\! =\! (k_i)\!\in\! A$.\bigskip

\noindent $\bullet$ So we found, for each $\alpha\!\in\!\omega^\omega$, 
$h_\alpha\!\in\! {\cal C}(2^\omega ,\omega^\omega )$ such that 
${\cal U}^{{\bf\Pi}^1_n,2^\omega}_\alpha\! =\! (\pi\circ h_\alpha )^{-1}(C)
\! =\!\big(\mu (h_\alpha )\big)^{-1}(C)$. It remains to see that the map
$\psi\! :\!\alpha\!\mapsto\! h_\alpha$, from $\omega^\omega$ into
${\cal C}(2^\omega ,\omega^\omega )$, can be $\Borel$-recursive (then $f$ will be 
$\mu\circ\psi$). By the previous discussion, it is enough to see that the relation ``$k_i\! =\! k$" is 
$\Borel$ in 
$\big(\alpha ,(\varepsilon_i),i,k\big)\!\in\!\omega^\omega\!\times\! 2^\omega\!\times\!\omega^2$. \bigskip

\noindent $\bullet$ We will define, by primitive recursion, a $\Borel$-recursive map 
$\tilde f\! :\!\omega\!\times\!\omega^\omega\!\times\! 2^\omega\!\rightarrow\! 
2^\omega\!\times\!\omega$ such that 
$\tilde f\big( n,\alpha ,(\varepsilon_i)\big)$ will be of the form 
$(<{\tilde\beta}^0,...,{\tilde\beta}^n,{\tilde\beta}^n,...>,<{\tilde k}^0,...,{\tilde k}^n>)$ and can play the role of $F(\varepsilon_0,...,\varepsilon_{n-1})$. We first set\bigskip

\leftline{$P\! :=\!\Big\{\big(\alpha ,(\varepsilon_i),\beta,k\big)\!\in\!
\omega^\omega\!\times\! (2^\omega )^2\!\times\!\omega\mid$}\smallskip

\rightline{$\forall i\!\in\!\omega ~~(\beta )_i\! =\! (\beta )_0\!\in\! {\cal D}\cap\Borel (\alpha )\wedge
\forall\delta\!\in\!\big[ {\cal O}\big( (\beta )_0\big)\big] ~~G(<(\beta )_0,\delta >)(0)\! =\! k
\Big\} .$}\bigskip
\noindent Note that $P$ is $\Ca$ and for any 
$\big(\alpha ,(\varepsilon_i)\big)\!\in\!\omega^\omega\!\times\! 2^\omega$ there is 
$(\beta ,k)\!\in\! 2^\omega\!\times\!\omega$ such that 
$\big(\alpha ,(\varepsilon_i),\beta ,k\big)\!\in\! P$. The uniformization lemma gives a $\Borel$-recursive map 
$\tilde g\! :\!\omega^\omega\!\times\! 2^\omega\!\rightarrow\! 2^\omega\!\times\!\omega$ such that 
$$\Big(\alpha ,(\varepsilon_i),\tilde g\big(\alpha ,(\varepsilon_i)\big)\Big)\!\in\! P$$ 
for each $\big(\alpha ,(\varepsilon_i)\big)\!\in\!\omega^\omega\!\times\! 2^\omega$. Then we set 
$$D\! :=\!\Big\{\big(\beta ,p,n,\alpha ,(\varepsilon_i)\big)\!\in\! 
2^\omega\!\times\!\omega^2\!\times\! \omega^\omega\!\times\! 2^\omega\mid
\mbox{Seq}(p)\wedge\mbox{lh}(p)\! =\! n\! +\! 1\wedge
\forall q\!\in\!\omega ~~(\beta )_q\!\in\! {\cal D}\cap\Borel (\alpha )\Big\} .$$ 
Note that $D$ is $\Ca$, as well as\bigskip
 
\rightline{$R\! :=\!\Big\{\big(\beta ,p,n,\alpha ,(\varepsilon_i),\beta',k'\big)\!\in\! 
D\!\times\! 2^\omega\!\times\!\omega\ \mid\ \forall i\! >\! n~~(\beta')_i\! =\! (\beta')_{n+1}\!\in\! {\cal D}\cap\Borel (\alpha )\ \wedge$}\smallskip

\rightline{$\mbox{Seq}(k')\wedge\mbox{lh}(k')\! =\! n\! +\! 2\wedge
\forall i\!\leq\! n~~(\beta')_i\! =\! (\beta )_i\wedge (k')_i\! =\! (p)_i\ \wedge$}\smallskip

\rightline{$\forall\delta\!\in\! 2^\omega ~~
\Big(\exists i\!\leq\! n\! +\! 1~~(\delta )_i\!\notin\!\big[ {\cal O}\big( (\beta')_i\big)\big]\ \vee\ 
\exists i\!\leq\! n~~(\beta')^*_i\big( (\delta )_i\big)\!\not=
<\overline{\varepsilon_i},(\beta')_{i+1},(\delta )_{i+1}>\vee$}\smallskip

\rightline{$\forall i\!\leq\! n\! +\! 1~~G\big( <(\beta')_0,(\delta )_0>\big)(i)\! =\! (k')_i\Big)\Big\} .$}\bigskip 

\noindent Moreover, for each $\big(\beta ,p,n,\alpha ,(\varepsilon_i)\big)\!\in\! 
D\! =\!\Pi_{2^\omega\times\omega^2\times \omega^\omega\times 2^\omega}[R]$ there is 
$(\beta',k')\!\in\!\big( 2^\omega\cap\Borel (\alpha )\big)\!\times\!\omega$ such that 
$\big(\beta ,p,n,\alpha ,(\varepsilon_i),\beta',k'\big)\!\in\! R$. The uniformization lemma gives a partial map $$\tilde h\! :\! 2^\omega\!\times\!\omega^2\!\times\!\omega^\omega\!\times\! 2^\omega
\!\rightarrow\! 2^\omega\!\times\!\omega$$ 
which is $\Ca$-recursive on its domain $D$, and such that 
$\Big(\beta ,p,n,\alpha ,(\varepsilon_i),\tilde h\big(\beta ,p,n,\alpha ,(\varepsilon_i)\big)\Big)\!\in\! R$ if 
$\big(\beta ,p,n,\alpha ,(\varepsilon_i)\big)\!\in\! D$. This implies that the partial map $\tilde f$ defined by
$$\left\{\!\!\!\!\!\!\!
\begin{array}{ll}
& \tilde f\big( 0,\alpha ,(\varepsilon_i)\big)\! :=\!\tilde g\big(\alpha ,(\varepsilon_i)\big)\mbox{,}\cr
& \tilde f\big( n\! +\! 1,\alpha ,(\varepsilon_i)\big)\! :=\!
\tilde h\Big(\tilde f\big( n,\alpha ,(\varepsilon_i)\big) ,n,\alpha ,(\varepsilon_i)\Big)\mbox{,}
\end{array}
\right.$$
is $\Ca$-recursive.

\vfill\eject

 Moreover, an induction shows that 
$\Big(\tilde f\big( n,\alpha ,(\varepsilon_i)\big) ,n,\alpha ,(\varepsilon_i)\Big)\!\in\! D$ for each 
$\big( n,\alpha ,(\varepsilon_i)\big)$, so that $\tilde f$ is in fact total, and thus $\Borel$-recursive. More precisely, $\tilde f\big( n,\alpha ,(\varepsilon_i)\big)$ is of the form 
$$(<\beta^0,...,\beta^n,\beta^n,...>,<k_0,...,k_n>)\mbox{,}$$ 
where $(\varepsilon_0,...,\varepsilon_{n-1})\!\mapsto\! (\beta^0,k_0,...,\beta^n,k_n)$ satisfies the properties (1)-(3) of $F$. It remains to note that 
$k_i\! =\!\tilde f\big( i,\alpha ,(\varepsilon_i)\big) (1)(i)$.\hfill{$\square$}\bigskip

\noindent - We now prove the consequences of our main tool.

\begin{defi} Let $\Gamma$ be a class of subsets of recursively presented Polish spaces, $\bf\Gamma$ be the corresponding boldface class, $X,Y$ be recursively presented Polish spaces, and
${\cal U}\!\in\!\Gamma (Y\!\times\! X)$. We say that ${\cal U}$ is {\bf effectively uniformly} 
$Y$-{\bf universal for the} $\bf\Gamma$ {\bf subsets of} $X$ if the following hold:\smallskip

(1) ${\bf\Gamma}(X)\! =\!\{ {\cal U}_y\mid y\!\in\! Y\}$,\smallskip

(2) $\Gamma(X)\! =\!\{ {\cal U}_y\mid y\!\in\! Y\ \Borel\mbox{-recursive}\}$,\smallskip

(3) for each $S\!\in\! {\bf\Gamma}(\omega^\omega\!\times\! X)$, there is a Borel map
$b\! :\!\omega^\omega\!\rightarrow\! Y$ such that $S_\alpha\! =\! {\cal U}_{b(\alpha )}$ for each
$\alpha\!\in\!\omega^\omega$,\smallskip

(4) for each $S\!\in\!\Gamma (\omega^\omega\!\times\! X)$, there is a $\Borel$-recursive map
$b\! :\!\omega^\omega\!\rightarrow\! Y$ such that $S_\alpha\! =\! {\cal U}_{b(\alpha )}$ for each
$\alpha\!\in\!\omega^\omega$.\end{defi}

\noindent\bf Notation.\rm\ Let ${\cal U}^{{\bf\Pi}^1_1,2^\omega}\!\in\!\Ca$ be a good
$\omega^\omega$-universal for the $\ca$ subsets of $2^\omega$, $X_1$ be a recursively presented Polish space, and ${\cal C}_1$ be a ${\it\Pi}^1_1$ subset of $X_1$ for which there is a $\Borel$-recursive map
$b\! :\!\omega^\omega\!\rightarrow\! X_1$ such that
$$(\alpha ,\beta )\!\in\! {\cal U}^{{\bf\Pi}^1_1,2^\omega}\Leftrightarrow
b(<\alpha ,\beta >)\!\in\! {\cal C}_1$$
if $(\alpha ,\beta )\!\in\!\omega^\omega\!\times\! 2^\omega$. We define, for each natural number $n\!\geq\! 1$,\bigskip

$\bullet$ $X_{n+1}\! :=\! {\cal C}(2^\omega ,X_n)$ (inductively),\bigskip

$\bullet$ ${\cal C}_{n+1}\! :=\!\{ h\!\in\! X_{n+1}\mid\forall\beta\!\in\! 2^\omega ~~
h(\beta )\!\notin\! {\cal C}_n\}$ (inductively),\bigskip

$\bullet$ ${\cal U}_n\! :=\!
\{ (h,\beta )\!\in\! X_{n+1}\!\times\! 2^\omega\mid h(\beta )\!\in\! {\cal C}_n\}$.

\begin{thm} \label{proj} Let $n\!\geq\! 1$ be a natural number. Then\smallskip

\noindent (a) the set ${\cal U}_n$ is effectively uniformly $X_{n+1}$-universal for the ${\bf\Pi}^1_n$ subsets of $2^\omega$,\smallskip

\noindent (b) the set ${\cal C}_n$ is ${\bf\Pi}^1_n$-complete.\end{thm}

\noindent\bf Proof.\rm\ We argue by induction on $n$.\bigskip

\noindent (a) Assume first that $n\! =\! 1$, and fix
$S\!\in\! {\bf\Pi}^1_1(\omega^\omega\!\times\! 2^\omega )$. Our assumption gives
$b_1\! :\!\omega^\omega\!\rightarrow\! X_1$. As ${\cal U}^{{\bf\Pi}^1_1,2^\omega}\!\in\!\Ca$ is a good $\omega^\omega$-universal for the $\ca$ subsets of $2^\omega$, there is by Theorem
\ref{tool} a $\Borel$-recursive map $f_1\! :\!\omega^\omega\!\rightarrow\! {\cal C}(2^\omega ,X_1)$ such that $(\alpha ,\beta )\!\in\! {\cal U}^{{\bf\Pi}^1_1,2^\omega}\Leftrightarrow
f_1(\alpha )(\beta )\!\in\! {\cal C}_1$ if
$(\alpha ,\beta )\!\in\!\omega^\omega\!\times\! 2^\omega$. Let $\alpha_S\!\in\!\omega^\omega$ with 
$S\! =\! {\cal U}^{\ca ,\omega^\omega\times 2^\omega}_{\alpha_S}$. Note that
$$\begin{array}{ll}
(\alpha ,\beta )\!\in\! S\!\!\!\!\!
& \Leftrightarrow\big( {\cal R}(\alpha_S,\alpha ),\beta\big)\!\in\! {\cal U}^{\ca ,2^\omega}
\Leftrightarrow f_1\big( {\cal R}(\alpha_S,\alpha )\big)(\beta )\!\in\! {\cal C}_1
\Leftrightarrow\Big( f_1\big( {\cal R}(\alpha_S,\alpha )\big),\beta\Big)\!\in\! {\cal U}_1.
\end{array}$$
As ${\cal C}_1$ is $\Ca$, ${\cal U}_1$ too. If $A\!\in\!\ca (2^\omega )$, then
$A\! =\! {\cal U}^{\ca ,2^\omega}_\alpha$ for some $\alpha\!\in\!\omega^\omega$. Applying the previous discussion to $S\! :=\! {\cal U}^{\ca ,2^\omega}$, we get
$A\! =\! ({\cal U}_1)_{f_1({\cal R}(\alpha_S,\alpha ))}$, so that
${\cal U}_1$ is $X_2$-universal for the $\ca$ subsets of $2^\omega$, effectively and uniformly.

\vfill\eject

 We now study ${\cal U}_{n+1}$. Fix
$S\!\in\! {\bf\Pi}^1_{n+1}(\omega^\omega\!\times\! 2^\omega )$. Let
${\cal U}^{{\bf\Pi}^1_n,2^\omega}$ be a good $\omega^\omega$-universal  for the ${\bf\Pi}^1_n$ subsets of $2^\omega$. We set ${\cal V}^{{\bf\Pi}^1_{n+1},2^\omega}\! :=\!
\big\{ (\alpha ,\beta )\!\in\!\omega^\omega\!\times\! 2^\omega\mid\forall\delta\!\in\! 2^\omega ~~
\big( {\cal R}(\alpha ,\beta ),\delta\big)\!\notin\! {\cal U}^{{\bf\Pi}^1_n,2^\omega}\big\}$, so that
${\cal V}^{{\bf\Pi}^1_{n+1},2^\omega}$ is a suitable $\omega^\omega$-universal for the
${\bf\Pi}^1_{n+1}$ subsets of $2^\omega$. Moreover, the induction assumption gives a $\Borel$-recursive map $b_{n+1}\! :\!\omega^\omega\!\rightarrow\! X_{n+1}$ such that
$$\begin{array}{ll}
(\alpha ,\beta )\!\in\! {\cal V}^{{\bf\Pi}^1_{n+1},2^\omega}\!\!\!\!\!
& \Leftrightarrow\forall\delta\!\in\! 2^\omega ~~
\big( {\cal R}(\alpha ,\beta ),\delta\big)\!\notin\! {\cal U}^{{\bf\Pi}^1_n,2^\omega}
\Leftrightarrow\forall\delta\!\in\! 2^\omega ~~
\Big( b_{n+1}\big( {\cal R}(\alpha ,\beta )\big),\delta\Big)\!\notin\! {\cal U}_n\cr
& \Leftrightarrow\forall\delta\!\in\! 2^\omega ~~
b_{n+1}\big( {\cal R}(\alpha ,\beta )\big)(\delta )\!\notin\! {\cal C}_n
\Leftrightarrow b_{n+1}\big( {\cal R}(\alpha ,\beta )\big)\!\in\! {\cal C}_{n+1}
\end{array}$$
Theorem \ref{tool} gives a $\Borel$-recursive map $f_{n+1}$ such that
$(\alpha ,\beta )\!\in\! {\cal V}^{{\bf\Pi}^1_{n+1},2^\omega}\Leftrightarrow
f_{n+1}(\alpha )(\beta )\!\in\! {\cal C}_{n+1}$ if
$(\alpha ,\beta )\!\in\!\omega^\omega\!\times\! 2^\omega$. Let
$$Q\!\in\! {\bf\Pi}^1_n(\omega^\omega\!\times\! 2^\omega\!\times\! 2^\omega)
\!\subseteq\! {\bf\Pi}^1_n(\omega^\omega\!\times\!\omega^\omega\!\times\! 2^\omega)$$ such that
$(\alpha ,\beta )\!\in\! S\Leftrightarrow\forall\delta\!\in\! 2^\omega ~~
(\alpha ,\beta ,\delta )\!\notin\! Q$, and $\alpha_Q\!\in\!\omega^\omega$ such that
$Q\! =\! {\cal U}^{{\bf\Pi}^1_n,\omega^\omega\times\omega^\omega\times 2^\omega}_{\alpha_Q}$. Note that
$$\begin{array}{ll}
(\alpha ,\beta )\!\in\! S\!\!\!\!\!
& \Leftrightarrow\forall\delta\!\in\! 2^\omega ~~
\Big( {\cal R}\big( {\cal R}'(\alpha_Q,\alpha ),\beta\big) ,\delta\Big)\!\notin\!
{\cal U}^{{\bf\Pi}^1_n,2^\omega}
\Leftrightarrow\big( {\cal R}'(\alpha_Q,\alpha ),\beta\big)\!\in\! {\cal V}^{{\bf\Pi}^1_{n+1},2^\omega}\cr
& \Leftrightarrow f_{n+1}\big( {\cal R}'(\alpha_Q,\alpha )\big)(\beta )\!\in\! {\cal C}_{n+1}\Leftrightarrow\Big( f_{n+1}\big({\cal R}'(\alpha_Q,\alpha )\big) ,\beta\Big)\!\in\! {\cal U}_{n+1}.
\end{array}$$
As ${\cal C}_n\!\in\! {\it\Pi}^1_n$, ${\cal C}_{n+1}\!\in\! {\it\Pi}^1_{n+1}$ and
${\cal U}_{n+1}\!\in\! {\it\Pi}^1_{n+1}$. If $A\!\in\! {\bf\Pi}^1_{n+1}(2^\omega )$, then
$A\! =\! {\cal U}^{{\bf\Pi}^1_{n+1},2^\omega}_\alpha$ for some $\alpha\!\in\!\omega^\omega$. Applying the previous discussion to $S\! :=\! {\cal U}^{{\bf\Pi}^1_{n+1},2^\omega}$, we get
$A\! =\! ({\cal U}_{n+1})_{f_{n+1}({\cal R}'(\alpha_Q,\alpha ))}$, so that ${\cal U}_{n+1}$ is
$X_{n+2}$-universal for the analytic subsets of $2^\omega$, effectively and uniformly.\bigskip

\noindent (b) By definition, ${\cal C}_1\!\in\! {\it\Pi}^1_1$, and
${\cal C}_{n+1}\!\in\! {\it\Pi}^1_{n+1}$ if ${\cal C}_n\!\in\! {\it\Pi}^1_n$. Assume first that
$E\!\in\! {\bf\Pi}^1_n(2^\omega )$. Then $E\! =\! ({\cal U}_n)_h$ for some
${h\!\in\! {\cal C}(2^\omega ,X_n)}$, by (a). Thus $E\! =\! h^{-1}({\cal C}_n)$. If $Z$ is a zero-dimensional Polish space and ${D\!\in\! {\bf\Pi}^1_n(Z)}$, then we may assume that $Z$ is a
$G_\delta$ subset of $2^\omega$ by 7.8 in [K2], so that $D\!\in\! {\bf\Pi}^1_n(2^\omega )$. The previous discussion gives $g\!\in\! {\cal C}(2^\omega ,X_n)$ with $D\! =\! g^{-1}({\cal C}_n)$. Thus
$D\! =\! (g_{\vert Z})^{-1}({\cal C}_n)$ and ${\cal C}_n$ is ${\bf\Pi}^1_n$-complete.\hfill{$\square$}\bigskip

\noindent\bf Proof of Theorem \ref{projapp}.\rm\ By Theorem \ref{proj}, it is enough to show that if
${\cal U}^{{\bf\Pi}^1_1,2^\omega}\!\in\!\Ca$ is a good $\omega^\omega$-universal set for the $\ca$ subsets of $2^\omega$, then there is a $\Borel$-recursive map 
$b\! :\!\omega^\omega\!\rightarrow\! [0,1]^{2^{<\omega}}$ such that 
$(\alpha ,\beta )\!\in\! {\cal U}^{{\bf\Pi}^1_1,2^\omega}\Leftrightarrow
b(<\alpha ,\beta >)\!\in\! {\cal P}$ 
if $(\alpha ,\beta )\!\in\!\omega^\omega\!\times\! 2^\omega$. Let
$H\!\in\!\Bormtwo (\omega^\omega\!\times\! 2^\omega\!\times\! 2^\omega )$ such that
$\neg {\cal U}^{\ca ,2^\omega}\! =\!\Pi_{\omega^\omega\times 2^\omega}[H]$. We set 
$G\! :=\!\big\{ (\alpha ,\beta )\!\in\!\omega^\omega\!\times\! 2^\omega\mid
\big( (\alpha )_0,(\alpha )_1,(\beta )_1\big)\!\in\! H~\wedge ~\beta\!\in\! {\cal K}\big\}$, 
so that $G\!\in\!\Borel (\omega^\omega\!\times\! 2^\omega )$, has $G_\delta$ vertical sections and
$G\!\subseteq\!\omega^\omega\!\times\! {\cal K}$. Lemma \ref{Borel} gives a $\Borel$-recursive map $F\! :\!\omega^\omega\!\rightarrow\! [0,1]^{2^{<\omega}}$, taking values in $\cal M$, and such that $G_\alpha\! =\! {\cal V}_{b(\alpha )}$ for each $\alpha\!\in\!\omega^\omega$. If 
$(\alpha ,\beta )\!\in\!\omega^\omega\!\times\! 2^\omega$, then
$$\begin{array}{ll}
(\alpha ,\beta )\!\notin\! {\cal U}^{\ca ,2^\omega}\!\!\!\!\!
& \Leftrightarrow\exists\delta\!\in\! 2^\omega ~~(\alpha ,\beta ,\delta )\!\in\! H
\Leftrightarrow\exists\delta\!\in\! 2^\omega ~~(<\alpha ,\beta >,\delta )\!\in\! G\cr
& \Leftrightarrow\exists\delta\!\in\! 2^\omega ~~\big( b(<\alpha ,\beta >),\delta\big)\!\in\! {\cal V}
\Leftrightarrow b(<\alpha ,\beta >)\!\notin\! {\cal P}.
\end{array}$$
This finishes the proof.\hfill{$\square$}\bigskip

\noindent\bf Questions.\rm\ Let $U$ be a $\Bormtwo$ subset of
$\omega^\omega\!\times\! 2^\omega$ which is universal for $\bormtwo (2^\omega )$. We set
$$G\! :=\!\big\{ (\alpha ,\beta )\!\in\!\omega^\omega\!\times\! {\cal K}\mid
\big(\alpha ,(\beta )_1\big)\!\in\! U\big\} .$$
Note that $G$ is a $\Bormtwo$ subset of $\omega^\omega\!\times\! 2^\omega$ contained in
$\omega^\omega\!\times\! {\cal K}$ which is universal for $\bormtwo ({\cal K})$. Indeed, fix
$H\!\in\!\bormtwo ({\cal K})$. Then
$H'\! :=\!\{\gamma\!\in\! 2^\omega\mid <0^\infty ,\gamma >\in\! H\}$ is $\bormtwo$, which gives
$\alpha_0\!\in\!\omega^\omega$ with $H'\! =\! U_{\alpha_0}$. Then $H\! =\! G_{\alpha_0}$.\bigskip

 Let $\alpha\!\mapsto\!\big( (\alpha )_k\big)_{k\in\omega}$ be a homeomorphism between
$\omega^\omega$ and $(\omega^\omega )^\omega$, with inverse map
$$(\alpha_k)_{k\in\omega}\!\mapsto <\alpha_0,\alpha_1,...>.$$
We set $S'\! :=\!\{\alpha\!\in\!\omega^\omega\mid\exists\gamma\!\in\!\omega^\omega ~~
\forall i\!\in\!\omega ~~\forall\beta\!\in\! 2^\omega ~~\beta\!\notin\! G_{(\alpha )_{\gamma (i)}}\}$. Note that $S'$ is ${\it\Sigma}^1_2$.\bigskip

\noindent (1) \it Is $S'$ a Borel ${\bf\Sigma}^1_2$-complete set?\rm\bigskip

 Assume that this is the case. Then the set
${\cal S}_2\! :=\!\{ (f_k)_{k\in\omega}\!\in\! {\cal M}^\omega\mid\exists\gamma\!\in\!\omega^\omega ~~\forall i\!\in\!\omega ~~f_{\gamma (i)}\!\in\! {\cal P}\}$ of sequences of martingales having a subsequence made of everywhere converging martingales is Borel ${\bf\Sigma}^1_2$-complete. Indeed, Lemma \ref{Borel} gives a Borel map $F\! :\!\omega^\omega\!\rightarrow\! {\cal M}$ such that $G_\alpha\! =\! {\cal V}_{F(\alpha)}$ for each $\alpha\!\in\!\omega^\omega$. The map
$\tilde F\! :\!\omega^\omega\!\rightarrow\! {\cal M}^\omega$ defined by
$\tilde F(\alpha )(k)\! :=\! F\big( (\alpha )_k\big)$ is Borel. Moreover,
$$\begin{array}{ll}
\tilde F(\alpha )\!\in\! {\cal S}_2\!\!\!\!
& \Leftrightarrow\exists\gamma\!\in\!\omega^\omega ~~
\forall i\!\in\!\omega ~~\forall\beta\!\in\! 2^\omega ~~
\beta\!\notin\! D\Big( F\big( (\alpha )_{\gamma (i)}\big)\Big)\cr
& \Leftrightarrow\exists\gamma\!\in\!\omega^\omega ~~
\forall i\!\in\!\omega ~~\forall\beta\!\in\! 2^\omega ~~
\beta\!\notin\! {\cal V}_{F((\alpha )_{\gamma (i)})}\cr
& \Leftrightarrow\exists\gamma\!\in\!\omega^\omega ~~
\forall i\!\in\!\omega ~~\forall\beta\!\in\! 2^\omega ~~
\beta\!\notin\! G_{(\alpha )_{\gamma (i)}}\cr
& \Leftrightarrow\alpha\!\in\! S'\mbox{,}
\end{array}$$
so that $S'\! =\! {\tilde F}^{-1}({\cal S}_2)$.\bigskip

\noindent (2) \it Is there a Borel map $f\! :\! {\cal C}(2^\omega ,[0,1])\!\rightarrow\!\omega^\omega$ such that, for each $(h_k)_{k\in\omega}\!\in\!\big( {\cal C}(2^\omega ,[0,1])\big)^\omega$ and each
$\beta\!\in\! 2^\omega$, the following are equivalent:\smallskip

(a) $\mbox{lim}_{k\rightarrow\infty}~h_k(\beta )\! =\! 0$,\smallskip

(b) $\forall k\!\in\!\omega ~~\beta\!\notin\! G_{f(h_k)}$?\rm\bigskip

 Assume that this is the case. Then $S'$ (and therefore ${\cal S}_2$) is Borel ${\bf\Sigma}^1_2$-complete, and thus ${\bf\Sigma}^1_2$-complete (see [P]).
We define $F\! :\!\big( {\cal C}(2^\omega ,[0,1])\big)^\omega\!\rightarrow\!\omega^\omega$ by
$F\big( (h_k)_{k\in\omega}\big)\! :=<f(h_0),f(h_1),...>$, so that $F$ is Borel. Note that
$$\begin{array}{ll}
F\big( (h_k)_{k\in\omega}\big)\!\in\! S'\!\!\!\!
& \Leftrightarrow\exists\gamma\!\in\!\omega^\omega ~~\forall i\!\in\!\omega ~~
\forall\beta\!\in\! 2^\omega ~~\beta\!\notin\! G_{f(h_{\gamma (i)})}\cr
& \Leftrightarrow\exists\gamma\!\in\!\omega^\omega ~~\forall\beta\!\in\! 2^\omega ~~
\mbox{lim}_{i\rightarrow\infty}~h_{\gamma (i)}(\beta )\! =\! 0\cr
& \Leftrightarrow (h_k)_{k\in\omega}\!\in\! S\mbox{,}
\end{array}$$
so that $S\! =\! F^{-1}(S')$.

\vfill\eject

\section{$\!\!\!\!\!\!$ References}

\noindent [B-Ka-L]\ \ H. Becker, S. Kahane and A. Louveau, Some complete ${\bf\Sigma}^1_2$ sets in harmonic analysis,\ \it Trans. Amer. Math. Soc.\ \rm 339, 1 (1993), 323-336

\noindent [C]\ \ M. Chaika, The Lusin-Menchoff theorem in metric space,~\it Indiana Univ. Math. J.~\rm 21 (1971/72), 351-354

\noindent [D]\ \ J. L. Doob,~\it Measure theory,~\rm volume 143 of Graduate Texts in Mathematics, Springer-Verlag, New York, 1994

\noindent [K1]\ \ A. S. Kechris, Measure and category in effective descriptive set theory,~\it Ann. Math. Logic~\rm 5 (1973), 337-384

\noindent [K2]\ \ A. S. Kechris,~\it Classical Descriptive Set Theory,~\rm Springer-Verlag, 1995

\noindent [K3]\ \ A. S. Kechris, On the concept of $\ca$-completeness,~\it Proc. Amer. Math. Soc.
\ \rm 125, 6 (1997), 1811-1814

\noindent [L]\ \ A. Louveau, A separation theorem for $\Ana$ sets,\ \it Trans. Amer. Math. Soc.\ \rm 260 (1980), 363-378

\noindent [Lu]\ \ J. Luke\v s,~The Lusin-Menchoff property of fine topologies,\ \it  Comment. Math. Univ. Carolinae,\ \rm 18, 3 (1977), 515-530

\noindent [Lu-Ma-Z]\ \ J. Luke\v s, J. Mal\'y, L. Zaji\v cek,~\it Fine Topology Methods in Real Analysis and Potential Theory,~\rm Lecture Notes in Mathematics 1189, Springer-Verlag, Berlin, 1986

\noindent [M]\ \ Y. N. Moschovakis,~\it Descriptive set theory,~\rm North-Holland, 1980

\noindent [P]\ \ J. Pawlikowski, On the concept of analytic hardness,~\it Proc. Amer. Math. Soc.\ \rm (2014)

\noindent [T1]\ \ H. Tanaka, Some results in the effective descriptive set theory,~\it Publ. Res. lnst. Math. Sci. Ser. A\ \rm 3 (1967), 11-52

\noindent [T2]\ \ H. Tanaka, A basis result for the $\ca$ sets of positive measure,~\it Comment. Math. Univ. St. Paul.\ \rm 16 (1968), 115-127

\noindent [Za]\ \ Z. Zahorski, Sur l'ensemble des points de non-d\'erivabilit\'e d'une fonction continue,~\it Bull. Soc. Math. France\ \rm 74 (1946), 147-178

\end{document}